\documentclass[11pt]{article}
\usepackage[margin=1in]{geometry}
\usepackage{amsmath,amssymb,amsthm,mathtools}
\usepackage[mathscr]{euscript}
\usepackage{enumitem}
\usepackage{graphicx}
\usepackage{booktabs,longtable}
\usepackage[section]{placeins}
\usepackage[pdfencoding=auto,psdextra,colorlinks=true,linkcolor=blue,citecolor=blue,urlcolor=blue]{hyperref}

\newtheorem{theorem}{Theorem}[section]
\newtheorem{klemma}[theorem]{Lemma}
\newtheorem{lemma}[theorem]{Lemma}
\newtheorem{proposition}[theorem]{Proposition}
\newtheorem{fact}[theorem]{Fact}
\newtheorem{corollary}[theorem]{Corollary}
\newtheorem{abstractthm}[theorem]{Theorem}
\theoremstyle{definition}
\newtheorem{definition}[theorem]{Definition}
\newtheorem{assumption}{Assumption}

\theoremstyle{remark}
\newtheorem{remark}[theorem]{Remark}
\numberwithin{equation}{section}

\newcommand{\HH}{L^2}
\newcommand{\R}{\mathbb R}
\newcommand{\Z}{\mathbb Z}
\newcommand{\N}{\mathbb N}
\newcommand{\C}{\mathbb C}
\newcommand{\T}{\mathbb T}
\newcommand{\BB}{\mathscr B(\HH)}
\newcommand{\KK}{\mathscr K(\HH)}
\newcommand{\Sch}{\mathscr S}
\newcommand{\jap}[1]{\langle #1\rangle}
\newcommand{\ip}[2]{\langle #1,#2\rangle}
\newcommand{\Mfr}{\mathfrak M}
\newcommand{\Ffr}{\mathfrak F_\rho}
\newcommand{\Acal}{\mathcal A_\rho}
\newcommand{\slim}{\mathop{\textrm{s-lim}}}
\newcommand{\Ran}{\operatorname{Ran}}
\newcommand{\supp}{\operatorname{supp}}
\newcommand{\dv}{\operatorname{div}}

\newcommand{\Op}{\mathrm{Op}}

\title{\bf Scattering theory and uniform-in-time convergence\\ of splitting methods for the Schr\"odinger equation}
\author{Guillaume Bal\thanks{Departments of Statistics and Mathematics and CCAM, University of
Chicago, Chicago, IL 60637; guillaumebal@uchicago.edu}}
\date{\today}

\providecommand{\one}{\mathbf 1}
\begin{document}
\maketitle
\phantomsection\pdfbookmark[0]{Paper}{Mainstart}

\begin{abstract}\noindent
We study the long-time behavior of splitting methods for the Schr\"odinger equation $i\partial_t\psi=(-\Delta_x+V)\psi$ on $\R^d$, with smooth potentials decaying like $\jap{x}^{-\rho}$, for Lie and Strang splitting and for compositions of any order $r$, with time step $\Delta$. We treat the one-step operator of the scheme as a unitary dynamics and develop its scattering theory.  Using a Mourre theory for unitary operators, we show for $\rho>1$ that its wave operators exist and are complete. Independently, for every $\rho>0$, we show that the error on scattering states localized in energy and space is $O(\Delta^r)$ uniformly in time. Consequently, the wave and scattering operators of the scheme converge to those of the Schr\"odinger equation, with the scattering operator converging at the effective order of the method, which is two for Lie splitting. Bound states behave differently, as they are embedded in the continuous spectrum of the scheme. We show that they survive for times longer than any power of $1/\Delta$. If $V$ has compact support or decays like $e^{-|x|^\beta}$ with $\beta>2$, they are in fact eventually radiated away by Lie and Strang splitting, for every $\Delta>0$.
\end{abstract}

\noindent\textbf{Keywords:} Splitting methods; Schr\"odinger equation; scattering theory; Mourre theory for unitary operators; uniform-in-time error estimates; long-range potentials; effective order; embedded eigenvalues; ionization; Fermi golden rule.

\medskip\noindent\textbf{MSC codes:} 65M12, 65M15, 35Q41, 47A40, 81U05.

%=====================================================================
\section{Introduction}\label{sec:intro}
%=====================================================================

This paper considers the long-time behavior of splitting methods for the Schr\"odinger equation $i\partial_t\psi=H\psi$, $H=H_0+V$, $H_0=-\Delta_x$, on $\R^d$. Splitting methods are standard time integrators for such equations, since each substep is explicit, the free flow being a Fourier multiplier and the potential flow a multiplication operator~\cite{MQ02,HLW06}. On a fixed interval $[0,T]$, the global error of a method of order $r$ is $O(T\Delta^r)$ under regularity assumptions \cite{JL00,Lub08,BGSvL26}, and such bounds grow linearly in $T$. Solutions in the absolutely continuous subspace of $H$, however, disperse. They leave every bounded region and, for short-range $V$, become asymptotically free. Since splitting errors are generated where the potential acts, one expects them not to accumulate on such states. The main objective of this paper is to show that this is indeed the case, with rates equal to the order of the method, uniformly in time and in the step $\Delta$. To describe the long-time behavior of the scheme itself, we also develop its scattering theory.

We first analyze the one-step operator $U$ of the method, e.g., the Lie step $U_{\rm L}=e^{-i\Delta H_0}e^{-i\Delta V}$. It is unitary, and $U_0=e^{-i\Delta H_0}$ samples the free flow exactly, $U_0^n=e^{-in\Delta H_0}$. We call the wave operators of the pair $(U,U_0)$ the {\em splitting wave operators}, and the associated scattering operator $S_U$ the {\em splitting scattering operator}. Two features distinguish this problem from scattering for a self-adjoint Hamiltonian. The first is aliasing: $U_0$ does not distinguish the energies $E$ and $E+2\pi m/\Delta$, so that $\sigma(U_0)=\T$ with infinite multiplicity and spectral projections of $U_0$ and $U$ do not localize in energy. The second is the lack of compactness: $U-U_0$ is a sum of terms such as $U_0(e^{-i\Delta V}-1)$, which are not compact, so that the perturbative results for unitary operators of \cite{FRT13} do not apply. We use the Mourre theory for unitary operators of \cite{ABCF06,FRT13}, with a conjugate operator whose commutator with the free step is positive away from zero kinetic energy. This positivity controls the high kinetic energies, and for the remaining terms it suffices that $(U-U_0)\eta(H_0)$ is compact for $\eta\in C_0([0,\infty))$. We obtain the absence of singular continuous spectrum for every $\rho>0$ and, for $\rho>1$, the existence and asymptotic completeness of the splitting wave operators at every fixed $\Delta$ (Theorems~\ref{thm:A} and~\ref{thm:B}). Only the free flow involves the Laplacian, and the theory extends to semibounded elliptic operators $\omega(p)$ with constant coefficients (Remark~\ref{rem:scalar}).

We next analyze the convergence. The global error is a sum of one-step defects along the exact solution. A defect of a method of order $r$ is $\Delta^{r+1}$ times an operator built from commutators of $H_0$ and $V$, and therefore decays in space one order faster than $V$. Along a scattering state of $H$, the local decay estimates of \cite[Theorem~5.1]{JMP84}, which hold for smooth potentials with any $\rho>0$, make these defects summable in time, and the error is $O(\Delta^r)$ uniformly in time (Theorem~\ref{thm:conv}). This argument does not use the scattering theory of $U$, and it applies to long-range potentials as well. Since wave and scattering operators are limits at large times, a reduction principle (Proposition~\ref{prop:reduction}) transfers the bound to them (Theorems~\ref{thm:wo} and~\ref{thm:S}), and it also shows that scattering states of $H$ have only an $O(\Delta^r)$ component on the point spectrum of $U$. The scattering operator can converge faster than the solution: if the method is conjugate to a method of higher order, as Lie splitting is to Strang splitting, its rate is this higher, effective order. For long-range potentials, $\rho\le1$, the free comparison dynamics must be modified, but the scheme can still be compared with the exact flow, through the relative wave operators $\Xi_\pm(\Delta)=\slim U^{-n}e^{-in\Delta H}P_{\rm ac}(H)$. These exist for every $\rho>0$ and are $O(\Delta^r)$-close to the identity on energy-localized data; composed with modified wave operators of $H$ \cite{Dol64,Hor76,DG97}, they give modified splitting wave operators, which are treated in the supplementary material (SM).

Uniform-in-time convergence fails on bound states of $H$, whose energies are embedded in the continuous spectrum $\sigma_{\rm ac}(U)=\T$ of the scheme. A bound state $H\phi_0=E_0\phi_0$ nevertheless persists for very long times: its survival probability stays close to $1$ up to times $\Delta^{-N}$ for every $N$, while its phase drifts with a shifted quasi-energy, $E_0+\frac{\Delta^2}{12}\ip{\phi_0}{|\nabla V|^2\phi_0}+O(\Delta^4)$ for Lie and Strang splitting (Theorem~\ref{thm:meta}). The error on a bound state therefore becomes of order one through its phase at times of order $\Delta^{-2}$, long before any mass is lost. When $\hat V$ decays exponentially, the survival time is even exponential in $\Delta^{-1/2}$ (Theorem~\ref{thm:expmeta}). Eventually, however, the scheme may {\em ionize} the bound state: its mass leaks into the continuous spectrum and disperses. The mechanism is again aliasing. The scheme couples the energy $E_0$ to the continuum energies $E_0+2\pi m/\Delta$, at momenta of order $(2\pi/\Delta)^{1/2}$, with a strength set by the Fourier transform of $V\phi_0$ at these momenta. We follow here the terminology of the literature on time-periodic forcing \cite{CLR00,CCLR01}; in the dispersive PDE literature, the same mechanism is called radiation damping \cite{SW99,KW01}. We prove that ionization occurs for every $\Delta>0$ when $V$ has compact support or decays like $e^{-a|x|^\beta}$ with $\beta>2$, since $U$ has then no eigenvalues at all, and for Gaussian potentials at large steps (Theorem~\ref{thm:noeig} and Corollary~\ref{cor:ion}). Numerically, the ionization rate follows a Fermi golden rule up to a factor close to one.

\medskip\noindent{\bf Related work.} Uniform-in-time error estimates for the Lie--Trotter splitting of defocusing NLS were obtained in \cite{CS24} by a vector-field method that provides time decay for both the exact and numerical solutions; see also \cite{Car26,CM26}. Discrete dispersive estimates for numerical schemes are studied in \cite{Ign11}, and improved long-time bounds for small potentials in \cite{BCF23}. On compact domains, where the spectrum is discrete, the long-time behavior of splitting methods is governed by resonances between eigenvalues of $H_0$, which occur only for exceptional time steps \cite{DF07,FGP10a,FJ15}. In the continuum, aliasing occurs for every time step, since every point of the unit circle is the image of infinitely many energies of $H_0$. We show that it does not prevent uniform-in-time convergence on scattering states, while it is responsible for the ionization of bound states. In quantum simulation, state-dependent Trotter bounds for unbounded Hamiltonians are obtained in \cite{BFHJY24}, and Trotterized dynamics is viewed as a Floquet system in \cite{HHZ19}. Scattering theory for pairs of unitary operators is presented in \cite{Yaf92}, and commutator methods for unitary operators are developed in \cite{ABCF06,FRT13,ABC15a} and applied to quantum walks in \cite{RST19,Tie20}. The latter framework would also apply to the fully discrete scheme obtained after a spatial discretization, a unitary operator on $\ell^2$ whose wave operators we do not consider here. Formally, $U$ is the monodromy of the time-periodic Hamiltonian $H_0+\Delta V\sum_k\delta(t-k\Delta)$, so that the problem is related to scattering for time-periodic Hamiltonians \cite{Yaj77,How79,Yok98} and to ionization by time-periodic fields \cite{CLR00,CCLR01,SW98,KW01}. Time dependence through Dirac masses is, however, outside the scope of those works, and we are not aware of an analysis of the wave operators of splitting schemes.

\medskip\noindent{\bf Outline.} Section~\ref{sec:results} states the main results. Sections~\ref{sec:scattering} and~\ref{sec:convergence} develop the scattering and convergence theories, and section~\ref{sec:bound} treats bound states. The appendix collects the abstract tools, the construction of the conjugate operator and the commutator computations. SM contains the proof of Theorem~\ref{thm:expmeta} and of the extension of Theorem~\ref{thm:noeig} used in Corollary~\ref{cor:ion}, together with an overview of the methods (section~\ref{sec:Sbackground}), numerical illustrations, further results on bound states, the long-range theory, the extension to $\omega(p)$, technical facts and an index of notation (section~\ref{sec:Snotation}).

\medskip\noindent{\bf Notation.} We work in the Hilbert space $L^2=L^2(\R^d)$, $d\ge1$, with scalar product antilinear in the first argument; $x$ is the position, $p=-i\nabla$ the momentum, $\jap{\xi}=(1+|\xi|^2)^{1/2}$, $\Sch=\Sch(\R^d)$, $\mathcal H^k$ the Sobolev spaces. For a function $b$, $b(p)$ is the Fourier multiplier, with $\hat u(\xi)=\int e^{-ix\cdot\xi}u(x)\,dx$. $\BB$, $\KK$ denote the bounded and compact operators, $\T$ the unit circle, $C_0(\R^d)$ the continuous functions vanishing at infinity, and $C_b^\infty$ the smooth functions with bounded derivatives. For a unitary or self-adjoint operator $T$, $E^T$ is its spectral measure and $P_{\rm ac}(T)$, $P_{\rm pp}(T)$ the projections onto the absolutely continuous and pure point subspaces $L^2_{\rm ac}(T)$ and $L^2_{\rm pp}(T)$. For an operator $T$, $\mathrm{Dom}(T)$ is its domain. We write $\Delta_x$ for the Laplacian and reserve $\Delta>0$ for the time step. Constants $C,c$ may change from line to line. An index of notation, which also lists the symbols used with several meanings, is given in section~\ref{sec:Snotation} of SM.

%=====================================================================
\section{Setting and main results}\label{sec:results}
%=====================================================================

This section presents the main results of the paper, namely the scattering theory and convergence of splitting schemes, and the behavior of bound states.

\subsection{Setting}\label{sec:setting}

\begin{assumption}\label{ass:V}
$V\in C^\infty(\R^d;\R)$ and, for some $\rho>0$ and every multi-index $\alpha$, $|\partial^\alpha V(x)|\le C_\alpha\jap{x}^{-\rho-|\alpha|}$.
\end{assumption}

Let $H_0=-\Delta_x$ and $H=H_0+V$, self-adjoint on $\mathcal H^2$. Under \eqref{ass:V}, $H$ has a Mourre estimate on every compact subinterval of $(0,\infty)$ with respect to the dilation generator, and $H$ is of class $C^\infty$ with respect to it \cite{Mou81,ABG96}, \cite[\S5]{JMP84}. Consequently $\sigma_{\rm ess}(H)=[0,\infty)$, $H$ has no singular continuous spectrum and no positive eigenvalues \cite{FH82}, and
\[
\HH_{\rm ac}(H)=E^H((0,\infty))\HH .
\]
For $\rho>1$ the wave operators $\Omega^c_\pm=\slim_{t\to\pm\infty}e^{itH}e^{-itH_0}$ exist and are complete, $\Ran\Omega^c_\pm=\HH_{\rm ac}(H)$ \cite{RS79}. We write $W^c_\pm=(\Omega^c_\pm)^*$, so that $W^c_\pm=\slim e^{itH_0}e^{-itH}P_{\rm ac}(H)$, and $S=(\Omega^c_+)^*\Omega^c_-$ for the scattering operator. Each result uses only finitely many derivatives of $V$; see section~\ref{sec:Stools} of SM.

\medskip\noindent{\bf Splitting methods.} For $\Delta>0$ we consider compositions
\begin{equation}\label{eq:comp}
U=F_sK_s\cdots F_1K_1,\quad F_j=e^{-ib_j\Delta H_0},\quad K_j=e^{-ia_j\Delta V},\quad a_j,b_j\in\R,\quad\sum_ja_j=\sum_jb_j=1 .
\end{equation}
Examples are the Lie step $U_{\rm L}=e^{-i\Delta H_0}e^{-i\Delta V}$, the reversed Lie step $U_{\rm R}=e^{-i\Delta V}e^{-i\Delta H_0}$, and the Strang step $U_{\rm S}=e^{-i\Delta V/2}e^{-i\Delta H_0}e^{-i\Delta V/2}$, as well as the triple-jump and optimized compositions of \cite{Yos90,MQ02}. The coefficients may be negative; for real methods of order $\ge3$ this is unavoidable, as shown in \cite{BC05}. All such $U$ are unitary, and
\[
U_0:=e^{-i\Delta H_0}=F_s\cdots F_1 .
\]
We write $U(s)$ for the composition with $\Delta$ replaced by $s$. The kick $K_j$ acts at the time $\beta_j\Delta$ of the step, $\beta_j=b_1+\dots+b_{j-1}$, and we call $\beta_1\Delta,\dots,\beta_s\Delta$ the kick times of the method. For $a\in\R$ we write $K(a)=e^{-ia\Delta V}$ for the kick of strength $a$, so that $K_j=K(a_j)$ in \eqref{eq:comp} and
\begin{equation}\label{eq:UminusU0}
U_{\rm L}-U_0=U_0\,(K(1)-1).
\end{equation}

\begin{definition}[Order, processing, effective order]\label{def:order}
\textup{(i)} The composition \eqref{eq:comp} has \emph{order $r$} if, as a formal power series in $s$ in the free associative algebra generated by $H_0$ and $V$, $U(s)=e^{-isH}+O(s^{r+1})$. \textup{(ii)} A \emph{processor} is a finite product $P$ of kicks $e^{-ic\Delta V}$ and free flows $e^{-ic'\Delta H_0}$; $C_P\Delta$ denotes the sum of the durations of its free flows. \textup{(iii)} $U$ has \emph{effective order $q$} if $U=P\tilde UP^{-1}$ for a processor $P$ and a composition $\tilde U$ of order $q$.
\end{definition}

Lie splitting has order $1$ and effective order $2$, since $U_{\rm L}=K(\tfrac12)^*U_{\rm S}K(\tfrac12)$, while Strang splitting has order $2$. Processing and effective order go back to \cite{But69}; see \cite{BCR99} for splitting methods.

\emph{Aliasing notation.} We write $\varepsilon_\Delta(E)=e^{-i\Delta E}$ and, for $0\le a<b\le2\pi/\Delta$, $\Theta_{[a,b]}=\varepsilon_\Delta([a,b])\subset\T$. Every compact subset of $\T\setminus\{1\}$ lies in some $\Theta_{[a,b]}$ with $0<a<b<2\pi/\Delta$. The point $1=\varepsilon_\Delta(0)$ is the image of the threshold $E=0$ of $H_0$ and of its aliases $2\pi k/\Delta$, $k\in\Z$. Every neighborhood of $1$ in $\T$ therefore contains the image of energies close to $0$, where the commutator estimate of Lemma~\ref{lem:free}(c) degenerates, and $1$ plays the role of a threshold of $U_0$. This is why $1$ is removed from $\T$ in the spectral statements below; see Figure~\ref{fig:aliasing} of SM.

\medskip\noindent{\bf Splitting wave operators.} For a composition $U$ we define, whenever the limits exist, the splitting wave operators $\Omega^U_\pm$, their inverses $W^U_\pm$, and the splitting scattering operator $S_U$ by
\[
\Omega^U_\pm(\Delta)=\slim_{n\to\pm\infty}U^{-n}U_0^n,\qquad W^U_\pm(\Delta)=\slim_{n\to\pm\infty}U_0^{-n}U^nP_{\rm ac}(U),\qquad S_U=\Omega^U_+(\Delta)^*\,\Omega^U_-(\Delta).
\]
For the Lie step we write simply $\Omega_\pm(\Delta)$ and $S_\Delta$.

\medskip\noindent{\bf Data classes.} We use two classes of states, and all quantitative statements are operator-norm bounds on them.
\emph{Energy-localized initial data:} $g(H)\jap{x}^{-s}$, with $g\in C_c^\infty((0,\infty))$ and $s>1$. For every fixed $s>1$, their ranges are dense in $\HH_{\rm ac}(H)$ when $g$ varies. \emph{Asymptotic data:} $J_\chi=\chi(p)\jap{x}^{-\rho}$, with $\chi\in C_c^\infty(\R^d\setminus\{0\})$, recalling that $p=-i\nabla$. Every $f$ with $\hat f\in C_c^\infty(\R^d\setminus\{0\})$ is of the form $J_\chi(\jap{x}^\rho f)$ with $\chi=1$ on $\supp\hat f$, so that the ranges of the $J_\chi$ span a dense subspace.
All constants below depend on $V$, on the method, and on $g,s$ or $\chi$, but not on $\Delta\in(0,1]$ or on the data. They degenerate as $\supp g$ or $\supp\chi$ approaches $0$, the threshold of $H_0$ and $H$, whose image under $\varepsilon_\Delta$ is the threshold $1\in\T$ of $U_0$, and they grow with the largest energy in $\supp g$ or $\supp\chi$.

\subsection{Scattering theory and convergence}\label{sec:mainresults}

\begin{theorem}[Spectral theory of splitting schemes]\label{thm:A}
Assume \eqref{ass:V}, and let $U$ be a composition \eqref{eq:comp}. For every $\Delta>0$:
\\[1mm]
\textup{(a)} every eigenvalue of $U$ in $\T\setminus\{1\}$ has finite multiplicity, and $\sigma_p(U)\setminus\{1\}$ has no accumulation point in $\T\setminus\{1\}$;
\\[1mm]
\textup{(b)} $U$ has no singular continuous spectrum;
\\[1mm]
\textup{(c)} for $s\in(1/2,1]$ and every compact $\Theta'\subset\T\setminus(\{1\}\cup\sigma_p(U))$, there is $C$ such that $\sum_{n\in\Z}\|\jap{x}^{-s}U^nE^U(\Theta')f\|^2\le C\|f\|^2$ for all $f\in\HH$; in the terminology of Definition~\ref{def:smooth}, $\jap{x}^{-s}$ is locally $U$-smooth on $\T\setminus(\{1\}\cup\sigma_p(U))$.
\end{theorem}
Part (c) is the propagation estimate from which the completeness of the splitting wave operators follows (Theorem~\ref{thm:B}(b)). Part (a) allows eigenvalues, but in many cases the point spectrum is empty (section~\ref{sec:bound-results}). The proof is given in section~\ref{sec:proofA}.

\begin{theorem}[Scattering theory of splitting schemes]\label{thm:B}
Assume \eqref{ass:V} with $\rho>1$, and let $U$ be a composition.
\\[1mm]
\textup{(a)} \emph{(Existence, uniformly in $\Delta$.)} For every $\Delta>0$, $\Omega^U_\pm(\Delta)$ exist, are isometric, satisfy $U\Omega^U_\pm=\Omega^U_\pm U_0$, and $\Ran\Omega^U_\pm\subset\HH_{\rm ac}(U)$. For every $\chi$ there is $C$ such that $\|(\Omega^U_\pm(\Delta)-U^{\mp N}U_0^{\pm N})J_\chi\|\le C\jap{N\Delta}^{1-\rho}$ for all $\Delta\in(0,1]$ and $N\ge0$.
\\[1mm]
\textup{(b)} \emph{(Asymptotic completeness.)} For every $\Delta>0$, $W^U_\pm(\Delta)$ exist, $\Ran\Omega^U_\pm(\Delta)=\HH_{\rm ac}(U)$ and $W^U_\pm=(\Omega^U_\pm)^*$. Thus $\HH=\HH_{\rm ac}(U)\oplus\HH_{\rm pp}(U)$, $\sigma_{\rm ac}(U)=\T$, and $S_U$ is unitary and commutes with $U_0$.
\\[1mm]
\textup{(c)} \emph{(Processing.)} If $U'=PUP^{-1}$ for a processor $P$, then $\Omega^{U'}_\pm=P\,\Omega^U_\pm\,\Pi$ and $S_{U'}=\Pi^*S_U\Pi$ with $\Pi=e^{iC_P\Delta H_0}$. In particular $U_{\rm L}$, $U_{\rm R}$, $U_{\rm S}$ are unitarily equivalent and have the same splitting scattering operator.
\end{theorem}
The bound in (a) is in physical time $N\Delta$, uniformly in $\Delta$: the limits defining $\Omega^U_\pm$ are reached after a time that does not depend on the step. By (b), every state of $\HH_{\rm ac}(U)$ is asymptotically free under the scheme. Parts (a) and (c) are proved in section~\ref{sec:discreteWO}, and part (b) in section~\ref{sec:proofA}.

\begin{theorem}[Uniform-in-time convergence]\label{thm:conv}
Assume \eqref{ass:V} and $U$ a composition of order $r$.
\\[1mm]
\textup{(a)} For $g\in C_c^\infty((0,\infty))$ and $s>1$ there is $C$ such that for all $\Delta\in(0,1]$
\[
\sup_{n\in\Z}\big\|(U^n-e^{-in\Delta H})\,g(H)\jap{x}^{-s}\big\|\le C\Delta^r .
\]
\par\vspace{1mm}\noindent \textup{(b)} For every $\psi\in\HH_{\rm ac}(H)$, $\lim_{\Delta\to0}\sup_{n\in\Z}\|U^n\psi-e^{-in\Delta H}\psi\|=0$.
\\[1mm]
\textup{(c)} If $\rho>2$, then for every $\chi$ there is $C$ such that for all $\Delta\in(0,1]$
\[
\sup_{\pm n\ge0}\big\|(U^n-e^{-in\Delta H})\,\Omega^c_\pm J_\chi\big\|\le C\Delta^r .
\]
\par\vspace{1mm}\noindent \textup{(d)} \emph{(Relative wave operators.)} The limits $\Xi_\pm(\Delta)=\slim_{n\to\pm\infty}U^{-n}e^{-in\Delta H}P_{\rm ac}(H)$ exist, are isometric on $\HH_{\rm ac}(H)$, satisfy $U\Xi_\pm=\Xi_\pm e^{-i\Delta H}$ and $\Ran\Xi_\pm\subset\HH_{\rm ac}(U)$, and, with $g,s$ as in (a), $\|(\Xi_\pm-1)g(H)\jap{x}^{-s}\|\le C\Delta^r$. For Lie splitting, $\Xi^{\rm L}_\pm=K(\frac12)^*\Xi^{\rm S}_\pm$ and $\|(\Xi^{\rm L}_\pm-1-\frac{i\Delta}2V)g(H)\jap{x}^{-s}\|\le C\Delta^2$.
\end{theorem}

On energy-localized data the error is thus $O(\Delta^r)$ uniformly in $n\in\Z$, whereas the standard bound $Cn\Delta^{r+1}$ grows linearly in time. The estimate also holds for data described by their free asymptotics rather than by their decay at time $0$: the states $\Omega^c_\pm J_\chi u$ behave as $e^{-itH_0}J_\chi u$ as $t\to\pm\infty$, and for them the error is uniform on the half-line $\pm n\ge0$, at the price of the stronger decay $\rho>2$. The error does not vanish as $n\to\pm\infty$ but converges, to $\|(\Xi_\pm-1)\psi\|$. For Lie splitting this limit is $\frac\Delta2\|V\psi\|+O(\Delta^2)$, and applying $K(\frac12)^*=e^{i\Delta V/2}$ once to the initial datum reduces it to the Strang limit, since $U_{\rm L}^nK(\frac12)^*=K(\frac12)^*U_{\rm S}^n$.

\begin{theorem}[Spectral stability and convergence of splitting wave operators]\label{thm:wo}
Assume \eqref{ass:V} and $U$ a composition of order $r$. Let $g\in C_c^\infty((0,\infty))$ and $s>1$. There is $C$ such that for all $\Delta\in(0,1]$:
\\[1mm]
\textup{(a)} \emph{(Spectral stability.)} $\|P_{\rm pp}(U)\,g(H)\jap{x}^{-s}\|\le C\Delta^r$. 
\\[1mm]
\textup{(b)} If $\rho>1$: $\Omega^U_\pm=\Xi_\pm\Omega^c_\pm$, $\Ran\Xi_\pm=\HH_{\rm ac}(U)$, and $\|(W^U_\pm(\Delta)-W^c_\pm)g(H)\jap{x}^{-s}\|+\|(\Omega^U_\pm(\Delta)-\Omega^c_\pm)W^c_\pm g(H)\jap{x}^{-s}\|\le C\Delta^r$. Consequently $\Omega^U_\pm(\Delta)\to\Omega^c_\pm$ and $W^U_\pm(\Delta)P_{\rm ac}(H)\to W^c_\pm$ strongly as $\Delta\to0$. 
\\[1mm]
\textup{(c)} If $\rho>2$: $\|(\Omega^U_\pm(\Delta)-\Omega^c_\pm)J_\chi\|\le C\Delta^r$ for all $\chi$.
\end{theorem}
By (a), a scattering state of $H$ has only an $O(\Delta^r)$ component on the eigenvectors of $U$; by (b) and (c), the asymptotic states of the scheme converge at the order of the method.

\begin{theorem}[Splitting scattering operator]\label{thm:S}
Assume \eqref{ass:V}, and let $U$ have effective order $q$, $U=P\tilde UP^{-1}$ with $\tilde U$ of order $q$ (Definition~\ref{def:order}). Let $G=g(H)\jap{x}^{-s}$ with $g\in C_c^\infty((0,\infty))$ and $s>1$. There is $C$ such that for all $\Delta\in(0,1]$:
\\[1mm]
\textup{(a)} The relative scattering operator $\Sigma_U=\Xi_+^*\Xi_-$ satisfies $\|G^*(\Sigma_U-P_{\rm ac}(H))G\|\le C\Delta^q$.
\\[1mm]
\textup{(b)} If $\rho>1$, then $S_U=(\Omega^c_+)^*\Sigma_U\Omega^c_-$ and, for $u,u'\in\HH$, $f'=W^c_+Gu'$ and $f=W^c_-Gu$,
\[
|\ip{f'}{(S_U-S)f}|\le C\Delta^q\|u\|\,\|u'\| .
\]
If $\rho>2$, also $\|J_\chi^*(S_U-S)J_\chi\|\le C\Delta^q$.
\\[1mm]
In particular, the relative scattering operator of Lie splitting is of second order for every $\rho>0$, and its splitting scattering operator for $\rho>1$. Moreover, if $\rho>1$, then for Lie splitting,
\[
\big\|\big(\Omega_\pm(\Delta)-(1+\tfrac{i\Delta}2V)\Omega^c_\pm\big)W^c_\pm g(H)\jap{x}^{-s}\big\|\le C\Delta^2.
\]
\end{theorem}
The rate $q$ can exceed the order $r$, as for Lie splitting, because on scattering states the processor $P$ acts as a free flow (Theorem~\ref{thm:B}(c)), which commutes with $S$. The proofs of Theorems~\ref{thm:conv}--\ref{thm:S} are given in section~\ref{sec:summation}.

\medskip\noindent{\bf Aliased channels.} $S$ commutes with $H_0$, but $S_U$ commutes only with $U_0$ (Theorem~\ref{thm:B}(b)), so it can map the energy $E$ to $E+2\pi m/\Delta$, $m\ne0$. This part of $S_U-S$ is included in the bound of Theorem~\ref{thm:S}; numerically it can dominate the error of high-order methods at moderate steps (section~\ref{sec:Snumerics} of SM).

\medskip\noindent{\bf Short- and long-range potentials.} The results on the spectrum of $U$, on $\Xi_\pm$ and $\Sigma_U$, and on the error for energy-localized data hold for every $\rho>0$. Those involving $\Omega^U_\pm$ or $\Omega^c_\pm$ require $\rho>1$, and $\rho>2$ for asymptotic data $J_\chi$. For $\rho\le1$, modified splitting wave operators are obtained from $\Xi_\pm$ (section~\ref{sec:Slongrange} of SM).

\begin{remark}[Other dispersive operators]\label{rem:scalar}
The results extend to $H_0=\omega(p)$, with $\omega$ a real polynomial of degree $2m$ whose principal part is positive away from $0$, and with $\T\setminus\{1\}$ replaced by $\T\setminus\varepsilon_\Delta(\mathcal T(\omega))$, where $\mathcal T(\omega)$ is the finite set of critical values of $\omega$. Theorems~\ref{thm:A} and~\ref{thm:B} extend, and the convergence results extend with the same rates under a local decay hypothesis for $\omega(p)+V$, which for $-\Delta_x$ is Proposition~\ref{prop:localdecay}. The changes in the proofs are given in section~\ref{sec:Sscalar} of SM.
\end{remark}

\subsection{Bound states}\label{sec:bound-results}

Uniform-in-time convergence fails on bound states of $H$ (Remark~\ref{rem:boundstates}). Let $H\phi_0=E_0\phi_0$ with $E_0<0$ and $\|\phi_0\|=1$. We say that $\phi_0$ is {\em ionized} by the scheme if $\phi_0\in\HH_{\rm ac}(U)$. Since $U$ has no singular continuous spectrum (Theorem~\ref{thm:A}(b)), this holds if and only if the survival amplitude $\ip{\phi_0}{U^n\phi_0}$ tends to $0$ in the mean, i.e., $\frac1N\sum_{n<N}|\ip{\phi_0}{U^n\phi_0}|^2\to0$, by Wiener's theorem (Proposition~\ref{prop:boundproved}). The numerical solution is then fully radiated, in the sense that $U^n\phi_0\rightharpoonup0$ weakly, so that no mass remains on $\phi_0$ or on any other bound state of $H$, and, for $\rho>1$, $\|U^n\phi_0-U_0^nW^U_+\phi_0\|\to0$, so that $U^n\phi_0$ leaves every bounded region. A bound state of $H$ nevertheless survives under the scheme for very long times, with an explicit drift of its phase.

\begin{theorem}[Metastability and quasi-energy]\label{thm:meta}
Assume \eqref{ass:V}, and let $E_0<0$ be a simple eigenvalue of $H$ with normalized eigenfunction $\phi_0$. Let $U$ be a composition of order $r$. There are real numbers $e_r,e_{r+1},\dots$, depending only on $V$, $\phi_0$ and the method, such that for every $N\ge r$, with $E_N(\Delta)=E_0+\sum_{k=r}^Ne_k\Delta^k$, there is $C_N$ with
\begin{equation}\label{eq:meta}
\big|\ip{\phi_0}{U^n\phi_0}-e^{-in\Delta E_N(\Delta)}\big|\le C_N\big(\Delta^{2r}+|n|\Delta^{N+2}\big)\qquad(n\in\Z,\ \Delta\in(0,1]).
\end{equation}
For the Strang step, $e_k=0$ for odd $k$ and
\begin{equation}\label{eq:quasiE}
e_2=\frac1{12}\ip{\phi_0}{|\nabla V|^2\phi_0}.
\end{equation}
If $U=P\tilde UP^{-1}$ for a processor $P$, then $U$ satisfies \eqref{eq:meta} with the numbers $e_k$ of $\tilde U$, provided $\Delta^{2r}$ is replaced by $\Delta^2$. In particular, the Lie and reversed Lie steps satisfy \eqref{eq:meta} with the numbers $e_k$ of Strang splitting.
\end{theorem}

Thus the survival probability $|\ip{\phi_0}{U^n\phi_0}|^2$ stays $1-O(\Delta^{2r})$ for $n\Delta\le\Delta^{2r-N-1}$, for every $N$, while the phase drifts away from $e^{-in\Delta E_0}$ on the time scale $\Delta^{-2}$ for Lie and Strang splitting. The rate of any decay of the survival amplitude is therefore $O(\Delta^N)$ for every $N$.

For potentials whose Fourier transform decays exponentially, which is the case for potentials analytic in a strip, the survival time is at least exponentially long in $\Delta^{-1/2}$. We assume
\begin{equation}\label{ass:A}
M_\alpha:=(2\pi)^{-d}\int_{\R^d}|\hat V(\zeta)|\,e^{\alpha|\zeta|}\,d\zeta<\infty\qquad\text{for some }\alpha>0.
\end{equation}

\begin{theorem}[Exponential metastability]\label{thm:expmeta}
Assume \eqref{ass:V} and \eqref{ass:A}, let $E_0<0$ be a simple eigenvalue of $H$ with normalized eigenfunction $\phi_0$, and let $U$ be the Lie, reversed Lie or Strang step. There are $c,C,\Delta_0>0$ and, for $\Delta\in(0,\Delta_0]$, real numbers $\tilde E(\Delta)$ with $\tilde E(\Delta)-E_N(\Delta)=O(\Delta^{N+1})$ for every $N$, where $E_N$ is as in Theorem~\ref{thm:meta}, such that
\begin{equation}\label{eq:expmeta}
\big|\ip{\phi_0}{U^n\phi_0}-e^{-in\Delta\tilde E(\Delta)}\big|\le C\big(\Delta^2+|n|e^{-c/\sqrt\Delta}\big)\qquad(n\in\Z,\ \Delta\in(0,\Delta_0]).
\end{equation}
\end{theorem}
The proof, based on a compression of the scheme to momenta $|\xi|\le(\delta/\Delta)^{1/2}$, is given in section~\ref{sec:Sexpmeta} of SM. The bound state is eventually ionized when $V$ decays sufficiently fast, since the point spectrum of $U$ can then be excluded altogether. Indeed, the eigenvalue equation forces the Fourier transform of $(e^{-i\Delta V}-1)w$ to vanish on the infinitely many resonant spheres $|\xi|^2=E+2\pi m/\Delta$, whose radii grow only like $m^{1/2}$, and an entire function of controlled growth cannot vanish on all of them.

\begin{theorem}[Absence of eigenvalues]\label{thm:noeig}
Let $V\in L^\infty(\R^d;\R)$, $\Delta>0$, and let $U$ be the Lie, reversed Lie or Strang step. Assume one of the following:
\\[1mm]
\textup{(a)} $V$ has compact support;
\\[1mm]
\textup{(b)} $|V(x)|\le Ce^{-a|x|^\beta}$ for some $a>0$ and $\beta>2$;
\\[1mm]
\textup{(c)} $|V(x)|\le Ce^{-a|x|^2}$ for some $a>0$, and $\Delta>\pi/(4a)$.
Then $U$ has no eigenvalues, including at $1$.
\end{theorem}

\begin{corollary}[Ionization of all bound states]\label{cor:ion}
Assume \eqref{ass:V} and one of the conditions of Theorem~\ref{thm:noeig}, and let $U$ be the Lie, reversed Lie or Strang step, or, under condition (a) or (b), a composition of any order whose kick times are pairwise distinct modulo $\Delta$, possibly after processing. Then $U$ is purely absolutely continuous, the splitting wave operators $\Omega^U_\pm(\Delta)$ are unitary, and for every $f\in\HH$, $\ip{f}{U^nf}\to0$ and $\|U^nf-U_0^nW^U_+f\|\to0$ as $n\to\infty$. In particular, every bound state $\phi_0$ of $H$ is ionized, and $\|\one_BU^n\phi_0\|\to0$ for every bounded set $B$.
\end{corollary}

The proofs of Theorems~\ref{thm:meta} and~\ref{thm:noeig} and of Corollary~\ref{cor:ion} are given in section~\ref{sec:bound}.

%=====================================================================
\section{Scattering theory of splitting schemes}\label{sec:scattering}
%=====================================================================

We now prove Theorems~\ref{thm:A} and~\ref{thm:B}: existence of the splitting wave operators by a Cook argument (section~\ref{sec:discreteWO}), a Mourre estimate for $U$ (section~\ref{sec:kick}), and completeness from local $U$-smoothness (section~\ref{sec:proofA}). Throughout this section $U$ is a composition \eqref{eq:comp}. We write $U=T_m\cdots T_1$ with factors $T_i\in\{F_j,K_j\}$.

\subsection{Existence of the splitting wave operators}\label{sec:discreteWO}

For the composition \eqref{eq:comp}, removing the kicks one at a time from the right gives
\begin{equation}\label{eq:telescope-comp}
U-U_0=\sum_{j=1}^sX_j\,\nu_j\,Y_j,\qquad X_j=F_sK_s\cdots F_{j+1}K_{j+1}F_j,\quad \nu_j=K_j-1,\quad Y_j=F_{j-1}\cdots F_1 ,
\end{equation}
where $Y_j=e^{-i\beta_j\Delta H_0}$ is a free flow and $|\nu_j|\le|a_j|\Delta|V|$. Indeed, if $U_{[j]}$ denotes $U$ with the kicks $K_1,\dots,K_{j-1}$ replaced by $1$, then $U_{[j]}-U_{[j+1]}=X_j\nu_jY_j$, $U_{[1]}=U$ and $U_{[s+1]}=U_0$. For Lie this is \eqref{eq:UminusU0}. The existence proof uses the following decay estimate for free waves.

\begin{lemma}[Free decay, operator form]\label{prop:decay}
Let $\chi\in C_c^\infty(\R^d\setminus\{0\})$ and $\mu\ge0$. There is $C$ such that for all $t\in\R$
\[
\big\|\jap{x}^{-\mu}e^{-itH_0}\chi(p)\jap{x}^{-\mu}\big\|\le C\jap{t}^{-\mu}.
\]
Moreover, $\|\jap{x}^{-\mu}\partial^\beta e^{-itH_0}J_\chi\|\le C\jap{t}^{-\min(\mu,\rho)}$ for any multi-index $\beta$.
\end{lemma}

\begin{proof}
For $|t|\le1$ the norm is at most $\|\chi\|_\infty$. Let $|t|\ge1$, $c_\chi=\min\{|\xi|:\xi\in\supp\chi\}>0$, $c=c_\chi/2$, $P_t=\mathbf 1_{\{|x|\le c|t|\}}$ and $Q_t=1-P_t$. Then ${\mathcal X} = Q_t{\mathcal X}  + P_t{\mathcal X} Q_t + P_t{\mathcal X} P_t$ for ${\mathcal X}=\jap{x}^{-\mu}e^{-itH_0}\chi(p)\jap{x}^{-\mu}$.
The first two terms have norm at most $\jap{ct}^{-\mu}\|\chi\|_\infty$. The last has the integral kernel $P_t(x)\jap{x}^{-\mu}K_t(x-y)\jap{y}^{-\mu}P_t(y)$ with
\[
K_t(z)=(2\pi)^{-d}\int e^{i\Phi(\xi)}\chi(\xi)\,d\xi,\qquad\Phi(\xi)=z\cdot\xi-t|\xi|^2 .
\]
On its support $|z|\le2c|t|=c_\chi|t|$, so for $\xi\in\supp\chi$, $|\nabla\Phi|=|z-2t\xi|\ge2c_\chi|t|-c_\chi|t|=c_\chi|t|$. With $Y=\nabla\Phi/|\nabla\Phi|^2$ and $L^{\mathsf T}w=i\nabla\cdot(wY)$ we have $K_t(z)=(2\pi)^{-d}\int e^{i\Phi}(L^{\mathsf T})^N\chi\,d\xi$. Since $\nabla\Phi$ is affine in $\xi$ with derivative $-2tI$, and $Y$ is homogeneous of degree $-1$ in $\nabla\Phi$, $|\partial^\beta_\xi Y|\le C_\beta|t|^{|\beta|}|\nabla\Phi|^{-1-|\beta|}\le C_\beta c_\chi^{-1-|\beta|}|t|^{-1}$. Hence $|K_t(z)|\le C_N|t|^{-N}$, with $C_N$ depending on finitely many derivatives of $\chi$, and the Hilbert--Schmidt norm of the last term is at most $C_N|t|^{-N}|B(0,c|t|)|\le C_N'|t|^{d-N}$. Take $N\ge d+\mu$. The last statement follows from $\partial^\beta e^{-itH_0}\chi(p)=e^{-itH_0}(i\xi)^\beta\chi(p)$, with $(i\xi)^\beta\chi\in C_c^\infty(\R^d\setminus\{0\})$, and from $\jap{x}^{-\mu},\jap{x}^{-\rho}\le\jap{x}^{-\min(\mu,\rho)}$, by the first statement with $\mu$ replaced by $\min(\mu,\rho)$.
\end{proof}

 The following proposition proves Theorem~\ref{thm:B}(a).

\begin{proposition}[Existence, uniformly in $\Delta$]\label{prop:discreteWO}
Assume \eqref{ass:V} with $\rho>1$ and let $U$ be a composition \eqref{eq:comp}. For every $\Delta>0$ the limits $\Omega^U_\pm(\Delta)$ exist on $\HH$, are isometric, satisfy $U\Omega^U_\pm(\Delta)=\Omega^U_\pm(\Delta)U_0$, and $\Ran\Omega^U_\pm(\Delta)\subset\HH_{\rm ac}(U)$. Moreover, for every $\chi$ there is $C$ such that for all $\Delta\in(0,1]$ and $N\ge0$,
\[
\big\|(\Omega^U_\pm(\Delta)-U^{\mp N}U_0^{\pm N})J_\chi\big\|\le C\jap{N\Delta}^{1-\rho}.
\]
\end{proposition}

\begin{proof}
By \eqref{eq:telescope-comp},
$U^{-(n+1)}U_0^{n+1}-U^{-n}U_0^n=-U^{-(n+1)}(U-U_0)U_0^n=-U^{-(n+1)}\sum_jX_j\nu_jY_jU_0^n$, so
\[
\big\|(U^{-(n+1)}U_0^{n+1}-U^{-n}U_0^n)J_\chi\big\|\le\sum_j|a_j|\Delta\big\|V e^{-i(n+\beta_j)\Delta H_0}J_\chi\big\|\le C\Delta\jap{n\Delta}^{-\rho},
\]
since $|\nu_j|\le|a_j|\Delta|V|$ pointwise, $|V|\le C\jap{x}^{-\rho}$, Lemma~\ref{prop:decay} with $\mu=\rho$, and $\jap{(n+\beta_j)\Delta}\ge c\jap{n\Delta}$ for $\Delta\le1$. For $\Delta>1$ the same argument applies with constants depending on $\Delta$. Summing over $n\ge N$ and using $\Delta\sum_{n\ge N}\jap{n\Delta}^{-\rho}\le\Delta\jap{N\Delta}^{-\rho}+\int_{N\Delta}^\infty\jap{t}^{-\rho}dt$ gives the tail bound and convergence on the ranges of the $J_\chi$; density and $\|U^{-n}U_0^n\|=1$ give existence on $\HH$. For $n\to-\infty$ use $U^{-(n-1)}U_0^{n-1}-U^{-n}U_0^n=U^{-n}(U-U_0)U_0^{n-1}$. Isometry and intertwining are immediate. The intertwining gives $\Omega^*E^U(B)\Omega=E^{U_0}(B)$ for Borel $B\subset\T$, and $U_0$ is purely absolutely continuous (its spectral measures are images of those of $H_0$ under $\varepsilon_\Delta$, which is a local diffeomorphism), so $\Ran\Omega^U_\pm\subset\HH_{\rm ac}(U)$.
\end{proof}

The next proposition proves Theorem~\ref{thm:B}(c).
\begin{proposition}[Conjugate schemes]\label{prop:strang}
Assume \eqref{ass:V} with $\rho>1$. Let $U$ be a composition and $U'=PUP^{-1}$, where $P$ is a finite product of kicks $e^{-ic\Delta V}$ and free flows whose durations add up to $C_P\Delta$. Then
\[
\Omega^{U'}_\pm(\Delta)=P\,\Omega^U_\pm(\Delta)\,\Pi,\qquad S_{U'}=\Pi^*S_U\Pi,\qquad\Pi=e^{iC_P\Delta H_0}.
\]
In particular, since $U_{\rm S}=K(\tfrac12)U_{\rm L}K(\tfrac12)^*$ and $U_{\rm R}=K(1)U_{\rm L}K(1)^*$ (so $\Pi=1$), the Lie, reversed Lie and Strang steps are unitarily equivalent and have the same splitting scattering operator.
\end{proposition}

\begin{proof}
$U'^{-n}U_0^n=PU^{-n}P^{-1}U_0^n$. For $f$ with $\hat f\in C_c^\infty(\R^d\setminus\{0\})$, $P^{-1}U_0^nf-\Pi U_0^nf\to0$: replacing the inverse kicks in $P^{-1}$ by $1$ one at a time produces terms containing a factor $e^{ic\Delta V}-1$ acting on a free state $e^{-i(n+c')\Delta H_0}f$, which tends to $0$ by Lemma~\ref{prop:decay}. Since $\Pi$ commutes with $U_0$, $U'^{-n}U_0^nf=PU^{-n}U_0^n\Pi f+o(1)\to P\Omega^U_\pm\Pi f$. Density gives the first claim, and $S_{U'}=\Pi^*\Omega_+^{U*}P^*P\Omega^U_-\Pi=\Pi^*S_U\Pi$.
\end{proof}

\subsection{Commutators and the Mourre estimate}\label{sec:kick}
The Mourre estimate for $U$ is obtained from that of $U_0$, using that $U-U_0$ becomes compact after a low-energy cut-off (Lemma~\ref{lem:lowenergy}).

We define the conjugate operator
\[
A=\tfrac12\big(x\cdot v(p)+v(p)\cdot x\big),\qquad v(\xi)=\frac{\xi}{\jap{\xi}^2},
\]
and write $\gamma(\xi)=\xi\cdot v(\xi)=|\xi|^2\jap{\xi}^{-2}$ and $R=\jap{p}^{-2}$, so that $\gamma(p)=1-R$. We also use the symbol class $\Mfr$ of $b\in C^\infty(\R^d)$ with $|\partial^\beta b(\xi)|\le C_\beta\jap{\xi}^{-|\beta|}$, the class $\Ffr$ of $a\in C^\infty(\R^d)$ with $|\partial^\alpha a(x)|\le C_\alpha\jap{x}^{-\rho-|\alpha|}$, and the algebra $\Acal$ spanned by finite products of operators $b(p)$, $b\in\Mfr$, and $a(x)$, $a\in\Ffr$; their properties are collected in appendix~\ref{sec:conjugate}. The regularity classes $C^k(A)$, the Mourre estimate for unitary operators and local smoothness are recalled in appendix~\ref{sec:abstract}.

\begin{klemma}[Free flows]\label{lem:free}
\textup{(a)} For $\tau\in\R$, $F_\tau=e^{-i\tau\Delta H_0}\in C^k(A)$ for every $k$, and $F_\tau^*[A,F_\tau]=2\tau\Delta\,\gamma(p)$.
\\[1mm]
\textup{(b)} $U_0$ is purely absolutely continuous and $\sigma(U_0)=\T$.
\\[1mm]
\textup{(c)} If $0<\theta_1<\theta_2<2\pi/\Delta$ and $\Theta\subset\Theta_{[\theta_1,\theta_2]}$ is Borel, then $E^{U_0}(\Theta)U_0^*[A,U_0]E^{U_0}(\Theta)\ge\frac{2\Delta\theta_1}{1+\theta_1}E^{U_0}(\Theta)$.
\\[1mm]
\textup{(d)} If $B\in\BB$ extends to an element of $\mathscr B(\mathrm{Dom}(\jap{A}^s)^*,\HH)$ for some $s\in(1/2,1]$, then $B$ is locally $U_0$-smooth on $\T\setminus\{1\}$. This applies to $B=\jap{x}^{-s}$.
\end{klemma}

\begin{proof}
(a) $F_\tau=u_\tau(p)$ with $u_\tau(\xi)=e^{-i\tau\Delta|\xi|^2}$. Since $v\cdot\nabla u_\tau=-2i\tau\Delta\,\gamma\,u_\tau$ and $v\cdot\nabla$ maps $\Mfr$ into $\Mfr$, induction gives $(v\cdot\nabla)^ku_\tau=c_ku_\tau$ with $c_k\in\Mfr$. Proposition~\ref{prop:mult} gives $F_\tau\in C^k(A)$ and $[A,F_\tau]=i(v\cdot\nabla u_\tau)(p)=2\tau\Delta \gamma(p)F_\tau$.
(b) The spectral measures of $U_0$ are images of those of $H_0$, which are absolutely continuous on $[0,\infty)$, under $\varepsilon_\Delta$. On each interval $[2\pi k/\Delta,2\pi(k+1)/\Delta)$ the map $\varepsilon_\Delta$ is a bijection onto $\T$ with derivative of modulus $\Delta$, so null sets have null preimages.
(c) $E^{U_0}(\Theta)=\mathbf 1_{\varepsilon_\Delta^{-1}(\Theta)}(H_0)$ with $\varepsilon_\Delta^{-1}(\Theta)\cap[0,\infty)\subset[\theta_1,\infty)$, and $\gamma(p)=\gamma_0(H_0)$ with $\gamma_0(E)=E/(1+E)$ increasing.
(d) A compact $\Theta'\subset\T\setminus\{1\}$ lies in an open arc $\Theta$ with $1\notin\overline\Theta$, on which (a), (c) give a strict Mourre estimate; $\sigma_p(U_0)=\emptyset$, and Theorem~\ref{thm:FRT}(iii) applies. For $B=\jap{x}^{-s}$ use Proposition~\ref{prop:weights}.
\end{proof}

\begin{klemma}[Kicks]\label{lem:kick}
Assume \eqref{ass:V}.
\\[1mm]
\textup{(a)} For $a\in\R$ and $K(a)=e^{-ia\Delta V}$, $\nu_a=K(a)-1\in\Ffr$, $K(a)=1+\nu_a\in\Acal$, and hence $K(a)\in C^k(A)$ for every $k$; moreover $|\nu_a|\le|a|\Delta|V|$.
\\[1mm]
\textup{(b)} $[A,K(a)]\in\KK$ and $[\gamma(p),K(a)]\in\KK$.
\end{klemma}

\begin{proof}
(a) We have $|\nu_a|\le|a|\Delta|V|\le C|a|\Delta\jap{x}^{-\rho}$. For $|\alpha|\ge1$, $\partial^\alpha \nu_a=e^{-ia\Delta V}\,Q_\alpha$, where $Q_\alpha$ is a finite sum of products $\prod_i(-ia\Delta\,\partial^{\alpha_i}V)$ with $\sum\alpha_i=\alpha$ and each $|\alpha_i|\ge1$. Each factor is $O(\jap{x}^{-\rho-|\alpha_i|})\le O(\jap{x}^{-|\alpha_i|})$ and at least one factor carries the extra $\jap{x}^{-\rho}$, so $\partial^\alpha \nu_a=O(\jap{x}^{-\rho-|\alpha|})$. The rest follows from Proposition~\ref{prop:algebra}.

(b) Since $A1-1A=0$ on $\Sch$, $[A,K(a)]=[A,\nu_a]$, given on $\Sch$ by \eqref{eq:comm-a} applied to the function $\nu_a$ (with $h_j=x_j\nu_a$). Since $\rho>0$, $\nu_a$ and $\partial_lh_j=\delta_{jl}\nu_a+x_j\partial_l\nu_a$ are in $C_0(\R^d)$, as are $\dv v$, $\jap{\xi}^{-2}$ and $\xi_l\jap{\xi}^{-2}$. Each term of \eqref{eq:comm-a} is then compact by Fact~\ref{fact:compact}: $(\dv v)(p)\nu_a$, $\nu_a(\dv v)(p)$ and $(\partial_jh_j)R$ directly; $(p_jp_lR)\cdot[(\partial_lh_j)R]$ as bounded times compact; and $(p_jR)\cdot[(\partial_lh_j)(p_lR)]$ likewise. Finally $\gamma=1-\jap{\xi}^{-2}$ gives $[\gamma(p),K(a)]=-[R,\nu_a]=\nu_aR-R\nu_a\in\KK$.
\end{proof}

Only $|\nu_a|+\jap{x}|\nabla \nu_a|\to0$ at infinity is used for the compactness in (b); this is the discrete counterpart of the condition $\jap{x}|\nabla V|\to0$ familiar from the Mourre theory of $H_0+V$.

\begin{klemma}[Commutator of a composition]\label{lem:product}
Assume \eqref{ass:V}. Then $U\in C^2(A)\subset C^{1+0}(A)$, and with $Q_i=T_{i-1}\cdots T_1$,
\begin{equation}\label{eq:UAU2}
U^*[A,U]=\sum_iQ_i^*\big(T_i^*[A,T_i]\big)Q_i=2\Delta(1-R)+\mathsf K_1,\qquad \mathsf K_1\in\KK .
\end{equation}
\end{klemma}

\begin{proof}
Each factor is in $C^2(A)$ by Lemmas~\ref{lem:free} and \ref{lem:kick}, hence so is $U$ (Fact~\ref{fact:reg}(b),(d)). The first identity is the Leibniz rule. For $T_i=F_j$, $T_i^*[A,T_i]=2b_j\Delta \gamma(p)$ by Lemma~\ref{lem:free}(a) with $\tau=b_j$; for $T_i=K_j$, $T_i^*[A,T_i]\in\KK$. Moreover $Q_i^*\gamma(p)Q_i-\gamma(p)\in\KK$: free flows commute with $\gamma(p)$, and $K(a)^*\gamma(p)K(a)-\gamma(p)=K(a)^*[\gamma(p),K(a)]\in\KK$ for kicks; induct on the number of factors. Hence $U^*[A,U]=2\Delta\big(\sum_jb_j\big)\gamma(p)+\KK=2\Delta(1-R)+\KK$. Negative coefficients play no role, only $\sum_jb_j=1$ enters.
\end{proof}

The following lemma replaces the hypothesis ``$V-1\in\KK$'' of \cite[Cor.~2.10, Lemma~3.1]{FRT13}, which fails here.

\begin{klemma}[Compactness after a kinetic-energy cut-off]\label{lem:lowenergy}
Assume \eqref{ass:V}. For every $\varphi\in C(\T)$ and $\eta\in C_0([0,\infty))$, $(\varphi(U)-\varphi(U_0))\eta(H_0)\in\KK$ and $\eta(H_0)(\varphi(U)-\varphi(U_0))\in\KK$.
\end{klemma}

\begin{proof}
By \eqref{eq:telescope-comp}, $U-U_0=\sum_iX_i\nu_iY_i$ with $Y_i$ free flows, which commute with $\eta(H_0)$, and $\nu_i\eta(H_0)\in\KK$ by Fact~\ref{fact:compact}. Hence $(U-U_0)\eta(H_0)\in\KK$. For $m\ge1$, $U^m-U_0^m=\sum_{j<m}U^j(U-U_0)U_0^{m-1-j}$ gives $(U^m-U_0^m)\eta(H_0)\in\KK$. Moreover $(U^{-m}-U_0^{-m})\eta(H_0)=-U^{-m}\big[(U^m-U_0^m)\eta(H_0)\big]U_0^{-m}\in\KK$. Fej\'er approximation extends this to $\varphi\in C(\T)$, since $\|\varphi(U)-P(U)\|\le\|\varphi-P\|_\infty$. The second statement is the adjoint of the first.
\end{proof}

Without the cut-off $\eta(H_0)$, the statement fails since $U-U_0$ is not compact. The cut-off suffices because the commutator \eqref{eq:UAU2} degenerates only at $p=0$. In the Mourre estimate, the high kinetic energies are controlled by its positivity, and compactness is needed only for the remaining term, which contains a factor $\eta(H_0)$.

\begin{klemma}[Mourre estimate]\label{lem:mourre}
Assume \eqref{ass:V}. Let $0<\theta_1<\theta_2<2\pi/\Delta$ and let $\Theta\subset\Theta_{[\theta_1,\theta_2]}$ be Borel. Then for every $c_0<2\Delta\theta_1/(1+\theta_1)$ there is $\mathsf K\in\KK$ with
\begin{equation}\label{eq:mourreU}
E^U(\Theta)\,U^*[A,U]\,E^U(\Theta)\ \ge\ c_0\,E^U(\Theta)+\mathsf K.
\end{equation}
In particular, a Mourre estimate for $U$ holds on every open set $\Theta\subset\T$ whose closure does not contain $1$.
\end{klemma}

\begin{proof}
\emph{Choice of cut-offs.} Since $\theta\mapsto2\Delta\theta/(1+\theta)$ is increasing, we may fix $\theta_0\in(0,\theta_1)$ with $2\Delta\theta_0/(1+\theta_0)\ge c_0$. Choose $\theta_1'\in(\theta_0,\theta_1)$ and $\theta_2'\in(\theta_2,2\pi/\Delta)$. Pick $\varphi\in C(\T;[0,1])$ with $\varphi=1$ on $\Theta_{[\theta_1,\theta_2]}$ and $\supp\varphi\subset\Theta_{[\theta_1',\theta_2']}$, and $\eta\in C_c(\R;[0,1])$ with $\eta=1$ on $[0,\theta_0]$ and $\supp\eta\subset(-1,\theta_1')$.

\emph{Step 1: $\varphi(U_0)\eta(H_0)=0$.} We have $\varphi(U_0)\eta(H_0)=(\varphi\circ\varepsilon_\Delta\cdot\eta)(H_0)$. If $E\ge0$ and $\eta(E)\ne0$ then $E\in[0,\theta_1')$, so $\varepsilon_\Delta(E)\in\Theta_{[0,\theta_1')}$, which is disjoint from $\Theta_{[\theta_1',\theta_2']}$ because $\varepsilon_\Delta$ is injective on $[0,2\pi/\Delta)$. Hence $\varphi(\varepsilon_\Delta(E))\eta(E)=0$ for all $E\ge0$.

\emph{Step 2: splitting the commutator.} By \eqref{eq:UAU2},
\[
\varphi(U)\,U^*[A,U]\,\varphi(U)=2\Delta\,\varphi(U)(1-R)\varphi(U)+\varphi(U)\mathsf K_1\varphi(U).
\]

\emph{Step 3: high energies.} Write $R=\zeta(H_0)+R\,\eta(H_0)$ with $\zeta(E)=(1+E)^{-1}(1-\eta(E))$. For $E\in[0,\theta_0]$, $\zeta(E)=0$; for $E>\theta_0$, $0\le \zeta(E)\le(1+E)^{-1}<(1+\theta_0)^{-1}$. Hence $1-\zeta(H_0)\ge\theta_0/(1+\theta_0)$, and
\[
\varphi(U)(1-R)\varphi(U)\ \ge\ \frac{\theta_0}{1+\theta_0}\,\varphi(U)^2-\varphi(U)\eta(H_0)R\,\varphi(U),
\]
where we used that $R$ and $\eta(H_0)$ commute.

\emph{Step 4: the term with the cut-off is compact.} By Step 1 and Lemma~\ref{lem:lowenergy}, $\varphi(U)\eta(H_0)=\varphi(U_0)\eta(H_0)+(\varphi(U)-\varphi(U_0))\eta(H_0)\in\KK$. Hence $\varphi(U)\eta(H_0)R\varphi(U)\in\KK$.

\emph{Step 5: conclusion.} Steps 2--4 give
\[
\varphi(U)\,U^*[A,U]\,\varphi(U)\ \ge\ c_0\,\varphi(U)^2+\mathsf K_2,\qquad \mathsf K_2:=\varphi(U)\mathsf K_1\varphi(U)-2\Delta\,\varphi(U)\eta(H_0)R\varphi(U)\in\KK,
\]
using $2\Delta\theta_0/(1+\theta_0)\ge c_0$ and $\varphi(U)^2\ge0$. Since $\varphi=1$ on $\Theta$, $E^U(\Theta)\varphi(U)=\varphi(U)E^U(\Theta)=E^U(\Theta)$. Sandwiching the inequality between $E^U(\Theta)$ gives \eqref{eq:mourreU} with $\mathsf K=E^U(\Theta)\mathsf K_2E^U(\Theta)$. For the last claim, the closure of such an open $\Theta$ is contained in some $\Theta_{[\theta_1,\theta_2]}$.
\end{proof}

\subsection{Proofs of Theorems~\ref{thm:A} and~\ref{thm:B}(b)}\label{sec:proofA}

\medskip\noindent{\bf Proof of Theorem~\ref{thm:A}.} Assume \eqref{ass:V} with $\rho>0$. By Lemma~\ref{lem:product}, $U\in C^{1+0}(A)$. Let $\Theta\subset\T$ be open with $1\notin\overline\Theta$. By Lemma~\ref{lem:mourre}, \eqref{eq:mourreU} holds on $\Theta$ for some $c_0>0$, so Theorem~\ref{thm:FRT} applies on $\Theta$.

(a) By Theorem~\ref{thm:FRT}(i), $U$ has finitely many eigenvalues in $\Theta$, each of finite multiplicity. A compact subset of $\T\setminus\{1\}$ is contained in such a $\Theta$, so $\sigma_p(U)\setminus\{1\}$ has no accumulation point in $\T\setminus\{1\}$.

(b) By Theorem~\ref{thm:FRT}(ii), $E^U_{\rm sc}(\Theta)=0$ for every such $\Theta$. Since $\T\setminus\{1\}$ is a countable union of such arcs, $E^U_{\rm sc}(\T\setminus\{1\})=0$; and a continuous measure gives no mass to the point $\{1\}$. Hence $\HH_{\rm sc}(U)=\{0\}$.

(c) Let $s\in(1/2,1]$ and let $\Theta'$ be a compact subset of $\T\setminus(\{1\}\cup\sigma_p(U))$. Then $\Theta'$ is contained in an open arc $\Theta$ with $1\notin\overline\Theta$, and is a compact subset of $\Theta\setminus\sigma_p(U)$. By Proposition~\ref{prop:weights}, $\jap{x}^{-s}$ extends to an element of $\mathscr B(\mathrm{Dom}(\jap{A}^s)^*,\HH)$, so Theorem~\ref{thm:FRT}(iii) gives $\sum_n\|\jap{x}^{-s}U^nE^U(\Theta')f\|^2\le c_{\Theta'}\|f\|^2$. \qed

\medskip\noindent{\bf Completeness.} We prove Theorem~\ref{thm:B}(b). Fix $\Delta>0$ and $s\in(1/2,1]$ with $2s\le\rho$, possible since $\rho>1$. Let $\mathcal O=\T\setminus(\{1\}\cup\sigma_p(U))$, open by Theorem~\ref{thm:A}(a).

\begin{klemma}[Smooth building blocks]\label{lem:smoothblocks}
Let $X$ be a finite product of kicks and free flows. Then $\jap{x}^{-s}X$ is locally $U$-smooth on $\mathcal O$ and locally $U_0$-smooth on $\T\setminus\{1\}$.
\end{klemma}

\begin{proof}
$X^*\in C^1(A)$, so $X^*$ maps $\mathrm{Dom}(A)$ into itself boundedly in the graph norm \cite[Ch.~5]{ABG96}, hence by interpolation $\jap{A}^sX^*\jap{A}^{-s}\in\BB$. Thus, by Proposition~\ref{prop:weights},
\[
\|\jap{x}^{-s}Xu\|\le\|\jap{A}^sX^*\jap{A}^{-s}\|\,\|\jap{A}^s\jap{x}^{-s}\|\,\|\jap{A}^{-s}u\|,
\] i.e., $\jap{x}^{-s}X$ extends to $\mathscr B(\mathrm{Dom}(\jap{A}^s)^*,\HH)$. Apply Theorem~\ref{thm:FRT}(iii) with Lemma~\ref{lem:mourre} for $U$, and Lemma~\ref{lem:free}(d) for $U_0$.
\end{proof}

\begin{proposition}\label{prop:vectors}
Let $\Theta'\subset\mathcal O$ be compact and $f\in E^U(\Theta')\HH$. With $\nu_i,Y_i$ as in \eqref{eq:telescope-comp}:
\\[1mm]
\textup{(a)} $\sum_{n\in\Z}\|\jap{x}^{-s}Y_iU^nf\|^2<\infty$ for every $i$;
\\[1mm]
\textup{(b)} $\tau(N):=\sup_{|k|\ge N}\sum_i\|\nu_iY_iU^kf\|\to0$ as $N\to\infty$;
\\[1mm]
\textup{(c)} for every $\psi\in C^\infty(\T)$ with $\psi=1$ on $\Theta'$, $\|(1-\psi(U_0))U^nf\|\to0$ as $|n|\to\infty$.
\end{proposition}

\begin{proof}
(a) is Lemma~\ref{lem:smoothblocks}, and (b) follows since $|\nu_i|\le C\jap{x}^{-\rho}\le C\jap{x}^{-s}$. (c) Since $(1-\psi(U))U^nf=0$,
$(1-\psi(U_0))U^nf=\sum_m\hat\psi_m(U^m-U_0^m)U^nf$, with $\sum_m|m||\hat\psi_m|<\infty$. For $m\ge1$, $U^m-U_0^m=\sum_{j<m}U_0^{m-1-j}(U-U_0)U^j$ and \eqref{eq:telescope-comp} give $\|(U^m-U_0^m)U^nf\|\le\sum_{j<m}\sum_i\|\nu_iY_iU^{n+j}f\|\le m\,\tau(|n|-m)$, with $\tau(N)=\tau(0)$ for $N\le0$. For $m=-l<0$ use $U^{-l}-U_0^{-l}=-U_0^{-l}(U^l-U_0^l)U^{-l}$. Hence $\|(1-\psi(U_0))U^nf\|\le\sum_m|m||\hat\psi_m|\tau(|n|-|m|)\to0$ by dominated convergence.
\end{proof}

\begin{proof}[Proof of Theorem~\ref{thm:B}(b)]
Let $\Theta'\subset\mathcal O$ be compact, $f\in E^U(\Theta')\HH$, $\Theta''\subset\T\setminus\{1\}$ a closed arc with $\Theta'\subset\operatorname{int}\Theta''$, and $\psi\in C^\infty(\T;[0,1])$ with $\psi=1$ on $\Theta'$, $\supp\psi\subset\Theta''$. Put $f_n=U_0^{-n}U^nf$ and $\tilde f_n=\psi(U_0)f_n$. By Proposition~\ref{prop:vectors}(c), $f_n-\tilde f_n\to0$. By \eqref{eq:telescope-comp}, $\tilde f_{n+1}-\tilde f_n=\sum_i\psi(U_0)U_0^{-n-1}X_i\nu_iY_iU^nf$. Hence for $h\in\HH$ and $0\le N<M$, with $\nu_{i,s}=\jap{x}^{2s}\nu_i$ bounded,
\[
|\ip{h}{\tilde f_M-\tilde f_N}|\le\sum_i\|\nu_{i,s}\|_\infty\Big(\sum_n\|\jap{x}^{-s}X_i^*U_0^{n+1}E^{U_0}(\Theta'')\bar\psi(U_0)h\|^2\Big)^{1/2}\Big(\sum_{n=N}^{M-1}\|\jap{x}^{-s}Y_iU^nf\|^2\Big)^{1/2}.
\]
The first factor is $\le c\|h\|$ by Lemma~\ref{lem:smoothblocks} applied to $X_i^*$ and $U_0$, and the second tends to $0$ by Proposition~\ref{prop:vectors}(a). So $(f_n)$ converges as $n\to+\infty$, and similarly as $n\to-\infty$.

$\mathcal O$ is a countable union of closed arcs, $\sigma_p(U)$ is countable, and $\HH_{\rm sc}(U)=\{0\}$, so $E^U(\mathcal O)\HH=\HH_{\rm ac}(U)$. Finite sums of vectors as above are dense in $\HH_{\rm ac}(U)$, and $\|U_0^{-n}U^n\|=1$, so $W^U_\pm$ exist.

For $f\in\HH_{\rm ac}(U)$, $\|f-U^{-n}U_0^nW^U_\pm f\|=\|U_0^{-n}U^nf-W^U_\pm f\|\to0$, so $f=\Omega^U_\pm W^U_\pm f$ and $\HH_{\rm ac}(U)\subset\Ran\Omega^U_\pm$; with part (a), equality holds. Then $W^U_\pm=(\Omega^U_\pm)^*$, $\HH=\HH_{\rm ac}(U)\oplus\HH_{\rm pp}(U)$ by Theorem~\ref{thm:A}(b), $S_U$ is unitary with $S_UU_0=U_0S_U$, and $U|_{\HH_{\rm ac}(U)}$ is unitarily equivalent to $U_0$, so $\sigma_{\rm ac}(U)=\T$.
\end{proof}

%=====================================================================
\section{Convergence theory}\label{sec:convergence}
%=====================================================================

The proofs of Theorems~\ref{thm:conv}--\ref{thm:S} combine decay estimates for the exact solution (section~\ref{sec:decay}), localization of the one-step defects (section~\ref{sec:defects}), and summation and reduction principles (section~\ref{sec:summation}); weighted estimates for the flows are collected in Fact~\ref{fact:weights}.

\subsection{Decay of scattering states}\label{sec:decay}
The one-step defects are summed along the exact solution, and their decay in time comes from local decay of $e^{-itH}$ on energy-localized data.

\begin{proposition}[Local decay]\label{prop:localdecay}
Assume \eqref{ass:V}. Let $g\in C_c^\infty((0,\infty))$, $\mu\ge0$, $\epsilon>0$ and $k\ge0$. There is $C$ such that for all $t\in\R$
\[
\big\|\jap{p}^k\jap{x}^{-\mu}e^{-itH}g(H)\jap{x}^{-\mu}\big\|\le C\jap{t}^{-\mu+\epsilon}.
\]
\end{proposition}

\begin{proof}
For $\mu=0$ there is nothing to prove. Let $\mu>0$ and $k=0$. The result of \cite[Theorem~5.1]{JMP84} is stated for $-\frac12\Delta_x+V$; it applies to $H$ after rescaling time, since $e^{-itH}=e^{-i(2t)H/2}$ and $H/2=-\frac12\Delta_x+\frac12V$. Under \eqref{ass:V} its hypotheses hold for every $m\in\N$, with their $\varepsilon$ equal to $\rho$: $V\in C^{2m+1}$ and $|D^\alpha V(x)|\le C(1+|x|)^{-(|\alpha|+\rho)}$ for $|\alpha|\le2m+1$. Since $H$ has no positive eigenvalues \cite{FH82}, it gives, for $\eta>0$ and $s\in(0,2m-\frac12+\eta)$,
\[
\|\jap{x}^{-s}e^{-itH}g(H)\jap{x}^{-s}\|\le c_s(1+|t|)^{-\beta(s)},\qquad\beta(s)=\frac{s(2m-1)}{2m-\frac12+\eta}.
\]
Given $\mu$ and $\epsilon$, choose $\eta=\frac12$ and $m$ so large that $2m>\mu$ and $\mu\,\frac{2m-1}{2m}\ge\mu-\epsilon$; then $s=\mu$ is admissible and $\beta(\mu)\ge\mu-\epsilon$, which is the claim for $k=0$. For $k>0$ insert $\tilde g(H)$, $\tilde g\in C_c^\infty(\R)$ with $\tilde g=1$ on $\supp g$, and use that $\jap{p}^k\jap{x}^{-\mu}\tilde g(H)\jap{x}^\mu$ is bounded (Fact~\ref{fact:weights}(b)).
\end{proof}

\begin{klemma}[Outgoing scattering states, uniform version]\label{lem:scattstates}
Assume \eqref{ass:V} with $\rho>1$, and let $\chi\in C_c^\infty(\R^d\setminus\{0\})$ and $k\in\N$. Then $\Omega^c_\pm J_\chi$ exists in operator norm,
\[
\Omega^c_+J_\chi=J_\chi+i\int_0^\infty e^{i\tau H}Ve^{-i\tau H_0}J_\chi\,d\tau ,
\]
and for $t\ge0$
\[
e^{-itH}\Omega^c_+J_\chi=e^{-itH_0}J_\chi+R_+(t),\qquad\|R_+(t)\|_{L^2\to \mathcal H^k}\le C\jap{t}^{1-\rho},
\]
with $C$ depending on $\chi$, $k$ and $V$. The same holds for $t\le0$ with $\Omega^c_-$ and $\jap{t}=\jap{|t|}$.
\end{klemma}

\begin{proof}
We use Cook's method~\cite{RS79}, based on $\frac{d}{d\tau}e^{i\tau H}e^{-i\tau H_0}J_\chi=ie^{i\tau H}Ve^{-i\tau H_0}J_\chi$. By the Leibniz rule, $\partial^\beta(Ve^{-i\tau H_0}J_\chi)$ is a sum of terms $(\partial^{\beta_1}V)\,e^{-i\tau H_0}(i\xi)^{\beta_2}\chi(p)\jap{x}^{-\rho}$ with $|\partial^{\beta_1}V|\le C\jap{x}^{-\rho}$, so Lemma~\ref{prop:decay} gives $\|Ve^{-i\tau H_0}J_\chi\|_{L^2\to \mathcal H^k}\le C\jap{\tau}^{-\rho}$, which is integrable. This gives the existence in operator norm and the formula. Then $R_+(t)=i\int_t^\infty e^{-i(t-\tau)H}Ve^{-i\tau H_0}J_\chi\,d\tau$, and Fact~\ref{fact:sobolev} gives $\|R_+(t)\|_{L^2\to \mathcal H^k}\le C\int_t^\infty\jap{\tau}^{-\rho}d\tau\le C'\jap{t}^{1-\rho}$.
\end{proof}

\subsection{One-step defects}\label{sec:defects}

We now give exact representations of the defect $U-e^{-i\Delta H}$, explicitly for Lie and Strang splitting, and structurally for a general composition. For the general case we write the composition with step $s$ as $U(s)=T_m(s)\cdots T_1(s)$, $T_i(s)=e^{-isc_iZ_i}$ with $Z_i\in\{H_0,V\}$ and $c_i\in\{a_j,b_j\}$, and put $L_i(s)=T_m(s)\cdots T_{i+1}(s)$ and $c_{\max}=\max_i|c_i|$.

\begin{klemma}[One-step defects]\label{lem:defect}
Assume \eqref{ass:V}; in (a) $U=U_{\rm L}$.
\\[1mm]
\textup{(a)} \emph{Lie.} For $\phi\in \mathcal H^2$,
\[
\|(U-e^{-i\Delta H})\phi\|\le\int_0^\Delta\!\!\int_0^s\big\|[H_0,V]\,e^{-i(s-\sigma)H_0}e^{-isV}\phi\big\|\,d\sigma\,ds,
\qquad [H_0,V]=-(\Delta_xV)-2\nabla V\cdot\nabla .
\]
\par\vspace{1mm}\noindent \textup{(b)} \emph{Strang.} Let $a=V/2$ and $C_a=[H_0,[H_0,a]]$. Then
\begin{equation}\label{eq:Ca}
C_a=4\sum_{j,k}(\partial_j\partial_ka)\,\partial_j\partial_k+4\,\nabla(\Delta_xa)\cdot\nabla+(\Delta_x^2a),
\end{equation}
and for $\phi\in \mathcal H^2$,
\[
\|(U_S-e^{-i\Delta H})\phi\|\le\int_0^\Delta\Big(s^2\big\||\nabla a|^2e^{-isH_0}e^{-isa}\phi\big\|+\int_0^s(s-\sigma)\big\|C_a\,e^{-i(s-\sigma)H_0}e^{-isa}\phi\big\|\,d\sigma\Big)ds.
\]
In particular, with suprema over $0\le\tau\le s\le\Delta$,
\begin{align}
\|(U-e^{-i\Delta H})\phi\|&\le\tfrac{\Delta^2}2\,\sup\big\|[H_0,V]e^{-i\tau H_0}e^{-isV}\phi\big\|,\label{eq:defL}\\
\|(U_S-e^{-i\Delta H})\phi\|&\le\tfrac{\Delta^3}3\,\sup\big\||\nabla a|^2e^{-i\tau H_0}e^{-isa}\phi\big\|+\tfrac{\Delta^3}6\,\sup\big\|C_a\,e^{-i\tau H_0}e^{-isa}\phi\big\|.\label{eq:defS}
\end{align}
\end{klemma}

\begin{proof}
We prove the identities on $\Sch$, which all operators below preserve; the inequalities extend to $\mathcal H^2$ by density. If $u\in C^1([0,\Delta];\HH)$ takes values in $\mathrm{Dom}(H)$ and $u'=-iHu-iD(s)u$, then $\frac{d}{ds}(e^{isH}u(s))=-ie^{isH}D(s)u(s)$, so that
\begin{equation}\label{eq:duhamel}
u(\Delta)-e^{-i\Delta H}u(0)=-i\int_0^\Delta e^{-i(\Delta-s)H}D(s)u(s)\,ds,\qquad \|u(\Delta)-e^{-i\Delta H}u(0)\|\le\int_0^\Delta\|D(s)u(s)\|\,ds.
\end{equation}

(a) For $u(s)=e^{-isH_0}e^{-isV}\phi$, $u'=-iHu-iD(s)u$ with $D(s)=e^{-isH_0}Ve^{isH_0}-V=-i\int_0^se^{-i\sigma H_0}[H_0,V]e^{i\sigma H_0}d\sigma$, and $e^{i\sigma H_0}u(s)=e^{-i(s-\sigma)H_0}e^{-isV}\phi$. Insert in \eqref{eq:duhamel}.

(b) For $u(s)=e^{-isa}e^{-isH_0}e^{-isa}\phi$, $u(\Delta)=U_{\rm S}\phi$ and, since $H=H_0+2a$, $u'=-iHu-iD_S(s)u$ with
\[
D_S(s)=\big(e^{-isa}H_0e^{isa}-H_0\big)+\big(e^{-isa}B(s)e^{isa}-a\big),\qquad B(s)=e^{-isH_0}ae^{isH_0}.
\]
Since $e^{-isa}pe^{isa}=p+s\nabla a$, the first bracket is $is[H_0,a]+s^2|\nabla a|^2$. Taylor's formula gives $B(s)=a-is[H_0,a]-\int_0^s(s-\sigma)e^{-i\sigma H_0}C_ae^{i\sigma H_0}d\sigma$, and $e^{-isa}[H_0,a]e^{isa}=[H_0,a]+is[[H_0,a],a]=[H_0,a]-2is|\nabla a|^2$. The first-order terms $\pm is[H_0,a]$ cancel, which reflects the second order of Strang splitting:
\begin{equation}\label{eq:DS}
D_S(s)=-s^2|\nabla a|^2-\int_0^s(s-\sigma)\,e^{-isa}e^{-i\sigma H_0}C_ae^{i\sigma H_0}e^{isa}\,d\sigma .
\end{equation}
Since $e^{i\sigma H_0}e^{isa}u(s)=e^{-i(s-\sigma)H_0}e^{-isa}\phi$, \eqref{eq:duhamel} gives (b). Formula \eqref{eq:Ca} follows by expanding $[-\Delta_x,-(\Delta_xa)-2\nabla a\cdot\nabla]$.
\end{proof}

\begin{klemma}[Defect of an order-$r$ composition]\label{lem:defect-comp}
Assume \eqref{ass:V} and let the composition have order $r\ge1$. Then for $\phi\in \mathcal H^{2r+2}$
\begin{equation}\label{eq:defect-comp}
U(\Delta)\phi-e^{-i\Delta H}\phi=-i\int_0^\Delta e^{-i(\Delta-s)H}D(s)U(s)\phi\,ds,\qquad D(s)=\frac1{(r-1)!}\int_0^s(s-\sigma)^{r-1}D^{(r)}(\sigma)\,d\sigma,
\end{equation}
where $D(s)=\sum_ic_iL_i(s)Z_iL_i(s)^{-1}-H$. The operator $D^{(r)}(\sigma)$ is a finite sum, with coefficients depending only on the method, of products of exactly $r+1$ factors from $\{H_0,V\}$ and of factors $e^{-ic\sigma H_0}$, $e^{-ic\sigma V}$ with $|c|\le c_{\max}$. Replacing $V$ by $0$ in every product gives products whose sum is $0$. Consequently, for $\sigma,s\in[0,1]$:
\\[1mm]
\textup{(a)} $\|D^{(r)}(\sigma)U(s)w\|\le C\|w\|_{\mathcal H^{2r+2}}$;
\\[1mm]
\textup{(b)} $D^{(r)}(\sigma)U(s)$ is a finite sum of integrals, over bounded sets of parameters, of products $R_1[H_0,V]R_2$, where $R_1$, $R_2$ are products of exponentials $e^{-icH_0}$, $e^{-icV}$ with $|c|\le C$ and of at most $r$ factors $H_0$ or $V$. The number of free flows $e^{-icH_0}$ in these products is bounded in terms of the method.
\end{klemma}

\begin{proof}
\emph{Formula.} On $\Sch$ (or $\mathcal H^\infty$), $U'(s)=\sum_iT_m\cdots T_{i+1}(-ic_iZ_i)T_i\cdots T_1$ and $T_i$ commutes with $Z_i$, so $iU'(s)U(s)^{-1}=\sum_ic_iL_iZ_iL_i^{-1}=H+D(s)$. Hence $u(s)=U(s)\phi$ solves $u'=-iHu-iD(s)u$, and \eqref{eq:duhamel} gives the first identity.

\emph{Vanishing of low-order derivatives.} The order condition says $\log U(s)=-isH+O(s^{r+1})$ in the free algebra; differentiating, $iU'U^{-1}=H+O(s^r)$ formally. The Taylor coefficients $D^{(k)}(0)$ are non-commutative polynomials in $H_0,V$ which vanish in the free algebra for $k<r$, hence vanish as operators on $\Sch$. Taylor's formula with integral remainder gives the expression for $D(s)$.

\emph{Structure.} Differentiating $L_iZ_iL_i^{-1}$ $r$ times by the Leibniz rule, each derivative produces a factor $-ic_lZ_l$ next to $T_l$ (or $+ic_lZ_l$ next to $T_l^{-1}$); with the factor $Z_i$ this gives $r+1$ factors from $\{H_0,V\}$. If $V=0$, then $U(s)=e^{-isH_0}$ because $\sum b_j=1$, so $D\equiv0$ and the sum of the products vanishes.

\emph{(a)} Each factor $H_0$ maps $\mathcal H^{k+2}\to \mathcal H^k$, $V$ maps $\mathcal H^k\to \mathcal H^k$, and the unitary factors are bounded on every $\mathcal H^k$ uniformly (Fact~\ref{fact:sobolev}).

\emph{(b)} We write $\mathrm{Ad}(T)X=TXT^{-1}$. Since $\sum_ja_j=\sum_jb_j=1$, we have $\sum_ic_iZ_i=H$, hence $D(\sigma)=\sum_ic_i\big(\mathrm{Ad}(L_i)Z_i-Z_i\big)$. Write $L_i=T_m\cdots T_{i+1}$ with $T_l=e^{-i\sigma c_lZ_l}$. Telescoping, and using $\mathrm{Ad}(T_l)Z-Z=-ic_l\int_0^\sigma\mathrm{Ad}(e^{-i\tau c_lZ_l})[Z_l,Z]\,d\tau$,
\[
\mathrm{Ad}(L_i)Z_i-Z_i=-i\sum_{l>i}c_l\int_0^\sigma\mathrm{Ad}\big(T_m\cdots T_{l+1}e^{-i\tau c_lZ_l}\big)[Z_l,Z_i]\,d\tau ,
\]
and $[Z_l,Z_i]$ is $0$ or $\pm[H_0,V]$. Hence $D(\sigma)$ is a finite sum of integrals of products $P[H_0,V]P^{-1}$, with $P$ a product of exponentials. Each derivative in $\sigma$ either hits a factor $e^{-i\sigma cZ}$, which adds one factor $Z$, or the upper limit of integration, which adds none. Multiplying on the right by $U(s)$ adds only exponentials.
\end{proof}

The three representations lead to a single localization statement, which is all that is used later.

\begin{proposition}[Localized defects]\label{prop:defectbounds}
Assume \eqref{ass:V}, let $U$ be a composition of order $r$, and let $0\le\mu\le\rho+1$, or $0\le\mu\le\rho+2$ for Strang splitting. For $k$ large enough, depending on the method, on $\rho$ and on $\mu$, there is $C$ such that for all $\Delta\in(0,1]$
\[
\big\|(U-e^{-i\Delta H})\,\jap{x}^{\mu}\jap{p}^{-k}\big\|\le C\Delta^{r+1}.
\]
\end{proposition}

\begin{proof}
Write $\phi=\jap{x}^\mu\jap{p}^{-k}u$. By Fact~\ref{fact:weights}(a), for $|\tau|,|c|\le T_0$ and a differential operator $Q=\sum_{|\beta|\le\ell}c_\beta(x)\partial^\beta$ with $|\partial^\alpha c_\beta|\le C\jap{x}^{-\mu}$,
\begin{equation}\label{eq:keyweight}
\big\|Q\,e^{-i\tau H_0}e^{-icV}\jap{x}^\mu\jap{p}^{-k}\big\|\le C\qquad(k\ge \mu+\ell).
\end{equation}
\emph{Lie}. By \eqref{eq:defL}, $\|(U-e^{-i\Delta H})\phi\|\le\frac{\Delta^2}2\sup\|[H_0,V]e^{-i\tau H_0}e^{-isV}\phi\|$, and $[H_0,V]=-(\Delta_xV)-2\nabla V\cdot\nabla$ has coefficients $O(\jap{x}^{-\rho-1})=O(\jap{x}^{-\mu})$ with all derivatives; apply \eqref{eq:keyweight} with $\ell=1$.
\emph{Strang}. By \eqref{eq:defS} and \eqref{eq:Ca}, the operators $C_a$ and $|\nabla a|^2$ have coefficients $O(\jap{x}^{-\rho-2})=O(\jap{x}^{-\mu})$ (note $2\rho+2\ge\rho+2$); apply \eqref{eq:keyweight} with $\ell=2$.
\emph{General}. By Lemma~\ref{lem:defect-comp}, $\|(U-e^{-i\Delta H})\phi\|\le\frac{\Delta^{r+1}}{(r+1)!}\sup_{\sigma\le s\le\Delta}\|D^{(r)}(\sigma)U(s)\phi\|$, and by part (b) of that lemma it suffices to bound $\|R_1[H_0,V]R_2\jap{x}^\mu\jap{p}^{-k}\|$ uniformly in the parameters. Since $\mu\le\rho+1$, $[H_0,V]=-(\Delta_xV)-2\nabla V\cdot\nabla=\jap{x}^{-\mu}Q$ with $Q$ of order $1$ with $C_b^\infty$ coefficients. Now $R_1$ is bounded $\mathcal H^{2r}\to L^2$ by Fact~\ref{fact:sobolev}, and, by Fact~\ref{fact:weights}(a), $\jap{p}^{2r}\jap{x}^{-\mu}Q\jap{x}^{\mu}\jap{p}^{-2r-1}$ and $\jap{p}^{2r+1}\jap{x}^{-\mu}R_2\jap{x}^{\mu}\jap{p}^{-k}$ are bounded. For the latter, the weight is moved through the factors of $R_2$ one at a time, at a cost of $\mu$ derivatives for each free flow, two for each factor $H_0$, and none for the factors $V$ and $e^{-icV}$, which commute with $\jap{x}$; this fixes how large $k$ must be.
\end{proof}

\subsection{Summation, reduction and proofs}\label{sec:summation}
The defect bounds of section~\ref{sec:defects} are summed in time by the following principle and transferred to wave operators by the reduction principle (Proposition~\ref{prop:reduction}).

\begin{proposition}[Summation principle]\label{prop:summation}
Let $U$ be unitary and $\psi\in\HH$. Then for $n\ge0$
\begin{gather*}
\|(U^n-e^{-in\Delta H})\psi\|\le\sum_{j=0}^{n-1}\|(U-e^{-i\Delta H})e^{-ij\Delta H}\psi\|,\\
\|(U^{-n}-e^{in\Delta H})\psi\|\le\sum_{j=1}^{n}\|(U-e^{-i\Delta H})e^{ij\Delta H}\psi\|.
\end{gather*}
In particular, if $\|(U-e^{-i\Delta H})e^{-ij\Delta H}\psi\|\le M\Delta^{r+1}\jap{j\Delta}^{-\mu}$ for all $j\in\Z$ with $\mu>1$ and $\Delta\le1$, then $\sup_{n\in\Z}\|(U^n-e^{-in\Delta H})\psi\|\le M\big(2+2\int_0^\infty\jap{t}^{-\mu}dt\big)\Delta^r$.
\end{proposition}

\begin{proof}
The first inequality follows from $U^n-e^{-in\Delta H}=\sum_jU^{n-1-j}(U-e^{-i\Delta H})e^{-ij\Delta H}$ and unitarity. The second follows from $U^{-1}-e^{i\Delta H}=-U^{-1}(U-e^{-i\Delta H})e^{i\Delta H}$ applied in the same way. For the last claim use $\Delta\sum_{j\ge0}\jap{j\Delta}^{-\mu}\le\Delta+\int_0^\infty\jap{t}^{-\mu}dt$.
\end{proof}

\begin{proposition}[Reduction principle]\label{prop:reduction}
Assume \eqref{ass:V} with $\rho>1$. For every composition $U$, $\Delta>0$ and $f\in\HH$, with $\psi_\pm=\Omega^c_\pm f$,
\[
\|(\Omega^U_\pm(\Delta)-\Omega^c_\pm)f\|=\lim_{n\to\pm\infty}\|U^n\psi_\pm-e^{-in\Delta H}\psi_\pm\|.
\]
\end{proposition}

\begin{proof}
$U^{-n}U_0^nf-\psi_\pm=U^{-n}(U_0^nf-e^{-in\Delta H}\psi_\pm)+U^{-n}(e^{-in\Delta H}\psi_\pm-U^n\psi_\pm)$. The first term tends to $0$ by definition of $\Omega^c_\pm$, since $U_0^n=e^{-in\Delta H_0}$, and $U^{-n}$ is unitary.
\end{proof}

\medskip\noindent{\bf Proof of Theorem~\ref{thm:conv}.} \emph{(a)} Let $\psi=g(H)\jap{x}^{-s}u$, fix $\mu\in(1,\min(\rho+1,s)]$ and $\epsilon\in(0,\mu-1)$. Apply Proposition~\ref{prop:defectbounds} with the weight $\mu$ and $k$ large enough. For $j\in\Z$,
\[
(U-e^{-i\Delta H})e^{-ij\Delta H}\psi=\big[(U-e^{-i\Delta H})\jap{x}^\mu\jap{p}^{-k}\big]\big[\jap{p}^k\jap{x}^{-\mu}e^{-ij\Delta H}g(H)\jap{x}^{-s}\big]u .
\]
The first bracket is $O(\Delta^{r+1})$ by Proposition~\ref{prop:defectbounds}. Since $\jap{x}^{-(s-\mu)}\le1$, the second is $O(\jap{j\Delta}^{-\mu+\epsilon})\|u\|$ by Proposition~\ref{prop:localdecay}. Proposition~\ref{prop:summation} concludes.

\emph{(b)} The ranges of $g(H)\jap{x}^{-s}$, $g\in C_c^\infty((0,\infty))$, span a dense subspace of $\HH_{\rm ac}(H)=E^H((0,\infty))\HH$. Combine (a) with $\|U^n-e^{-in\Delta H}\|\le2$.

\emph{(c)} Let $\psi=\Omega^c_+J_\chi u$. By Lemma~\ref{lem:scattstates}, $e^{-ij\Delta H}\psi=e^{-ij\Delta H_0}J_\chi u+R_+(j\Delta)u$ for $j\ge0$. Apply Proposition~\ref{prop:defectbounds} with $\mu=\rho$ and $k$ large enough. Then
\begin{gather*}
\|\jap{p}^k\jap{x}^{-\rho}e^{-ij\Delta H_0}J_\chi\|\le C\jap{j\Delta}^{-\rho},\\
\|\jap{p}^k\jap{x}^{-\rho}R_+(j\Delta)\|\le\|R_+(j\Delta)\|_{L^2\to \mathcal H^k}\le C\jap{j\Delta}^{1-\rho},
\end{gather*}
by Lemmas~\ref{prop:decay} (after commuting $\jap{p}^k$ through the weight) and \ref{lem:scattstates}. Hence the defects are $O(\Delta^{r+1}\jap{j\Delta}^{1-\rho})$, summable for $\rho>2$. Proposition~\ref{prop:summation} gives the bound for $n\ge0$; for $n\le0$ use $\Omega^c_-$ and the $t\le0$ part of Lemma~\ref{lem:scattstates}.

\emph{(d)} Let $\psi=g(H)\jap{x}^{-s}u$. Since $\|(U^n-e^{-in\Delta H})\psi\|=\|U^{-n}e^{-in\Delta H}\psi-\psi\|$ and
\[
U^{-m}e^{-im\Delta H}\psi-U^{-n}e^{-in\Delta H}\psi=\sum_{j=n}^{m-1}U^{-(j+1)}\big(e^{-i\Delta H}-U\big)e^{-ij\Delta H}\psi ,
\]
the summability of the defects shown in (a) makes $U^{-n}e^{-in\Delta H}\psi$ converge as $n\to\pm\infty$, to a limit $\Xi_\pm\psi$ with $\|(\Xi_\pm-1)\psi\|\le C\Delta^r\|u\|$. The general case follows by density, since $\|U^{-n}e^{-in\Delta H}\|=1$. Isometry and $U\Xi_\pm=\Xi_\pm e^{-i\Delta H}$ are immediate. The latter gives $\Xi_\pm^*E^U(B)\Xi_\pm=E^{e^{-i\Delta H}}(B)P_{\rm ac}(H)$ for Borel sets $B\subset\T$, and the spectral measures of $e^{-i\Delta H}$ on $\HH_{\rm ac}(H)$ are images of absolutely continuous measures under $\varepsilon_\Delta$, so $\Ran\Xi_\pm\subset\HH_{\rm ac}(U)$. For Lie splitting, $U_{\rm L}^{-n}=K(\tfrac12)^*U_{\rm S}^{-n}K(\tfrac12)$ and $\|(K(\tfrac12)-1)e^{-in\Delta H}\psi\|\le\frac\Delta2\|Ve^{-in\Delta H}\psi\|\to0$ by Proposition~\ref{prop:localdecay}. Hence $\Xi^{\rm L}_\pm=K(\tfrac12)^*\Xi^{\rm S}_\pm$, and, writing $K(\tfrac12)^*\Xi^{\rm S}_\pm-1-\frac{i\Delta}2V=(e^{i\Delta V/2}-1-\frac{i\Delta}2V)\Xi^{\rm S}_\pm+(1+\frac{i\Delta}2V)(\Xi^{\rm S}_\pm-1)$, we get $\|(\Xi^{\rm L}_\pm-1-\frac{i\Delta}2V)\psi\|\le\frac{\Delta^2}8\|V\|_\infty^2\|\psi\|+(1+\frac\Delta2\|V\|_\infty)\|(\Xi^{\rm S}_\pm-1)\psi\|\le C\Delta^2\|u\|$, which gives the Lie statement.

\emph{Weaker weights.} The condition $s>1$ in (a) comes from the summation of the defects in $\ell^1$ (Proposition~\ref{prop:summation}). For $0<s\le1$ and $s_1>1$, the three-lines theorem, applied to $z\mapsto\ip{h}{(U^n-e^{-in\Delta H})g(H)\jap{x}^{-zs_1}f}$ on $0\le\mathrm{Re}\,z\le1$, interpolates between (a) for $s_1$ and $\|U^n-e^{-in\Delta H}\|\le2$. It gives $\sup_n\|(U^n-e^{-in\Delta H})g(H)\jap{x}^{-s}\|\le C\Delta^{rs/s_1}$. Summing the defects in $\ell^2$ instead would require local smoothing estimates for $U$ uniform in $\Delta$, which the proofs of section~\ref{sec:scattering} do not give. \qed

\medskip\noindent{\bf Proof of Theorem~\ref{thm:wo}.} Let $\psi=g(H)\jap{x}^{-s}u$; by Theorem~\ref{thm:conv}(a), $\|(U^n-e^{-in\Delta H})\psi\|\le C\Delta^r\|u\|$ for all $n\in\Z$.

\emph{(a)} Let $\psi_{\rm pp}=P_{\rm pp}(U)\psi$. For every $n$, $\|\psi_{\rm pp}\|^2=\ip{U^n\psi_{\rm pp}}{U^n\psi}$, so
\[
\|\psi_{\rm pp}\|^2\le|\ip{U^n\psi_{\rm pp}}{e^{-in\Delta H}\psi}|+C\Delta^r\|u\|\,\|\psi_{\rm pp}\| .
\]
The orbit $\{U^n\psi_{\rm pp}\}_n$ is relatively compact, since $\psi_{\rm pp}$ is a norm limit of finite combinations of eigenvectors. Moreover $e^{-in\Delta H}\psi\rightharpoonup0$ as $n\to\infty$, by the Riemann--Lebesgue lemma, since $\psi\in\HH_{\rm ac}(H)$. A weakly null sequence converges to $0$ uniformly on compact sets, so the first term tends to $0$, and $\|\psi_{\rm pp}\|\le C\Delta^r\|u\|$.

\emph{(b)} \emph{Chain rule.} $U^{-n}U_0^nf=U^{-n}e^{-in\Delta H}(e^{in\Delta H}U_0^nf)$, where the bracket converges to $\Omega^c_\pm f\in\HH_{\rm ac}(H)$; hence $\Omega^U_\pm=\Xi_\pm\Omega^c_\pm$, and $\Ran\Xi_\pm=\Ran\Omega^U_\pm=\HH_{\rm ac}(U)$ by Theorem~\ref{thm:B}(b). For $\rho\le1$ we do not know whether $\Ran\Xi_\pm=\HH_{\rm ac}(U)$.
\emph{Wave operators.} Let $f=W^c_+\psi$, so that $\Omega^c_+f=\psi$. By Proposition~\ref{prop:reduction}, $\|(\Omega^U_+-\Omega^c_+)f\|=\lim_n\|(U^n-e^{-in\Delta H})\psi\|\le C\Delta^r\|u\|$.
\emph{Inverse splitting wave operators.} $\|U_0^{-n}U^n\psi-U_0^{-n}e^{-in\Delta H}\psi\|\le C\Delta^r\|u\|$ and $U_0^{-n}e^{-in\Delta H}\psi\to W^c_+\psi$. By Theorem~\ref{thm:B}(b), $U_0^{-n}U^nP_{\rm ac}(U)\psi\to W^U_+\psi$. Also $U_0^{-n}U^n\psi_{\rm pp}\rightharpoonup0$, because $U_0^nh\rightharpoonup0$ ($U_0$ is absolutely continuous) and $\{U^n\psi_{\rm pp}\}$ is relatively compact. Hence $U_0^{-n}U^n\psi\rightharpoonup W^U_+\psi$, and weak lower semicontinuity of the norm gives $\|W^U_+\psi-W^c_+\psi\|\le C\Delta^r\|u\|$. The $-$ case is identical.
\emph{Strong convergence.} $W^c_\pm$ maps $\HH_{\rm ac}(H)$ unitarily onto $\HH$, so the vectors $W^c_\pm g(H)\jap{x}^{-s}u$ are dense in $\HH$; with $\|\Omega^U_\pm\|=1$ this gives $\Omega^U_\pm\to\Omega^c_\pm$ strongly, and similarly $W^U_\pm P_{\rm ac}(H)\to W^c_\pm$.

\emph{(c)} Proposition~\ref{prop:reduction} and Theorem~\ref{thm:conv}(c). \qed

\medskip\noindent{\bf Proof of Theorem~\ref{thm:S}.} Let $\psi=Gu$, $\psi'=Gu'$ and $\Pi=e^{iC_P\Delta H_0}$.

(a) Let $\varphi\in\HH_{\rm ac}(H)$. Replace the factors of $P^{-1}$ one at a time, starting from the right, a kick $e^{ic\Delta V}$ by $1$ and a free flow $e^{i\tau H_0}$ by $e^{i\tau H}$. Since $\|(e^{ic\Delta V}-1)e^{-itH}\varphi\|\le|c|\Delta\|Ve^{-itH}\varphi\|$ and $\|(e^{i\tau H_0}-e^{i\tau H})e^{-itH}\varphi\|\le\int_0^{|\tau|}\|Ve^{-i(t-\tau+\sigma)H}\varphi\|\,d\sigma$, which tend to $0$ as $|t|\to\infty$ by Proposition~\ref{prop:localdecay}, we get $\|(P^{-1}-e^{iC_P\Delta H})e^{-itH}\varphi\|\to0$. With $U^{-n}e^{-in\Delta H}=P\tilde U^{-n}P^{-1}e^{-in\Delta H}$ this gives $\Xi^U_\pm=P\,\Xi^{\tilde U}_\pm e^{iC_P\Delta H}$, hence $\Sigma_U=e^{-iC_P\Delta H}\Sigma_{\tilde U}e^{iC_P\Delta H}$. With $\tilde g\in C_c^\infty((0,\infty))$ equal to $1$ on $\supp g$, $e^{i\tau H}G=GB_\tau$ with $B_\tau=\jap{x}^se^{i\tau H}\tilde g(H)\jap{x}^{-s}$, and $\sup_{|\tau|\le|C_P|}\|B_\tau\|<\infty$ by Fact~\ref{fact:weights}(c). So it suffices to treat $U=\tilde U$, for which
\[
\ip{\psi'}{(\Sigma_{\tilde U}-1)\psi}=\ip{(\Xi^{\tilde U}_+-1)\psi'}{\Xi^{\tilde U}_-\psi}+\ip{\psi'}{(\Xi^{\tilde U}_--1)\psi}=O(\Delta^q)\|u\|\,\|u'\|
\]
by Theorem~\ref{thm:conv}(d).

(b) By Theorem~\ref{thm:wo}(b), $\Omega^U_\pm=\Xi_\pm\Omega^c_\pm$, so $S_U=(\Omega^c_+)^*\Sigma_U\Omega^c_-$ and $\ip{f'}{(S_U-S)f}=\ip{\psi'}{(\Sigma_U-P_{\rm ac}(H))\psi}$, and (a) gives the bound. For the $J_\chi$ version, $S_U=\Pi^*S_{\tilde U}\Pi$ by Theorem~\ref{thm:B}(c), $S$ commutes with $\Pi$, and $|\ip{a'}{a}-\ip{b'}{b}|\le\|a'-b'\|\|a\|+\|b'\|\|a-b\|$ reduce it to Theorem~\ref{thm:wo}(c) for $\tilde U$, with $J_\chi$ replaced by $\Pi J_\chi=e^{iC_P\Delta H_0}J_\chi$. The proof of Theorem~\ref{thm:conv}(c) is unchanged by this replacement, since $\jap{t-\tau}\ge c\jap{t}$ for bounded $\tau$.

For Lie, $\Omega_\pm(\Delta)=K(\tfrac12)^*\Omega^{\rm S}_\pm(\Delta)$ by Theorem~\ref{thm:B}(c), with $K(\tfrac12)^*=e^{i\Delta V/2}$, and $\|e^{i\Delta V/2}-1-\frac{i\Delta}2V\|\le\Delta^2\|V\|_\infty^2/8$. Hence
\[
\Omega_\pm-(1+\tfrac{i\Delta}2V)\Omega^c_\pm=\big(e^{i\Delta V/2}-1-\tfrac{i\Delta}2V\big)\Omega^{\rm S}_\pm+(1+\tfrac{i\Delta}2V)(\Omega^{\rm S}_\pm-\Omega^c_\pm),
\]
and Theorem~\ref{thm:wo}(b) for Strang gives the expansion. \qed

\begin{remark}[Bound states]\label{rem:boundstates}
Theorem~\ref{thm:conv} cannot hold on $\HH_{\rm pp}(H)$. By Theorem~\ref{thm:meta}, the phase of a bound state generically drifts by $e_q\Delta^q$ per unit time, with $q$ the effective order, so that $\|U^n\phi_0-e^{-in\Delta E_0}\phi_0\|$ becomes of order one for $n\Delta\sim\Delta^{-q}$; see section~\ref{sec:bound-results}.
\end{remark}

%=====================================================================
\section{Proofs for bound states}\label{sec:bound}
%=====================================================================

This section proves the results of section~\ref{sec:bound-results}, with $\phi_0$ and $E_0$ as there, and characterizes ionization spectrally.

\begin{proposition}\label{prop:boundproved}
Assume \eqref{ass:V}, let $U$ be a composition, and let $\phi_0$ be as above. Part (c) and the last statement of (b), which involve the splitting wave operators, assume moreover $\rho>1$.
\\[1mm]
\textup{(a)} $e^{-i\Delta E_0}\in\sigma_{\rm ac}(U)=\T$; it is the image under $\varepsilon_\Delta$ of the continuum energies $E_0+2\pi m/\Delta>0$, $m\ge1$.
\\[1mm]
\textup{(b)} With $\lambda_j$ the eigenvalues of $U$,
\[
\lim_{N\to\infty}\frac1N\sum_{n<N}|\ip{\phi_0}{U^n\phi_0}|^2=\sum_j\ip{\phi_0}{E^U(\{\lambda_j\})\phi_0}^2 .
\]
Hence the survival amplitude tends to $0$ in the mean if and only if $\phi_0\perp\HH_{\rm pp}(U)$, i.e., if and only if $\phi_0$ is ionized. Moreover $\|U^nP_{\rm ac}(U)\phi_0-U_0^nW^U_+\phi_0\|\to0$ as $n\to\infty$.
\\[1mm]
\textup{(c)} The spectral measure of $W^U_+\phi_0$ for $U_0$, i.e., the $H_0$-energy distribution of the ionized part folded modulo $2\pi/\Delta$ by $\varepsilon_\Delta$, equals the spectral measure of $P_{\rm ac}(U)\phi_0$ for $U$.
\\[1mm]
\textup{(d)} For all $n\in\Z$, $|\ip{\phi_0}{U_{\rm L}^n\phi_0}-\ip{\phi_0}{U_{\rm S}^n\phi_0}|\le\Delta\|V\|_\infty$.
\end{proposition}

\begin{proof}
(a) By Theorem~\ref{thm:conv}(d), $\Ran\Xi_+$ is a closed subspace of $\HH_{\rm ac}(U)$, and it reduces $U$, since $U\Xi_+=\Xi_+e^{-i\Delta H}$ and $U^{-1}\Xi_+=\Xi_+e^{i\Delta H}$. On it, $U$ is unitarily equivalent to $e^{-i\Delta H}P_{\rm ac}(H)$, whose spectrum is $\varepsilon_\Delta([0,\infty))=\T$ and purely absolutely continuous. Hence $\sigma_{\rm ac}(U)=\T$. (b) The first identity is Wiener's theorem (see \cite{Kat04}) applied to the spectral measure of $\phi_0$ for $U$, whose atoms are the $\ip{\phi_0}{E^U(\{\lambda_j\})\phi_0}$; the right-hand side vanishes if and only if $E^U(\{\lambda_j\})\phi_0=0$ for all $j$. The last statement is the existence of $W^U_+$ (Theorem~\ref{thm:B}(b)), since $\|U^nP_{\rm ac}(U)\phi_0-U_0^nW^U_+\phi_0\|=\|U_0^{-n}U^nP_{\rm ac}(U)\phi_0-W^U_+\phi_0\|$. (c) $W^U_+$ is isometric on $\HH_{\rm ac}(U)$ and $W^U_+U=U_0W^U_+$, and the spectral measures of $U_0$ are the images of those of $H_0$ under $\varepsilon_\Delta$. (d) $U_{\rm L}^n=K(\tfrac12)^*U_{\rm S}^nK(\tfrac12)$ gives $\ip{\phi_0}{U_{\rm L}^n\phi_0}=\ip{K(\tfrac12)\phi_0}{U_{\rm S}^nK(\tfrac12)\phi_0}$, and $2\|(K(\tfrac12)-1)\phi_0\|\le\Delta\|V\|_\infty$.
\end{proof}

\begin{proof}[Proof of Theorem~\ref{thm:meta}]
\emph{Step 1: Taylor expansion.} Write $U(\Delta)=\sum_{k\ge0}\Delta^kA_k$ as a formal power series, where $A_k$ is a non-commutative polynomial of degree $k$ in $H_0$ and $V$. Since $V\in C_b^\infty$, such a polynomial maps $\mathcal H^{m+2k}$ to $\mathcal H^m$. The derivative of order $M+1$ of $\Delta\mapsto U(\Delta)f$ is a finite sum of products of $M+1$ factors $H_0$ or $V$ with free flows and kicks, which are bounded on every $\mathcal H^m$ uniformly for $\Delta\in[0,1]$ (Fact~\ref{fact:sobolev}). Taylor's formula with integral remainder therefore gives, for $M\ge0$ and $\Delta\in(0,1]$,
\begin{equation}\label{eq:taylorU}
\Big\|U(\Delta)f-\sum_{k\le M}\Delta^kA_kf\Big\|\le C_M\Delta^{M+1}\|f\|_{\mathcal H^{2M+2}}\qquad(f\in \mathcal H^{2M+2}).
\end{equation}

\emph{Step 2: modified Hamiltonian.} By the Baker--Campbell--Hausdorff formula, $U(\Delta)=\exp(-i\Delta\hat H(\Delta))$ as formal power series in the free algebra generated by $H_0$ and $V$, where $\hat H(\Delta)=H+\sum_{k\ge r}\Delta^kh_k$ and each $h_k$ is a Lie polynomial in $H_0$ and $V$ \cite{HLW06}. Formally $U(\Delta)^*U(\Delta)=1$, so that $\hat H(\Delta)^*=\hat H(\Delta)$. For Strang splitting, $U(-\Delta)=U(\Delta)^{-1}$, so that $\hat H$ is even in $\Delta$, and $h_2=-\frac1{12}[H_0,[H_0,V]]+\frac1{24}[V,[V,H_0]]$. On $\mathcal H^\infty=\cap_m\mathcal H^m$, each $h_k$ is a differential operator with $C_b^\infty$ coefficients, symmetric on $\mathcal H^\infty$.

\emph{Step 3: Rayleigh--Schr\"odinger series.} Since $\phi_0\in\mathrm{Dom}(H^j)=\mathcal H^{2j}$ for all $j$ (Fact~\ref{fact:sobolev}), $\phi_0\in \mathcal H^\infty$. Let $Q=1-\ip{\phi_0}{\cdot}\phi_0$ and $R_0=(H-E_0)^{-1}Q$. As a bounded function of $H$, $R_0$ is bounded on every $\mathcal H^m$, by the equivalence of norms in the proof of Fact~\ref{fact:sobolev}. Put $h_j=0$ for $j<r$, $\chi_0=\phi_0$, $e_0=E_0$, and recursively
\[
e_k=\sum_{j=1}^k\ip{\phi_0}{h_j\chi_{k-j}},\qquad \chi_k=R_0\sum_{j=1}^k(e_j-h_j)\chi_{k-j}\qquad(k\ge1).
\]
Then $\chi_k\in \mathcal H^\infty$, $\chi_k\perp\phi_0$ for $k\ge1$, and the choice of $e_k$ makes the argument of $R_0$ orthogonal to $\phi_0$, so that $(H-E_0)\chi_k=\sum_{j=1}^k(e_j-h_j)\chi_{k-j}$. Hence the formal series $\chi=\sum_k\Delta^k\chi_k$ and $e=\sum_k\Delta^ke_k$ satisfy $\hat H\chi=e\chi$. The coefficients $e_k$ are real: $e\ip{\chi}{\chi}=\ip{\chi}{\hat H\chi}$, both $\ip{\chi}{\chi}=1+O(\Delta)$ and $\ip{\chi}{\hat H\chi}$ have real coefficients, and so does their quotient. Since $h_j=0$ for $j<r$, $e_k=0$ and $\chi_k=0$ for $1\le k<r$. For Strang splitting, induction on $k$ gives $e_k=0$ and $\chi_k=0$ for odd $k$, and $e_2=\ip{\phi_0}{h_2\phi_0}$. Since $H\phi_0=E_0\phi_0$, $\ip{\phi_0}{[H_0,[H_0,V]]\phi_0}=-\ip{\phi_0}{[V,[H_0,V]]\phi_0}=\ip{\phi_0}{[V,[V,H_0]]\phi_0}$, and $[V,[V,H_0]]=-2|\nabla V|^2$. This gives \eqref{eq:quasiE}.

\emph{Step 4: quasi-mode.} From $\hat H\chi=e\chi$ we get $\hat H^j\chi=e^j\chi$ for all $j$, hence $U(\Delta)\chi=e^{-i\Delta e}\chi$ as formal series. Let $\chi^N=\sum_{k\le N}\Delta^k\chi_k$ and $\lambda_N=e^{-i\Delta E_N(\Delta)}$. Since $U-e^{-i\Delta e}=O(\Delta)$ and $e-E_N=O(\Delta^{N+1})$, the formal series of $(U-\lambda_N)\chi^N$ has no terms of order $\le N+1$. Then \eqref{eq:taylorU} with $M=N+1$, applied to the finitely many vectors $\chi_k$, gives $\|(U-\lambda_N)\chi^N\|\le C_N\Delta^{N+2}$.

\emph{Step 5: conclusion.} Let $f,\psi$ be normalized, $|\lambda|=1$ and $\|(U-\lambda)\psi\|\le\delta$. Then $\|U^n\psi-\lambda^n\psi\|\le|n|\delta$ for $n\in\Z$, by telescoping. Writing $f=\alpha\psi+\beta w$ with $w\perp\psi$ normalized, so that $|\beta|$ is the distance from $f$ to $\C\psi$, and expanding $\ip{f}{U^nf}$, we get $|\ip{f}{U^nf}-\lambda^n|\le2|\beta|^2+2|n|\delta$. For $f=\phi_0$, $\psi=\chi^N/\|\chi^N\|$ and $\lambda=\lambda_N$, $\delta\le C_N\Delta^{N+2}$, and $|\beta|^2=1-\|\chi^N\|^{-2}=O(\Delta^{2r})$, since $\ip{\chi^N}{\phi_0}=1$ and $\chi_k=0$ for $1\le k<r$. This gives \eqref{eq:meta}. If $U=P\tilde UP^{-1}$, then $\ip{\phi_0}{U^n\phi_0}=\ip{P^{-1}\phi_0}{\tilde U^nP^{-1}\phi_0}$ and $\|P^{-1}\phi_0-\phi_0\|\le C\Delta$, since $\|(e^{ic\Delta V}-1)f\|\le|c|\Delta\|V\|_\infty\|f\|$, $\|(e^{ic'\Delta H_0}-1)f\|\le|c'|\Delta\|f\|_{\mathcal H^2}$, and kicks are bounded on $\mathcal H^2$. The bound for $\tilde U$, with $f=P^{-1}\phi_0$ and $|\beta|\le\|P^{-1}\phi_0-\phi_0\|+\|\phi_0-\psi\|\le C\Delta$, gives \eqref{eq:meta} for $U$ with $\Delta^{2r}$ replaced by $\Delta^2$. For Lie and reversed Lie splitting, $\tilde U=U_{\rm S}$ and $P=K(\frac12)^*$ or $P=K(\frac12)$.
\end{proof}

\begin{proof}[Proof of Theorem~\ref{thm:noeig}]
Since $U_{\rm S}=K(\tfrac12)U_{\rm L}K(\tfrac12)^*$ and $U_{\rm R}=K(1)U_{\rm L}K(1)^*$, it suffices to treat $U=U_{\rm L}$. Let $U w=\zeta w$ with $w\in\HH$ and $\zeta=e^{-i\Delta E}$, $E\in\R$, and put $\nu=e^{-i\Delta V}-1$, so that $|\nu|\le\Delta|V|$.

\emph{Step 1: the equation in Fourier variables.} The equation reads $(1+\nu)w=\zeta e^{i\Delta H_0}w$, i.e.,
\begin{equation}\label{eq:eigF}
\big(\zeta e^{i\Delta|\xi|^2}-1\big)\hat w(\xi)=F(\xi),\qquad F=\widehat{\nu w}.
\end{equation}

\emph{Step 2: growth of $F$.} For $z=\xi+i\eta\in\C^d$, $|F(z)|\le\int e^{|\eta||x|}|\nu w|\,dx\le\|w\|\,\|e^{|\eta||x|}\nu\|$. In case (a), with $\supp V\subset\{|x|\le R\}$, this gives $|F(z)|\le Ce^{R|\eta|}$. In case (b), $\|e^{|\eta||x|}\nu\|\le C\big(\int e^{2|\eta||x|-2a|x|^\beta}dx\big)^{1/2}\le Ce^{C|\eta|^{\beta'}}$ with $\beta'=\beta/(\beta-1)<2$. In case (c), maximizing $2|\eta|r-2ar^2$ at $r=|\eta|/(2a)$ gives $|F(z)|\le C\jap{\eta}^{(d-1)/2}e^{|\eta|^2/(4a)}$. In all cases $F$ is entire on $\C^d$.

\emph{Step 3: $F$ vanishes on the resonant spheres.} Let $\rho_m=(E+2\pi m/\Delta)^{1/2}$ for $m\ge m_0$, where $E+2\pi m_0/\Delta>0$. The factor $\zeta e^{i\Delta|\xi|^2}-1$ vanishes on $\{|\xi|=\rho_m\}$, and $|\zeta e^{i\Delta|\xi|^2}-1|\le C_m\,\big||\xi|-\rho_m\big|$ nearby. If $F(\xi_0)\ne0$ for some $|\xi_0|=\rho_m$, then by continuity $|\hat w(\xi)|\ge c\,||\xi|-\rho_m|^{-1}$ on a neighborhood of $\xi_0$, which contradicts $\hat w\in\HH$. Hence $F=0$ on all these spheres.

\emph{Step 4: counting zeros.} Fix $\theta\in S^{d-1}$ and let $g(t)=F(t\theta)$, $t\in\C$, an entire function vanishing at $t=\pm\rho_m$, $m\ge m_0$. Suppose $g\not\equiv0$ and write $g(t)=t^kh(t)$ with $h(0)\ne0$. By Jensen's formula \cite{Lev96},
\[
\sum_{|t_j|<r}\log\frac r{|t_j|}=\frac1{2\pi}\int_0^{2\pi}\log|h(re^{i\phi})|\,d\phi-\log|h(0)|,
\]
the sum running over the zeros $t_j\ne0$ of $h$. The zeros $\pm\rho_m$ alone contribute
\[
\sum_{m\ge m_0,\ \rho_m<r}\log\frac{r^2}{\rho_m^2}\ \ge\ \frac\Delta{2\pi}\int_{\rho_{m_0}^2}^{r^2}\log\frac{r^2}{\sigma}\,d\sigma-O(1)=\frac{\Delta}{2\pi}r^2-O(\log r),
\]
since the $\rho_m^2$ are spaced by $2\pi/\Delta$ and $\sigma\mapsto\log(r^2/\sigma)$ is decreasing. By Step 2, the right-hand side of Jensen's formula is $O(r)$ in case (a) and $O(r^{\beta'})$ in case (b), and in case (c) it is at most $\frac1{2\pi}\int_0^{2\pi}\frac{r^2\sin^2\phi}{4a}\,d\phi+O(\log r)=\frac{r^2}{8a}+O(\log r)$. In cases (a) and (b) this contradicts the quadratic lower bound; in case (c) it does so because $\Delta/(2\pi)>1/(8a)$. Hence $g\equiv0$.

\emph{Step 5: conclusion.} Since $\theta$ is arbitrary, $F=0$ on $\R^d$, so \eqref{eq:eigF} gives $(\zeta e^{i\Delta|\xi|^2}-1)\hat w=0$ almost everywhere. The zero set of the factor is a countable union of spheres, a null set, so $\hat w=0$.
\end{proof}

\begin{proof}[Proof of Corollary~\ref{cor:ion}]
Under the decay conditions of Theorem~\ref{thm:noeig}, \eqref{ass:V} holds for every $\rho$. Indeed, the derivatives of $V$ are bounded, and interpolation gives $|\partial^\alpha V(x)|\le C_\alpha\sup_{|y-x|\le1}|V(y)|^{\theta_\alpha}$ for some $\theta_\alpha>0$, which decays faster than any power. By Theorem~\ref{thm:A}(b), and by Theorem~\ref{thm:noeig} or Proposition~\ref{prop:Snoeig} of SM, $\HH=\HH_{\rm ac}(U)$, so $\Ran\Omega^U_\pm=\HH$ and $W^U_\pm=(\Omega^U_\pm)^*$ is unitary by Theorem~\ref{thm:B}(b). The survival amplitude tends to $0$ by the Riemann--Lebesgue lemma, and the second limit is the existence of $W^U_+$. Finally $\|\one_BU^n\phi_0\|\le\|\one_BU_0^nW^U_+\phi_0\|+o(1)$, and $\|\one_Be^{-itH_0}g\|\to0$ for every $g\in\HH$, by the dispersive bound $\|\one_Be^{-itH_0}g\|\le|B|^{1/2}(4\pi|t|)^{-d/2}\|g\|_{L^1}$ on the dense set $L^1\cap\HH$.
\end{proof}

\begin{remark}\label{rem:noeig}
\textup{(i)} Compact support and hypothesis \eqref{ass:A} are incompatible. For compactly supported $V$, we obtain survival at least up to times $\Delta^{-N}$ (Theorem~\ref{thm:meta}) and eventually ionization, but with no rate. Gaussian potentials satisfy \eqref{ass:A}, so that Theorem~\ref{thm:expmeta} applies to them, and Theorem~\ref{thm:noeig}(c) gives ionization when $\Delta>\pi/(4a)$. \textup{(ii)} The argument does not apply to potentials with exponential decay, since $F$ is then analytic in a strip only, and bounded analytic functions in a strip may vanish at the points $\pm\rho_m$; the P\"oschl--Teller well of the numerical tests in SM is of this type, and its bound state is ionized numerically. \textup{(iii)} Under conditions (a) and (b), Theorem~\ref{thm:noeig} extends to compositions whose kick times are pairwise distinct modulo $\Delta$, possibly after processing, which includes the triple-jump compositions (Proposition~\ref{prop:Snoeig} of SM); we have not examined condition (c) for general compositions.
\end{remark}

Section~\ref{sec:Sbound} of SM discusses the mechanism of ionization, proves Remark~\ref{rem:noeig}(iii), bounds the rate of ionization from below under a non-cancellation hypothesis, and gives numerical evidence. Numerically, the ionized part is emitted at the momenta $|\xi|^2=E_0+2\pi m/\Delta$ predicted by Proposition~\ref{prop:boundproved}(c), at a rate close to a Fermi golden rule, and the fourth-order triple jump ionizes about twice as fast as Lie splitting at the same step.

\appendix
%=====================================================================
\section{Tools}\label{sec:abstract}
%=====================================================================

This appendix collects the results on operators other than the scheme that are used in sections~\ref{sec:scattering}--\ref{sec:bound}. The proofs, or references, are given in section~\ref{sec:Stools} of SM.

\subsection{Abstract results}
\medskip\noindent{\bf Regularity classes.} Let $A$ be self-adjoint in $\HH$. Following \cite{ABG96,FRT13}, $S\in\BB$ belongs to $C^k(A)$ if $t\mapsto e^{-itA}Se^{itA}$ is strongly of class $C^k$.

\begin{fact}[{\cite[\S2.1]{FRT13}, \cite[Ch.~5--6]{ABG96}}]\label{fact:reg}
\textup{(a)} $S\in C^1(A)$ if and only if the form $\mathrm{Dom}(A)\ni f\mapsto\ip{Af}{Sf}-\ip{f}{SAf}$ is continuous for the topology of $\HH$; the associated bounded operator is denoted $[A,S]$.
\\[1mm]
\textup{(b)} $C^1(A)$ and $C^2(A)$ are $*$-subalgebras of $\BB$, and $[A,ST]=[A,S]T+S[A,T]$ for $S,T\in C^1(A)$.
\\[1mm]
\textup{(c)} $S\in C^2(A)$ if and only if $S\in C^1(A)$ and $[A,S]\in C^1(A)$ (cf.\ \cite[Lemma~5.2]{ABC15a} for unitaries).
\\[1mm]
\textup{(d)} $C^2(A)\subset C^{1+0}(A)\subset C^{1,1}(A)\subset C^1(A)$.
\end{fact}

\begin{proposition}[Core criterion]\label{prop:core}
Let $D\subset \mathrm{Dom}(A)$ be a core for $A$ and $S\in\BB$ with $SD\subset D$. Suppose the operator $AS-SA:D\to\HH$ is bounded. Then $S\in C^1(A)$ and $[A,S]$ is the closure of $(AS-SA)|_D$.
\end{proposition}

The proof is standard; see section~\ref{sec:Stools} of SM.

\medskip\noindent{\bf Mourre estimates and local smoothness for unitary operators.} \begin{definition}[\cite{Kat66,FRT13}]\label{def:smooth}
Let $U$ be unitary and $\Theta\subset\T$ open. $B\in\BB$ is \emph{locally $U$-smooth on $\Theta$} if for each compact $\Theta'\subset\Theta$ there is $c_{\Theta'}$ with
\[
\sum_{n\in\Z}\|BU^nE^U(\Theta')f\|^2\le c_{\Theta'}\|f\|^2\qquad(f\in\HH).
\]
\end{definition}

We say that $U\in C^1(A)$ satisfies a \emph{Mourre estimate} on a Borel set $\Theta\subset\T$ if there are $c_0>0$ and $\mathsf K\in\KK$ with
\begin{equation}\label{eq:mourre}
E^U(\Theta)\,U^*[A,U]\,E^U(\Theta)\ \ge\ c_0\,E^U(\Theta)+\mathsf K .
\end{equation}
It is \emph{strict} if $\mathsf K=0$. Note $U^*[A,U]=U^*AU-A$ in form sense. For $U=e^{-iH}$ one has $U^*[A,U]=\int_0^1e^{isH}[iH,A]e^{-isH}\,ds$, so \eqref{eq:mourre} is the time average over one step of the usual Mourre commutator.

\begin{abstractthm}[{\cite[Cor.~2.2, Thm.~2.7, Rem.~2.8, Prop.~2.9]{FRT13}}]\label{thm:FRT}
Let $U$ be unitary, $A$ self-adjoint, $U\in C^{1+0}(A)$, and suppose \eqref{eq:mourre} holds on an open set $\Theta\subset\T$. Then:
\textup{(i)} $U$ has at most finitely many eigenvalues in $\Theta$, each of finite multiplicity; \textup{(ii)} $U$ has no singular continuous spectrum in $\Theta$; if the estimate is strict, $U$ is purely absolutely continuous in $\Theta$; \textup{(iii)} every $B\in\BB$ which extends continuously to an element of $\mathscr B(\mathrm{Dom}(\jap{A}^s)^*,\HH)$ for some $s>1/2$ is locally $U$-smooth on $\Theta\setminus\sigma_p(U)$.
\end{abstractthm}

\begin{fact}\label{fact:compact}
If $a,b\in C_0(\R^d)$ then $a(x)b(p)\in\KK$ (and hence $b(p)a(x)\in\KK$).
\end{fact}

\begin{fact}[Heinz inequality, see, e.g., \cite{Hei51}]\label{fact:heinz}
Let $T_1,T_2\ge0$ be self-adjoint with $\mathrm{Dom}(T_2)\subset \mathrm{Dom}(T_1)$ and $\|T_1u\|\le\|T_2u\|$ for $u\in \mathrm{Dom}(T_2)$. Then for $\alpha\in[0,1]$, $\mathrm{Dom}(T_2^\alpha)\subset \mathrm{Dom}(T_1^\alpha)$ and $\|T_1^\alpha u\|\le\|T_2^\alpha u\|$.
\end{fact}

\medskip\noindent{\bf Sobolev and weighted estimates for the flows.} The following two facts are standard consequences of the Weyl calculus and of the Helffer--Sj\"ostrand functional calculus \cite[Ch.~7--8]{DS99}. Their proofs, with precise references to \cite{DS99}, are given in section~\ref{sec:Stools} of SM.

\begin{fact}[Sobolev regularity of the flows]\label{fact:sobolev}
Assume \eqref{ass:V}. For every $k\in\N$ and $c_{\max}>0$, $\sup_{t\in\R}\|e^{-itH}\|_{\mathcal H^k\to \mathcal H^k}<\infty$ and $\sup_{|c|\le c_{\max}}\|e^{-icV}\|_{\mathcal H^k\to \mathcal H^k}<\infty$; free flows are isometric on $\mathcal H^k$.
\end{fact}

\begin{fact}[Weighted estimates]\label{fact:weights}
Assume \eqref{ass:V}. Let $\mu\ge0$, $k\ge0$, $T_0>0$ and $\tilde g\in C_c^\infty(\R)$.
\\[1mm]
\textup{(a)} For $|\tau|,|c|\le T_0$ and every differential operator $Q$ of order $\ell$ with $C_b^\infty$ coefficients, the operators
\[
\jap{p}^k\jap{x}^{-\mu}e^{-i\tau H_0}\jap{x}^\mu\jap{p}^{-k-\mu},\qquad\jap{p}^ke^{-icV}\jap{p}^{-k},\qquad\jap{p}^k\jap{x}^{-\mu}Q\jap{x}^\mu\jap{p}^{-k-\ell}
\]
are bounded, uniformly in $\tau,c$.
\\[1mm]
\textup{(b)} For every $\mu\in\R$, $\jap{p}^k\jap{x}^{-\mu}\tilde g(H)\jap{x}^{\mu}$ is bounded.
\\[1mm]
\textup{(c)} $\sup_{|\tau|\le T_0}\|\jap{x}^\mu e^{i\tau H}\tilde g(H)\jap{x}^{-\mu}\|<\infty$.
\end{fact}

\subsection{The conjugate operator}\label{sec:conjugate}

Recall $v$, $\gamma=1-\jap{\xi}^{-2}$, $R=(1+H_0)^{-1}$ and the class $\Mfr$ of symbols of order $0$ from section~\ref{sec:kick}. The components $v_j$ are symbols of order $-1$, $\dv v=d\jap{\xi}^{-2}-2|\xi|^2\jap{\xi}^{-4}$ and $\jap{\xi}^{-2}$ are of order $-2$, and $\xi_j\jap{\xi}^{-2}$, $\xi_j\xi_l\jap{\xi}^{-2}$ are of order $-1$ and $0$; all belong to $\Mfr$. If $b\in\Mfr$ then $v\cdot\nabla b$ is a symbol of order $-2$, in particular $v\cdot\nabla b\in\Mfr$. Moreover $|v|\le1/2$, and $\dv v$, $\xi_j\jap{\xi}^{-2}$, $\jap{\xi}^{-2}$ belong to $C_0(\R^d)$.

\begin{proposition}[The group generated by $A$]\label{prop:group}
\textup{(a)} The ODE $\dot\xi=v(\xi)$ has a global flow $(\phi_t)_{t\in\R}$ of $C^\infty$ diffeomorphisms of $\R^d$ with $\phi_t^{-1}=\phi_{-t}$, $|\phi_t(\xi)-\xi|\le|t|/2$, and, for fixed $t$, all derivatives of order $\ge1$ of $\phi_t$ bounded on $\R^d$. The Jacobian $J_t=\det D\phi_t=\exp\int_0^t(\dv v)(\phi_r(\xi))\,dr$ is smooth, satisfies $e^{-c|t|}\le J_t\le e^{c|t|}$ with $c=\|\dv v\|_\infty$, and has bounded derivatives.
\\[1mm]
\textup{(b)} $\widehat{T_tu}=J_t^{1/2}\,(\hat u\circ\phi_t)$ defines a strongly continuous unitary group on $\HH$ with $T_t\Sch=\Sch$.
\\[1mm]
\textup{(c)} Let $A$ be the self-adjoint generator, $T_t=e^{-itA}$. Then $\Sch$ is a core for $A$, and on $\Sch$
\begin{equation}\label{eq:A}
A=\tfrac12\big(x\cdot v(p)+v(p)\cdot x\big)=x\cdot v(p)-\tfrac i2(\dv v)(p),\qquad \widehat{Au}=iv\cdot\nabla\hat u+\tfrac i2(\dv v)\hat u .
\end{equation}
\par\vspace{1mm}\noindent \textup{(d)} For every bounded Borel $b$, $e^{-itA}b(p)e^{itA}=(b\circ\phi_t)(p)$.
\end{proposition}

This construction is that of \cite[\S3.2]{FRT13}. The operator $A$ was introduced in \cite[Def.~1.1]{Yok98}, where its essential self-adjointness on $\Sch$ is also stated, with reference to \cite[Prop.~7.6.3(a)]{ABG96}.

\begin{proposition}[Fourier multipliers]\label{prop:mult}
Let $b\in C^\infty(\R^d)$ have polynomially bounded derivatives, so that $b(p)\Sch\subset\Sch$. Then on $\Sch$
\begin{equation}\label{eq:multcomm}
Ab(p)-b(p)A=i\,(v\cdot\nabla b)(p).
\end{equation}
If in addition $(v\cdot\nabla)^kb$ is bounded for every $k\ge0$, then $b(p)\in C^k(A)$ for all $k$, with $[A,b(p)]=i(v\cdot\nabla b)(p)$. This applies in particular to every $b\in\Mfr$.
\end{proposition}

\begin{proposition}[Comparison of weights]\label{prop:weights}
For $s\in[0,1]$, $\jap{x}^{-s}\HH\subset \mathrm{Dom}(\jap{A}^s)$ and $\jap{A}^s\jap{x}^{-s}\in\BB$. Consequently, $\jap{x}^{-s}$ extends to an element of $\mathscr B(\mathrm{Dom}(\jap{A}^s)^*,\HH)$:
\[
\|\jap{x}^{-s}u\|\le\|\jap{A}^s\jap{x}^{-s}\|\,\|\jap{A}^{-s}u\|\qquad(u\in\HH).
\]
\end{proposition}

\medskip\noindent{\bf Commutators with multiplication operators.} For the class $\Ffr$ of section~\ref{sec:kick} (defined for any $\rho\ge0$), note that $a\in\Ffr$ implies $\partial_l(x_ja)=\delta_{jl}a+x_j\partial_la\in\Ffr$, and that $\Ffr\subset C_0(\R^d)$ when $\rho>0$.

\begin{proposition}[Commutator with a multiplication operator]\label{prop:comm-a}
Let $a\in\Ffr$, $\rho\ge0$, and put $h_j=x_ja$. Then on $\Sch$
\begin{equation}\label{eq:comm-a}
\begin{split}
Aa-aA={}&\frac i2\big\{(\dv v)(p)\,a+a\,(\dv v)(p)\big\}\\
&+\sum_{j=1}^d\Big[-i(\partial_jh_j)R+i\sum_{l=1}^d\Big((p_jp_lR)(\partial_lh_j)R+(p_jR)(\partial_lh_j)(p_lR)\Big)\Big].
\end{split}
\end{equation}
\end{proposition}

\begin{proposition}[An algebra stable under $\mathrm{ad}_A$]\label{prop:algebra}
Let $\rho\ge0$, and let $\Acal$ be the algebra of section~\ref{sec:kick}.
\\[1mm]
\textup{(a)} $\Acal\subset\BB$, $\Acal\Sch\subset\Sch$ and $\Acal^*=\Acal$.
\\[1mm]
\textup{(b)} For $P\in\Acal$ there is $P'\in\Acal$ with $AP-PA=P'$ on $\Sch$.
\\[1mm]
\textup{(c)} $\Acal\subset C^k(A)$ for every $k$, and $\mathrm{ad}_A^k(P)\in\Acal$ for $P\in\Acal$.
\end{proposition}

%=====================================================================
\section*{Acknowledgments}
The author thanks R\'emi Carles for a careful reading of an earlier version of this paper and for helpful comments. This work was supported in part by NSF grant DMS-2606421 and ONR grant N00014-26-1-2017.

\section*{Declarations}
\noindent\emph{Competing interests.} The author has no competing interests to declare.
\emph{Data and code availability.} The scripts reproducing the numerical results, together with the computed data, are provided as ancillary files with the arXiv version of this paper; they are described in section~\ref{sec:Scode} of SM.

\section*{Use of generative AI}
Generative AI (Claude Opus 5.5) was used in the preparation of this work: to carry out the numerical computations and prepare the figures, and to write the scripts reproducing them; to help draft and revise the text and the proofs; to audit the manuscript for internal consistency and gaps in its arguments and proofs; and to help with bibliographical data. All such output was checked by the author, who is alone accountable for the correctness of all statements and proofs.

%=====================================================================
%  SUPPLEMENTARY MATERIAL
%=====================================================================
\clearpage
\setcounter{page}{1}
\renewcommand{\thepage}{S\arabic{page}}
\setcounter{section}{0}\setcounter{subsection}{0}\setcounter{equation}{0}
\setcounter{theorem}{0}\setcounter{figure}{0}\setcounter{table}{0}\setcounter{footnote}{0}
\renewcommand{\thesection}{S\arabic{section}}
\renewcommand{\theequation}{S\arabic{section}.\arabic{equation}}
\renewcommand{\thefigure}{S\arabic{figure}}
\renewcommand{\thetable}{S\arabic{table}}
\renewcommand{\theHsection}{S\arabic{section}}
\renewcommand{\theHequation}{S\arabic{section}.\arabic{equation}}
\renewcommand{\theHfigure}{S\arabic{figure}}
\renewcommand{\theHtable}{S\arabic{table}}
\renewcommand{\theHtheorem}{S\arabic{section}.\arabic{theorem}}
\renewcommand{\refname}{References for the supplementary material}
\begin{center}
\phantomsection\pdfbookmark[0]{Supplementary material}{SMstart}
{\Large\bf Supplementary material for\\ ``Scattering theory and uniform-in-time convergence\\ of splitting methods for the Schr\"odinger equation''\par}
\vspace{1em}
{\large Guillaume Bal\par}
\end{center}
\vspace{1em}

This supplement contains numerical illustrations, further results and proofs on bound states, long-range potentials and other elliptic operators, and standard technical facts. Section~\ref{sec:Sbackground} gives an overview of its methods and of their relation to the literature. Section~\ref{sec:Snumerics} presents numerical illustrations of the results on scattering states, including the aliased channels of the scheme. Section~\ref{sec:Sbound} is devoted to bound states: it describes the mechanism of ionization, proves Theorem~\ref{thm:expmeta}, extends Theorem~\ref{thm:noeig}(a,b) to general compositions, gives a lower bound on the rate of ionization under a non-cancellation hypothesis, and presents numerical evidence. Section~\ref{sec:Slongrange} treats modified splitting wave operators for long-range potentials, and section~\ref{sec:Sscalar} extends the results to elliptic operators $\omega(p)$. Section~\ref{sec:Stools} contains the proofs of the results stated in the appendix of the paper. An index of notation for the paper and the supplement and a description of the numerical codes are given in sections~\ref{sec:Snotation} and~\ref{sec:Scode}. References without the prefix S are to the main paper.

%=====================================================================
\section{Overview of the methods}\label{sec:Sbackground}
%=====================================================================

This section summarizes the methods of the paper and their relation to the literature.

\subsection{Splitting methods over long times}\label{sec:Sbg-split}

We consider $i\partial_t\psi=H\psi$ on $L^2=L^2(\R^d)$, with $H=H_0+V$, $H_0=-\Delta_x$, and a smooth real potential satisfying \eqref{ass:V}, $|\partial^\alpha V(x)|\le C_\alpha\jap{x}^{-\rho-|\alpha|}$ for some $\rho>0$. The two parts of $H$ have explicit flows: $e^{-itH_0}$ is a Fourier multiplier and $e^{-itV}$ a multiplication operator. A splitting method approximates $e^{-i\Delta H}$ by a product of such flows,
\[
U=F_sK_s\cdots F_1K_1,\qquad F_j=e^{-ib_j\Delta H_0},\quad K_j=K(a_j)=e^{-ia_j\Delta V},\qquad \sum_ja_j=\sum_jb_j=1,
\]
the basic examples being the Lie step $U_{\rm L}=e^{-i\Delta H_0}K(1)$, the reversed Lie step $U_{\rm R}=K(1)e^{-i\Delta H_0}$ and the Strang step $U_{\rm S}=K(\frac12)e^{-i\Delta H_0}K(\frac12)$ \cite{S:MQ02,S:HLW06,S:Fao12}. Every such $U$ is unitary, and the free part of the scheme is exact: $U_0=e^{-i\Delta H_0}=F_s\cdots F_1$ and $U_0^n=e^{-in\Delta H_0}$. The method has order $r$ if $U$ agrees with $e^{-i\Delta H}$ up to $O(\Delta^{r+1})$ as a formal power series in $\Delta$ with non-commuting $H_0$ and $V$; methods of order $r\ge3$ need coefficients of both signs \cite{S:BC05}. If $U=P\tilde UP^{-1}$ with a processor $P$ and $\tilde U$ of order $q$, then $U^n=P\tilde U^nP^{-1}$, and $q$ is the effective order of $U$ \cite{S:But69,S:BCR99}. Lie splitting has order $1$ and effective order $2$, since $U_{\rm L}=K(\frac12)^*U_{\rm S}K(\frac12)$.

Each step makes a defect $U-e^{-i\Delta H}$ of size $\Delta^{r+1}$, expressed through commutators of $H_0$ and $V$ \cite{S:JL00,S:Lub08}, and on an interval $[0,T]$ the defects add up to an error $O(T\Delta^r)$, also for unbounded operators \cite{S:BGSvL26}. Such bounds grow linearly in $T$. Long-time behavior is well understood on compact domains, where $H$ has discrete spectrum: the scheme is close to the flow of a modified Hamiltonian, and resonances between eigenvalues of $H_0$, which occur only for exceptional time steps, are the main obstruction \cite{S:DF07,S:Fao12}. On $\R^d$ the situation is different. Scattering states disperse, and since the defects are generated where $V$ acts, one expects them not to accumulate. Bound states do not disperse, and the question is whether the scheme preserves them. The paper addresses both questions.

\subsection{The one-step operator as a scattering problem}\label{sec:Sbg-scat}

Scattering theory compares $e^{-itH}$ with $e^{-itH_0}$ through the wave operators $\Omega^c_\pm$, the strong limits of $e^{itH}e^{-itH_0}$ as $t\to\pm\infty$, which for $\rho>1$ exist and are complete, i.e., $\Ran\Omega^c_\pm=\HH_{\rm ac}(H)$, and through the scattering operator $S=(\Omega^c_+)^*\Omega^c_-$, which maps incoming to outgoing free data and acts energy by energy \cite{S:RS79,S:Yaf92}. For $\rho\le1$ the free comparison dynamics is modified by a phase $\Lambda_\pm(t,p)$ \cite{S:Dol64,S:Hor76,S:DG97} (Hypothesis~(M) in section~\ref{sec:Slrresults}). The same questions arise for a discrete dynamics $n\mapsto U^n$. For unitary $U$ and $U_1$, the wave operators $\Omega_\pm(U,U_1)=\slim_{n\to\pm\infty}U^{-n}U_1^nP_{\rm ac}(U_1)$ are partial isometries intertwining $U_1$ and $U$, with the chain rule $\Omega_\pm(U,U_2)=\Omega_\pm(U,U_1)\Omega_\pm(U_1,U_2)$ \cite{S:Yaf92}. The paper studies the splitting wave operators $\Omega^U_\pm=\Omega_\pm(U,U_0)$, their inverses $W^U_\pm$ and the splitting scattering operator $S_U=(\Omega^U_+)^*\Omega^U_-$, as well as the relative wave operators $\Xi_\pm=\Omega_\pm(U,e^{-i\Delta H})$ and the relative scattering operator $\Sigma_U=\Xi_+^*\Xi_-$ (Theorems~\ref{thm:conv}(d) and~\ref{thm:S}(a)). For $\rho>1$, $\Omega^c_\pm=\Omega_\pm(e^{-i\Delta H},U_0)$, so that $\Omega^U_\pm=\Xi_\pm\Omega^c_\pm$ and $S_U=(\Omega^c_+)^*\Sigma_U\Omega^c_-$ (Theorems~\ref{thm:wo}(b) and~\ref{thm:S}(b)): the scheme affects scattering only through $\Xi_\pm$ and $\Sigma_U$.

\begin{figure}[ht]
\centering
\includegraphics[width=\textwidth]{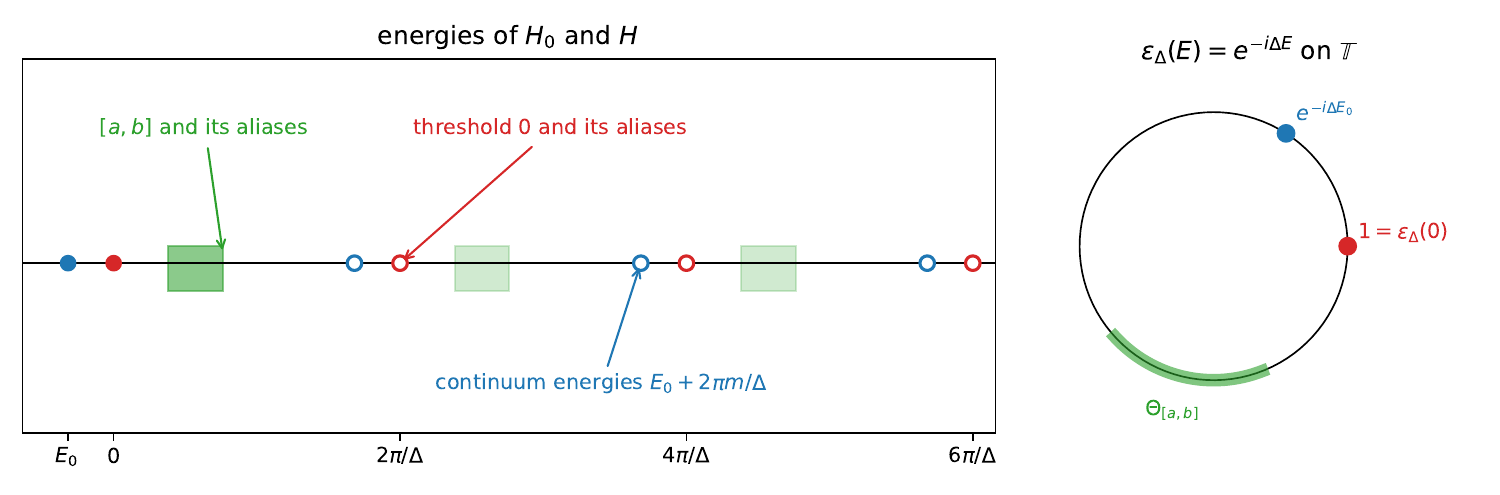}
\caption{The aliasing map $\varepsilon_\Delta$. Left: the threshold $0$ and its aliases $2\pi m/\Delta$ (red), a bound-state energy $E_0$ and the continuum energies $E_0+2\pi m/\Delta$ (blue), and an energy window $[a,b]$ with its aliases (green). Filled: the energies themselves; open or lighter: their aliases, $m\ge1$, which have the same image. Right: the images on $\T$.}
\label{fig:aliasing}
\end{figure}

Two features distinguish the one-step operator from a time-independent Hamiltonian. The first is \emph{aliasing}. Energies are defined only modulo $2\pi/\Delta$: the map $\varepsilon_\Delta(E)=e^{-i\Delta E}$ winds $[0,\infty)$ infinitely often around the unit circle $\T$ (Figure~\ref{fig:aliasing}), so that $\sigma(U_0)=\T$ with infinite multiplicity and spectral projections of $U_0$ on an arc do not localize in energy. The point $1=\varepsilon_\Delta(0)$ is the image of the threshold $0$ of $H_0$ and of its aliases $2\pi k/\Delta$, and plays the role of a threshold. A bound-state energy $E_0<0$ is mapped to the image of the continuum energies $E_0+2\pi m/\Delta>0$, so that bound states of $H$ are embedded in the continuous spectrum of $U$. The second feature is the \emph{lack of compactness}: $U_{\rm L}-U_0=U_0(e^{-i\Delta V}-1)$ is not compact, and in general $U-U_0$ is a sum of such terms.

The spectral theory of $U$ is obtained by the Mourre theory adapted to the setting of unitary operators \cite{S:ABCF06,S:FRT13}. The conjugate operator $A=\frac12(x\cdot v(p)+v(p)\cdot x)$, $v(\xi)=\xi\jap{\xi}^{-2}$ \cite{S:Yok98,S:FRT13}, satisfies $U_0^*[A,U_0]=2\Delta|p|^2\jap{p}^{-2}$, which is positive away from zero kinetic energy. The commutator is a time average of the usual Mourre commutator over one step, $U^*[A,U]=\int_0^1e^{isH}[iH,A]e^{-isH}ds$ for $U=e^{-iH}$. The perturbative results of \cite{S:FRT13} assume that $U-U_0$ is compact, which fails here. The key observation is that $(U-U_0)\eta(H_0)$ is compact for every $\eta\in C_0([0,\infty))$ (Lemma~\ref{lem:lowenergy}). This suffices because high kinetic energies are controlled by the positivity of the free commutator, and compact errors are needed only in terms containing $\eta(H_0)$ (Lemma~\ref{lem:mourre}). The Mourre estimate on arcs avoiding $1$ gives the absence of singular continuous spectrum and local smoothness (Theorem~\ref{thm:A}), and Kato's smooth perturbation argument then gives asymptotic completeness (Theorem~\ref{thm:B}). These results hold at every fixed $\Delta>0$, and the Laplacian enters only through the free flow, which is why they extend to elliptic operators $\omega(p)$ (section~\ref{sec:Sscalar}).

The scheme can also be viewed as a time-periodic problem: $U$ is the monodromy of $H_0+\Delta V\sum_k\delta(t-k\Delta)$, and by Poisson summation of $H+V\sum_{m\ne0}e^{2\pi imt/\Delta}$. Scattering theory for time-periodic Hamiltonians is developed in \cite{S:Yaj77,S:How79,S:Yok98}, and Trotterized dynamics is viewed as a Floquet system in \cite{S:HHZ19}. The harmonics of this perturbation are, however, of the same size as $V$ for all $m$, and the time dependence through Dirac masses is outside the scope of those works. The same commutator methods for unitary operators are used for quantum walks \cite{S:RST19,S:Tie20}; they would also apply to the fully discrete scheme obtained after a spatial discretization.

\subsection{Convergence and bound states}\label{sec:Sbg-results}

The convergence theory rests on three estimates. For a method of order $r$, $0\le\mu\le\rho+1$ ($\mu\le\rho+2$ for Strang splitting) and $k$ large, the one-step defect is localized in space (Proposition~\ref{prop:defectbounds}),
\[
\|(U-e^{-i\Delta H})\jap{x}^{\mu}\jap{p}^{-k}\|\le C\Delta^{r+1};
\]
energy-localized solutions of the exact flow decay locally in time \cite{S:JMP84} (Proposition~\ref{prop:localdecay}),
\[
\|\jap{p}^k\jap{x}^{-\mu}e^{-itH}g(H)\jap{x}^{-\mu}\|\le C\jap{t}^{-\mu+\epsilon};
\]
and, by telescoping, the global error is bounded by the defects along the exact solution (Proposition~\ref{prop:summation}),
\[
\|(U^n-e^{-in\Delta H})\psi\|\le\sum_{j=0}^{n-1}\|(U-e^{-i\Delta H})e^{-ij\Delta H}\psi\|.
\]
For $\psi=g(H)\jap{x}^{-s}u$ and $1<\mu-\epsilon<\mu\le\min(\rho+1,s)$, the first two estimates bound the $j$-th term by $C\Delta^{r+1}\jap{j\Delta}^{-\mu+\epsilon}\|u\|$, and the sum is $O(\Delta^r)\|u\|$ uniformly in $n$. A reduction principle (Proposition~\ref{prop:reduction}) then identifies the error of the wave operators with the limit of this error. The decay is thus provided by the linear scattering theory of $H$, at the level of local decay estimates, which allows any $\rho>0$ and methods of any order. The vector-field method of \cite{S:CS24} for the Lie--Trotter splitting of defocusing NLS provides instead decay for both the exact and the numerical solutions; other long-time results \cite{S:Car26,S:CM26,S:Ign11,S:BCF23,S:BFHJY24} are discussed in the introduction of the paper.

For bound states, the modified Hamiltonian of backward error analysis \cite{S:HLW06,S:Fao12} and Rayleigh--Schr\"odinger perturbation theory give quasi-modes of the scheme. On compact domains the modified Hamiltonian has true eigenvectors; on $\R^d$ it only provides a quasi-mode. The exponential bound of Theorem~\ref{thm:expmeta} uses a compression of the scheme to momenta $|\xi|\le c\Delta^{-1/2}$. Ionization is the numerical counterpart of ionization by time-periodic fields \cite{S:CCLR01,S:SW98}, also called radiation damping \cite{S:SW99}: the harmonic of frequency $2\pi m/\Delta$ couples $\phi_0$ to the continuum at energy $E_0+2\pi m/\Delta$. The absence of eigenvalues is proved by showing that an eigenvector would force an entire function to vanish on infinitely many spheres, in the spirit of the uniqueness results that exclude embedded eigenvalues of Schr\"odinger operators.

%=====================================================================
\section{Numerical illustrations: scattering states}\label{sec:Snumerics}
%=====================================================================

We illustrate in this section the results of the paper on scattering states: the saturation of the error in time (Theorem~\ref{thm:conv}), the aliased channels of the scheme, the limiting error and the splitting scattering operator (Theorems~\ref{thm:conv}(d) and~\ref{thm:S}), slowly decaying and long-range potentials, and the energy dependence of the error. All computations use a Fourier pseudo-spectral discretization in $d=1$, on grids fine enough that spatial errors are negligible and large enough that no significant mass reaches the boundary of the periodic box, and the reference solutions are computed with a much smaller step, by Strang splitting or by the Yoshida triple jump as stated for each test. Besides Lie and Strang splitting we use the fourth-order Yoshida triple jump \cite{S:Yos90}
\[
U_{\rm Y}=S(w_1\Delta)S(w_0\Delta)S(w_1\Delta),\qquad S(h)=e^{-ihV/2}e^{-ihH_0}e^{-ihV/2},\quad w_1=\frac1{2-2^{1/3}},\quad w_0=1-2w_1,
\]
a composition of three Strang steps whose middle step $w_0\Delta\approx-1.70\Delta$ is negative. Unless stated otherwise the initial datum is the normalized packet $\psi=c\,e^{-(x-x_0)^2/8+1.5ix}$, with mean momentum $1.5$ and group velocity $3$. The scripts are described in section~\ref{sec:Scode}.

\begin{figure}[ht]
\centering
\includegraphics[width=\textwidth]{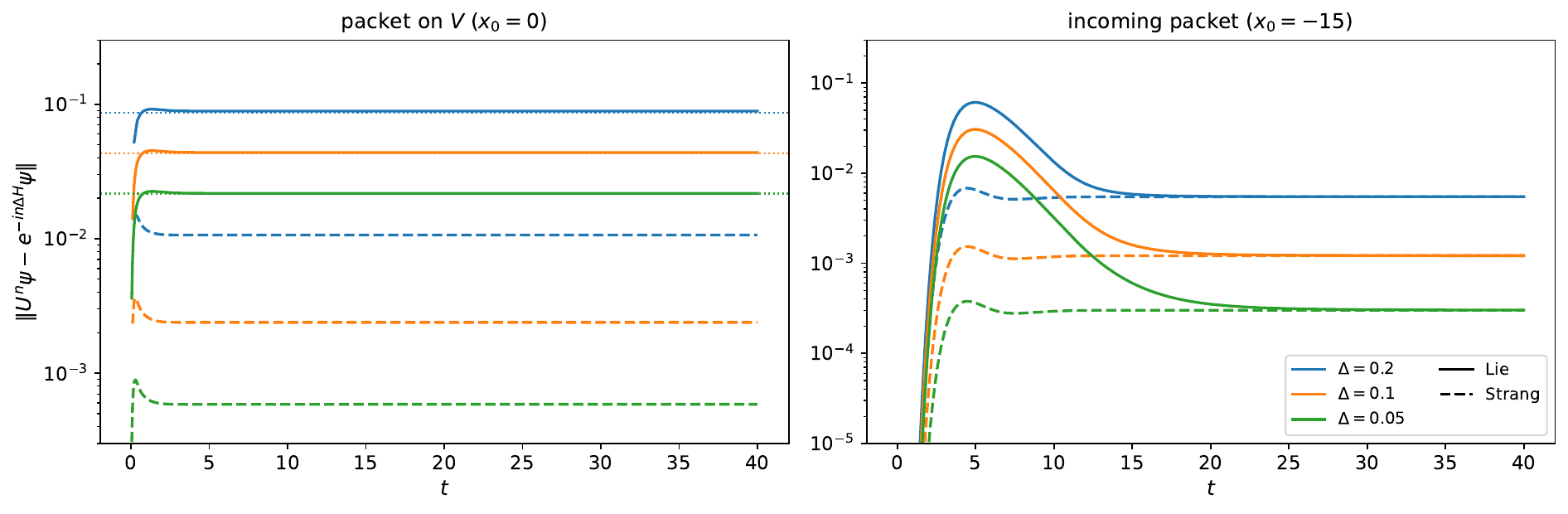}
\caption{Global error of Lie splitting (solid) and Strang splitting (dashed) for $V=1.5\,e^{-x^2}$. Left: packet on the potential, with the Lie plateaus $\frac\Delta2\|V\psi\|$ (dotted). Right: incoming packet.}
\label{fig:errors}
\end{figure}

\paragraph{Saturation of the error in time.} Figure~\ref{fig:errors} shows $\|U^n\psi-e^{-in\Delta H}\psi\|$ for $V=1.5\,e^{-x^2}$, a repulsive potential without bound states, on $[-400,400]$ with $2^{14}$ points and a reference step $0.0025$. For the packet on top of the potential ($x_0=0$, left), the Lie error saturates at once at the value $\frac\Delta2\|V\psi\|$ predicted by Theorem~\ref{thm:conv}(d) ($8.88\cdot10^{-2}$ against $8.66\cdot10^{-2}$ at $\Delta=0.2$, $2.17\cdot10^{-2}$ against $2.17\cdot10^{-2}$ at $\Delta=0.05$), and the Strang error is $O(\Delta^2)$ uniformly in time (maximum $1.54\cdot10^{-2}$, $3.62\cdot10^{-3}$, $8.91\cdot10^{-4}$ for $\Delta=0.2,0.1,0.05$). For the incoming packet ($x_0=-15$, right), the Lie error is $\frac\Delta2\|Ve^{-itH}\psi\|+O(\Delta^2)$, since $U_{\rm L}^n=K(\frac12)^*U_{\rm S}^nK(\frac12)$ and $K(\frac12)\psi\approx\psi$. It shows a transient bump of order $\Delta$ while the packet crosses the potential, and then drops to the Strang value ($5.47\cdot10^{-3}$, $1.21\cdot10^{-3}$ and $3.0\cdot10^{-4}$ for $\Delta=0.2,0.1,0.05$), as predicted by Theorem~\ref{thm:conv}(d), since $\Xi^{\rm L}_+=K(\frac12)^*\Xi^{\rm S}_+$ and $K(\frac12)$ acts trivially far from $V$. For the packet on $V$, applying $K(\frac12)^*=e^{i\Delta V/2}$ once to the initial datum gives Lie errors at $t=40$ equal to the Strang errors to four digits ($1.065\cdot10^{-2}$, $2.39\cdot10^{-3}$, $5.87\cdot10^{-4}$ for $\Delta=0.2,0.1,0.05$; script \texttt{lie\_prekick.py}). For the Yoshida triple jump, with $U_{\rm Y}$ at step $0.0025$ as reference, the packet on $V$ gives maximal errors over $[0,40]$ of $8.3\cdot10^{-2}$, $7.3\cdot10^{-3}$, $4.1\cdot10^{-4}$, $2.8\cdot10^{-5}$ and $1.7\cdot10^{-6}$ for $\Delta=0.4,0.2,0.1,0.05,0.025$, and the error saturates in time. The ratios $14.8$ and $15.9$ for $\Delta\le0.1$ show the fourth order; at $\Delta=0.4$ and $0.2$ the error is dominated by the aliased channels described next.

\begin{figure}[ht]
\centering
\includegraphics[width=\textwidth]{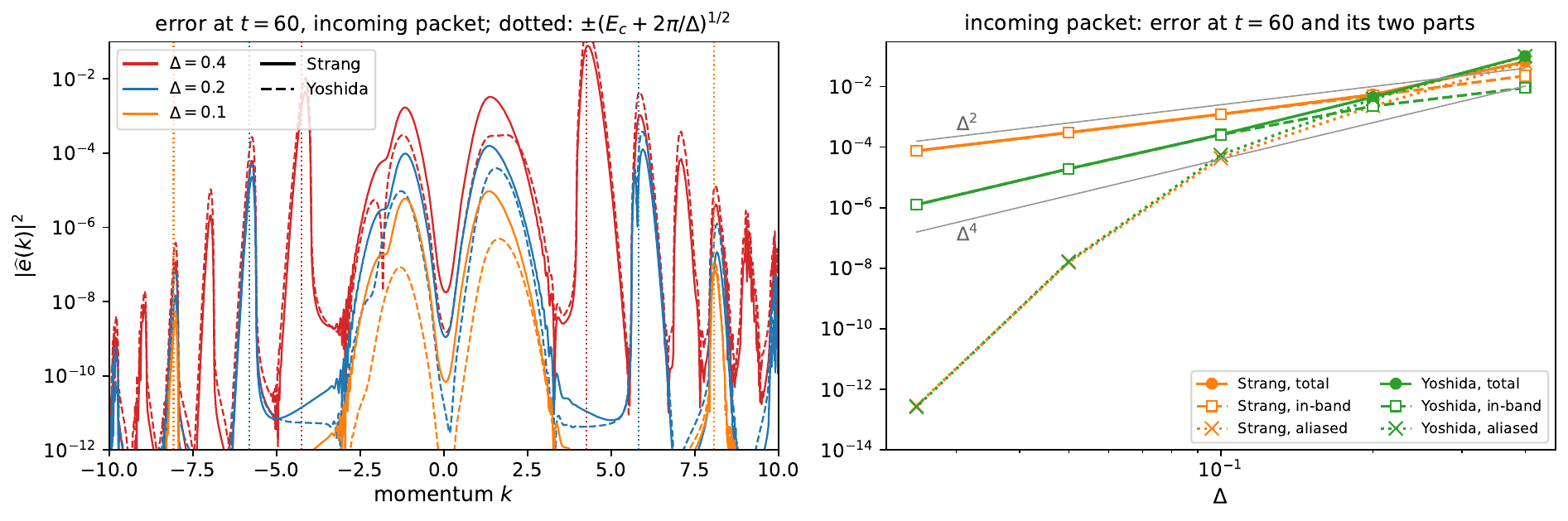}
\caption{Aliased channels for the incoming packet ($x_0=-30$), $V=1.5\,e^{-x^2}$, $t=60$. Left: momentum density of the error $e=U^n\psi-e^{-in\Delta H}\psi$ for Strang splitting (solid) and the Yoshida triple jump (dashed); dotted lines: $\pm(E_c+2\pi/\Delta)^{1/2}$. Right: the error and its in-band and aliased parts against $\Delta$. Grey lines: slopes $\Delta^2$ and $\Delta^4$.}
\label{fig:alias}
\end{figure}

\paragraph{Aliased channels.} Since $S_U$ commutes with $U_0$ but not with $H_0$, the scheme can scatter mass from the energy $E$ of the data to $E+2\pi m/\Delta$, $m\ge1$. Figure~\ref{fig:alias} shows this for the incoming packet ($x_0=-30$) at $t=60$, against the Yoshida triple jump with step $0.0025$. We split the error at the $H_0$-energy $E_c+\pi/\Delta$, $E_c=\ip{\psi}{H\psi}=2.375$, into an in-band and an aliased part. The in-band part is of order $\Delta^2$ for Strang splitting (ratios $4.5$, $4.1$, $4.0$, $4.0$ for $\Delta=0.4,\dots,0.025$) and approaches order $\Delta^4$ for the Yoshida triple jump (ratios $4.0$, $8.8$, $13.4$, $15.3$). The aliased part decays much faster ($6.2\cdot10^{-2}$, $2.3\cdot10^{-3}$, $4.3\cdot10^{-5}$, $1.6\cdot10^{-8}$, $3\cdot10^{-13}$ for Strang splitting, the last value being close to the accuracy of the computation); it dominates the error at $\Delta=0.4$, and for the Yoshida triple jump also at $\Delta=0.2$. At $\Delta=0.4$ and $0.2$ its size is given by single-kick transitions, governed by $\hat V$ at the momentum transfer $\approx(2\pi/\Delta)^{1/2}$ and by the coefficients $\hat\kappa_m$ of the kick comb of the method (section~\ref{sec:Smechanism}): the Yoshida and Strang aliased masses differ by the factors $2.6$ and $2.9$, close to $|\hat\kappa_1|^2=2.43$. At smaller steps this first-order approximation underestimates the aliased part, by factors $3$ and $270$ at $\Delta=0.1$ and $0.05$, and the ratio of the masses drifts to $1.6$ and $1.2$: since $\hat V$ decays like a Gaussian, transitions through several kicks, each with a smaller momentum transfer, dominate. Orders of convergence must therefore be read at steps where the aliased part is negligible, here $\Delta\le0.1$. For the packet on $V$, whose energy spread is larger, the split is less sharp.

\begin{figure}[ht]
\centering
\includegraphics[width=\textwidth]{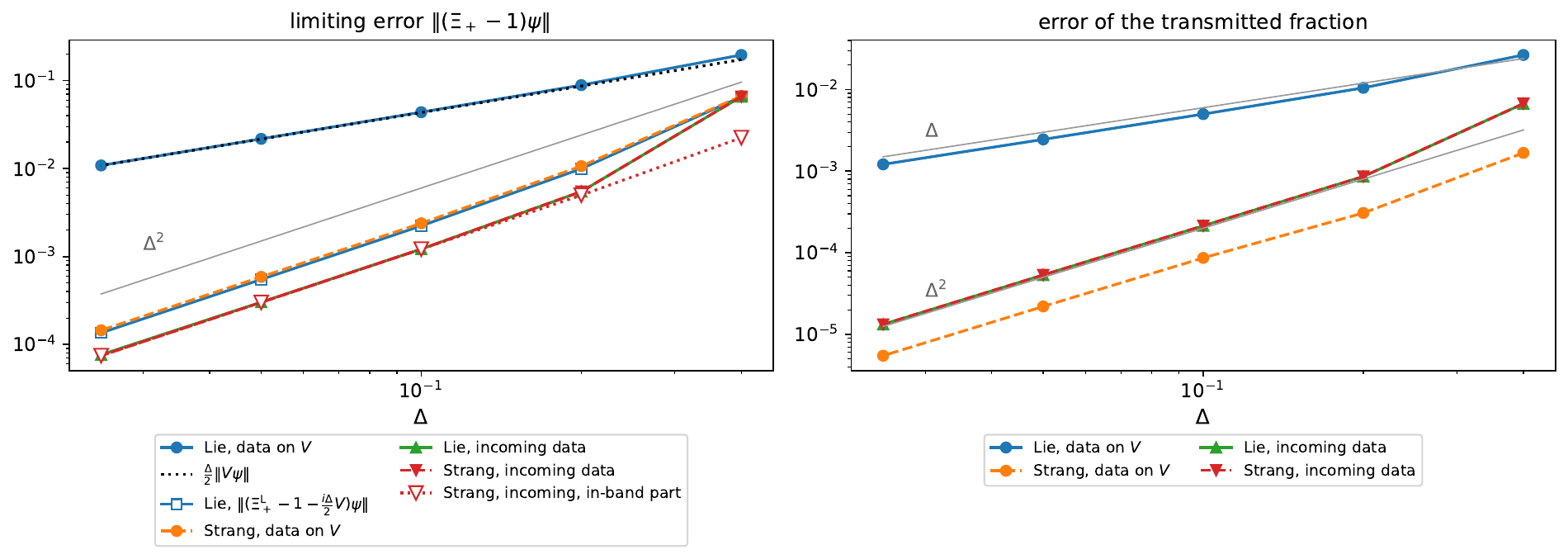}
\caption{Left: limiting error $\|(\Xi_+-1)\psi\|$ against $\Delta$, with the Lie plateau $\frac\Delta2\|V\psi\|$, the corrected Lie quantity of Theorem~\ref{thm:conv}(d), and the in-band part of the Strang error for incoming data (open triangles). Right: error of the transmitted fraction. Grey lines: slopes $\Delta$ and $\Delta^2$.}
\label{fig:rates}
\end{figure}

\paragraph{The limiting error and the splitting scattering operator.} Since $U^n$ is unitary, $\|U^n\psi-e^{-in\Delta H}\psi\|=\|U^{-n}e^{-in\Delta H}\psi-\psi\|$, which converges to $\|(\Xi_+-1)\psi\|$. Figure~\ref{fig:rates} shows this limit, evaluated at $t=60$ when the packets have left the potential, for the packet on $V$ ($x_0=0$) and for an incoming packet ($x_0=-30$). Three effects are visible. First, for data on $V$, the Lie error equals $\frac\Delta2\|V\psi\|$ within $3\%$ for $\Delta\le0.2$, and subtracting the first-order correction of Theorem~\ref{thm:conv}(d) leaves an error of order $\Delta^2$. Second, for incoming data, the Lie and Strang errors coincide and are of order $\Delta^2$: the processor $K(\frac12)$ of Lie splitting acts trivially on incoming and outgoing states, and only the second-order relative scattering operator of Theorem~\ref{thm:S}(a) remains; at $\Delta=0.4$ the error is dominated by the aliased channel, and its in-band part continues the $\Delta^2$ behavior. Third, the same holds for a physical observable, the transmitted fraction $\|\one_{x>0}U^n\psi\|^2$, which for incoming data is determined by the scattering operator: its error is of order $\Delta^2$ for both schemes, whereas for data on $V$ the Lie error is of order $\Delta$.

\begin{figure}[ht]
\centering
\includegraphics[width=\textwidth]{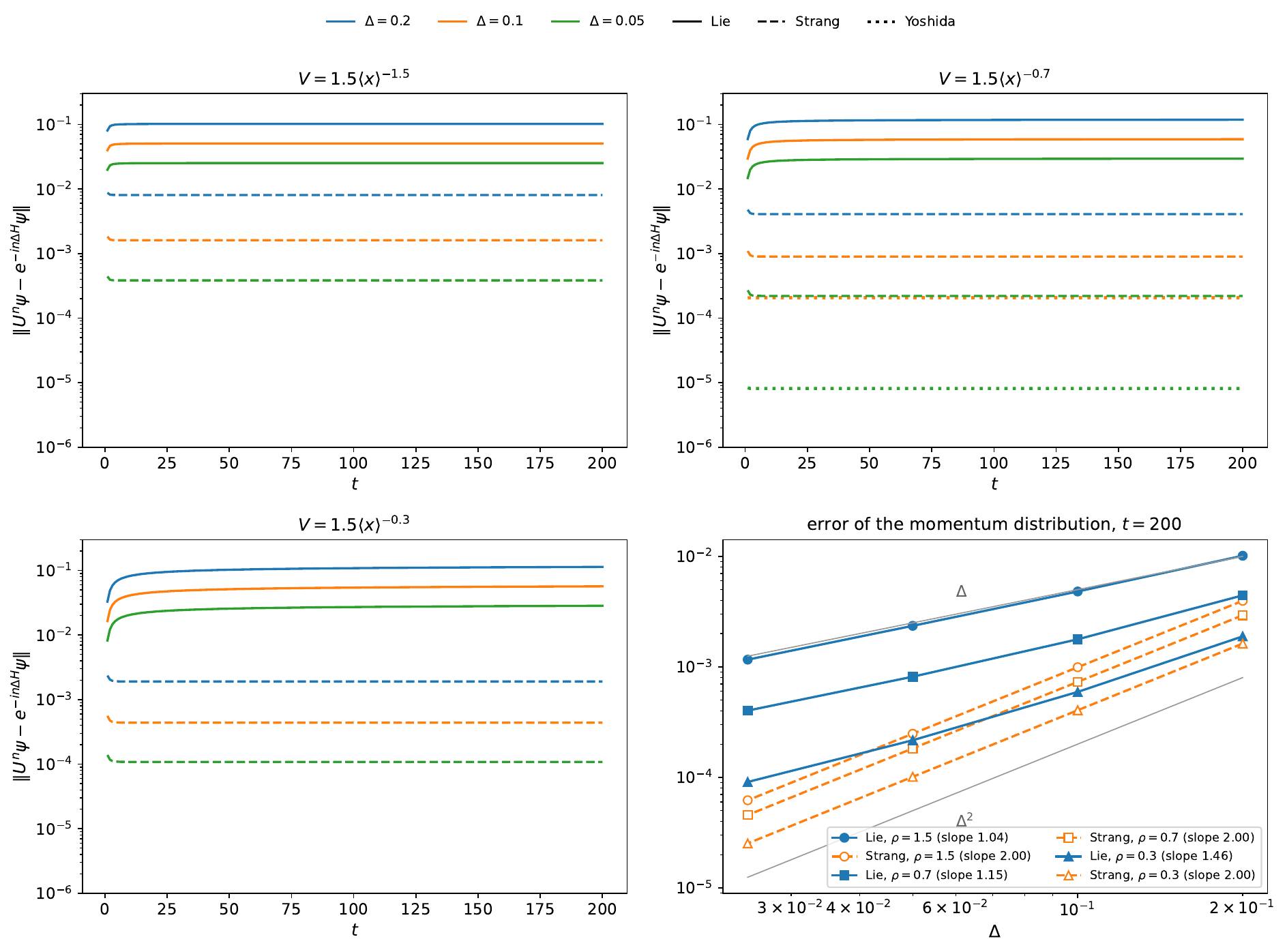}
\caption{$V=1.5\jap{x}^{-\rho}$, packet on $V$. First three panels: global error in time for $\rho=1.5$, $0.7$, $0.3$ (Lie solid, Strang dashed, and for $\rho=0.7$ Yoshida dotted, at $\Delta=0.1$ and $0.05$). Last panel: error of the momentum distribution at $t=200$ against $\Delta$, with least-squares slopes. Grey lines: slopes $\Delta$ and $\Delta^2$.}
\label{fig:slow}
\end{figure}

\paragraph{Slowly decaying and long-range potentials.} We take $V=1.5\jap{x}^{-\rho}$, repulsive, with $\rho=1.5$, $0.7$ and $0.3$, the packet on $V$, the box $[-1600,1600]$ with $2^{16}$ points (less than $10^{-4}$ of the mass lies in $|x|>1400$ at $t=200$), and the Yoshida triple jump with step $0.01$ as reference; $\rho=0.3$ lies below the value $\frac12$ required by the Dollard modifier (section~\ref{sec:Slrresults}). Figure~\ref{fig:slow} shows that the error saturates in time and scales with the order of the method, as predicted by Theorem~\ref{thm:conv}(a), which holds for every $\rho>0$. At $\Delta=0.2$ and $t=200$ the Lie errors are $0.102$, $0.118$ and $0.114$, against the limits $\frac\Delta2\|V\psi\|=0.101$, $0.119$ and $0.134$, which are approached slowly for $\rho<1$. For $\rho=0.7$ the Yoshida errors at $\Delta=0.1$ and $0.05$ are dominated by the aliased channels: since $\hat V(\zeta)$ decays only like $e^{-|\zeta|}$, the aliased part of the Strang error decays like $e^{-c/\sqrt\Delta}$, with $c$ between $2.4$ and $2.8$, close to $(2\pi)^{1/2}$.

The last panel shows the error $\int\big||\widehat{U^n\psi}|^2-|\hat\psi_{\rm ref}|^2\big|\frac{d\xi}{2\pi}$ of the momentum distribution at $t=200$, which is at most twice the error in norm. Its values at $t=100$ and $t=200$ agree to at least three digits, and to $2\%$ for Lie splitting at $\rho=0.3$. It is of order $\Delta^2$ for Strang splitting. For Lie splitting the local slopes are $1.07$, $1.03$, $1.02$ for $\rho=1.5$, and $1.32$, $1.12$, $1.02$ and $1.67$, $1.45$, $1.25$ for $\rho=0.7$ and $0.3$: the order $\Delta$ is approached slowly for $\rho<1$.

\begin{figure}[ht]
\centering
\includegraphics[width=0.62\textwidth]{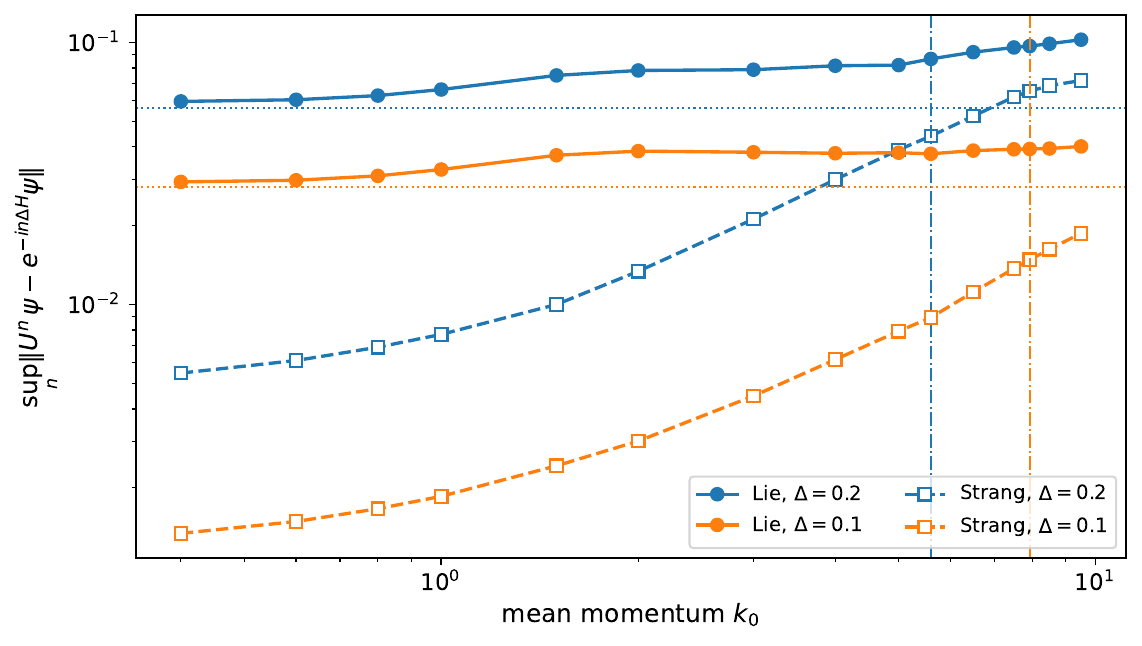}
\caption{Maximal error over time against the mean momentum $k_0$, for $\Delta=0.2$ (blue) and $0.1$ (orange). Dotted: $\frac\Delta2\|V\psi\|$. Dash-dotted: $k_0=\sqrt{2\pi/\Delta}$.}
\label{fig:energy}
\end{figure}

\paragraph{Energy dependence.} Figure~\ref{fig:energy} shows the maximal error over time for packets $\psi=c\,e^{-x^2/50+ik_0x}$ on top of $V=1.5\,e^{-x^2}$, with momentum spread $0.14$ and $0.4\le k_0\le9.5$. The time interval is $[0,\min(150,40/k_0+5)]$ and the reference step $0.1/\max(40,\lceil k_0^2\rceil)$. The Lie error stays close to $\frac\Delta2\|V\psi\|$, which does not depend on $k_0$, until the kinetic part of the error becomes comparable. The Strang error grows roughly linearly in $k_0$ while $\Delta k_0^2\lesssim1$, as expected from a defect of size $\Delta^3k_0^2$ per step (Lemma~\ref{lem:defect}(b)) during a time of order $1/k_0$ spent near the potential. The constants of Theorem~\ref{thm:conv} thus depend on the energy window $\supp g$. Nothing special happens at $k_0=\sqrt{2\pi/\Delta}$, where the energy is aliased onto the threshold, because the kicks couple $k_0$ to momenta near $0$ only through $\hat V$ at a momentum transfer of order $k_0$, which is exponentially small for this potential; this is the channel $m=-1$ of the aliased channels discussed above. No growth of the error is visible at low energies either, since slow packets are pushed out of the interaction region by this repulsive potential; the degeneration of the constants as $\supp g$ approaches $0$ is a limitation of the estimates that is not observed here.

%=====================================================================
\section{Bound states and ionization}\label{sec:Sbound}
%=====================================================================

This section considers a bound state $H\phi_0=E_0\phi_0$, $E_0<0$, $\|\phi_0\|=1$, and its survival amplitude $\ip{\phi_0}{U^n\phi_0}$ under the scheme. Theorem~\ref{thm:meta} shows that the survival probability stays close to $1$ at least up to times $\Delta^{-N}$ for every $N$, and Theorem~\ref{thm:noeig} that the bound state is eventually ionized when $V$ has compact support or decays like $e^{-a|x|^\beta}$ with $\beta>2$. Section~\ref{sec:Smechanism} describes the spectral meaning of ionization and the mechanism that drives it. Section~\ref{sec:Sexpmeta} proves the exponential metastability of Theorem~\ref{thm:expmeta}, section~\ref{sec:Snoeig} extends the absence of eigenvalues to general compositions under conditions (a) and (b) of Theorem~\ref{thm:noeig}, and section~\ref{sec:Srate} bounds the rate of ionization from below under a non-cancellation hypothesis, which we do not prove and which is supported by the numerical evidence of section~\ref{sec:Sevidence}. For ionization by time-periodic perturbations and the Fermi golden rule, we refer to \cite{S:SW98,S:CCLR01} and \cite[Ch.~XII]{S:RS78}.

\subsection{Ionization and its mechanism}\label{sec:Smechanism}

\paragraph{Ionization and the other bound states.} The definition of ionization in section~\ref{sec:bound-results} involves only $\phi_0$, but it controls the overlap of $U^n\phi_0$ with every vector. The survival amplitude $\ip{\phi_0}{U^n\phi_0}=\int e^{in\theta}d\mu_{\phi_0}(\theta)$ is a Fourier coefficient of the spectral measure $\mu_{\phi_0}$ of $\phi_0$ for $U$. By Wiener's theorem \cite[Ch.~I]{S:Kat04}, see also \cite{S:RS79}, its decay in the mean is equivalent to $\mu_{\phi_0}$ having no atoms, i.e., to $\phi_0$ being orthogonal to every eigenvector of $U$, and hence, since $U$ has no singular continuous spectrum, to $\phi_0\in\HH_{\rm ac}(U)$. For every $f$, the measure $\mu_{f,\phi_0}$ with $\ip{f}{U^n\phi_0}=\int e^{in\theta}d\mu_{f,\phi_0}$ is then absolutely continuous, and $\ip{f}{U^n\phi_0}\to0$ by the Riemann--Lebesgue lemma. In particular, the overlap with every other bound state $\phi_j$ of $H$ tends to $0$. A permanent transfer of mass from $\phi_0$ to $\phi_j$ would require an eigenvector $w$ of $U$ with $\ip{w}{\phi_0}\ne0$ and $\ip{w}{\phi_j}\ne0$, and then $\mu_{\phi_0}$ would have an atom of mass $|\ip{w}{\phi_0}|^2$. %Transfers over finite times are not excluded.

\paragraph{The quasi-energy.} Theorem~\ref{thm:meta} identifies the phase of the survival amplitude. For Lie and Strang splitting, $\ip{\phi_0}{U^n\phi_0}\approx e^{-in\Delta E(\Delta)}$ with $E(\Delta)=E_0+e_2\Delta^2+O(\Delta^4)$ and $e_2$ given by \eqref{eq:quasiE}, up to times $\Delta^{-N}$. It is the eigenvalue of the modified Hamiltonian of the Strang step, given by the Baker--Campbell--Hausdorff formula as in backward error analysis \cite[Ch.~IX]{S:HLW06}, obtained by Rayleigh--Schr\"odinger perturbation theory \cite[Ch.~XII]{S:RS78}.

\paragraph{Ionization by time-periodic fields.} A bound state of a Schr\"odinger operator subject to a time-periodic perturbation is generically ionized, in the sense that its survival amplitude tends to $0$ and its mass disperses. Complete ionization was proved for model systems with arbitrary amplitude and frequency in \cite{S:CCLR01}. For small, localized, time-periodic perturbations, metastability and decay at the rate given by the Fermi golden rule were proved in \cite{S:SW98}. In the dispersive PDE literature, the same mechanism is called radiation damping \cite{S:SW99}. In the setting of the paper, $U$ is the monodromy of the time-periodic Hamiltonian $H_0+\Delta V\sum_k\delta(t-k\Delta)$. Our setting differs from those works in two respects. The harmonics of the perturbation are not small, and the relevant limit is that of large frequency $2\pi/\Delta$. The rates are then exponentially small in $\Delta^{-1/2}$ for analytic potentials, since the smallness comes from the matrix elements at large momentum.

\paragraph{The golden rule.} By the Poisson summation formula \cite{S:Kat04}, $\Delta\sum_k\delta(t-k\Delta)=\sum_{m\in\Z}e^{2\pi imt/\Delta}$ as distributions, so that $U$ is the monodromy operator, i.e., the propagator over one period \cite{S:Yaj77,S:How79}, of the time-periodic Hamiltonian $H+V\sum_{m\ne0}e^{2\pi imt/\Delta}$. In the interaction picture with respect to $H$, in which the evolution generated by $H$ is factored out and only the time-periodic part drives transitions, the harmonic $e^{2\pi imt/\Delta}V$ couples $\phi_0$ to the generalized eigenfunctions of $H$ at energy $E_0+2\pi m/\Delta$, which lies in the continuum when it is positive. First-order time-dependent perturbation theory then gives the Fermi golden rule \cite[Ch.~XII]{S:RS78}, \cite{S:SW98}: the survival probability decays exponentially, at a rate proportional to the squared matrix elements of the perturbation between $\phi_0$ and the continuum at the resonant energies, here
\begin{equation}\label{eq:Sfgr}
\Gamma_{\rm FGR}(\Delta)=\sum_{m\ge1}\Gamma_m,\qquad\Gamma_m=2\pi\,\big\langle V\phi_0,\ \delta\big(H-E_0-2\pi m/\Delta\big)\,V\phi_0\big\rangle,
\end{equation}
a sum of the rates $\Gamma_m$ of the individual channels. To leading order in $V$, the generalized eigenfunctions may be replaced by plane waves, and since $V\phi_0=(E_0-H_0)\phi_0$, $\ip{e^{i\xi\cdot x}}{V\phi_0}=(E_0-|\xi|^2)\hat\phi_0(\xi)$ involves $\hat\phi_0$ on the spheres $|\xi|^2=E_0+2\pi m/\Delta$, of radius of order $\Delta^{-1/2}$. For a potential analytic in a strip, $\phi_0$ extends analytically to a strip and $\hat\phi_0$ decays exponentially, so that $\Gamma_{\rm FGR}(\Delta)$ is exponentially small in $\Delta^{-1/2}$.

\paragraph{Higher-order compositions.} In the interaction picture with respect to $H_0$, the composition \eqref{eq:comp} sees the $\Delta$-periodic kick comb $\kappa(\tau)=\Delta\sum_k\sum_ja_j\delta(\tau-\beta_j\Delta-k\Delta)$, with the kick times $\beta_j\Delta$ of section~\ref{sec:setting}. Its Fourier coefficients are $\hat\kappa_m=\sum_ja_je^{-2\pi im\beta_j}$, and the golden rule becomes
\begin{equation}\label{eq:fgr-comp}
\Gamma_U\approx\sum_{m\ge1}|\hat\kappa_m|^2\,\Gamma_m .
\end{equation}
For Lie and Strang splitting, $\hat\kappa_m=1$. For the Yoshida triple jump of section~\ref{sec:Snumerics}, $\hat\kappa_m=w_1+(1-w_1)\cos(2\pi mw_1)$, and $|\hat\kappa_1|^2=2.43$.

\subsection{Exponential metastability}\label{sec:Sexpmeta}
The counterpart of ionization is the survival of the bound state over long times. For potentials whose Fourier transform decays exponentially, i.e., under hypothesis \eqref{ass:A}, this survival time is at least exponentially long in $\Delta^{-1/2}$. We recall the result and give its proof, whose quasi-mode $\tilde\psi$ is used again in sections~\ref{sec:Srate} and~\ref{sec:Sevidence}.

\begin{theorem}[Theorem~\ref{thm:expmeta}]\label{thm:Sexpmeta}
Assume \eqref{ass:V} and \eqref{ass:A}, let $E_0<0$ be a simple eigenvalue of $H$ with normalized eigenfunction $\phi_0$, and let $U$ be the Lie, reversed Lie or Strang step. There are $c,C,\Delta_0>0$ and, for $\Delta\in(0,\Delta_0]$, real numbers $\tilde E(\Delta)$ with $\tilde E(\Delta)-E_N(\Delta)=O(\Delta^{N+1})$ for every $N$, where $E_N$ is as in Theorem~\ref{thm:meta}, such that
\begin{equation}\label{eq:Sexpmeta}
\big|\ip{\phi_0}{U^n\phi_0}-e^{-in\Delta\tilde E(\Delta)}\big|\le C\big(\Delta^2+|n|e^{-c/\sqrt\Delta}\big)\qquad(n\in\Z,\ \Delta\in(0,\Delta_0]).
\end{equation}
\end{theorem}

\begin{proof}
We treat Strang splitting; Lie and reversed Lie splitting follow as in Step 5 of the proof of Theorem~\ref{thm:meta}. Write $a=V/2$, $\Pi_K=\one_{|p|\le K}$ and $L^2_K=\Ran\Pi_K$.

\emph{Step 0: momentum weights.} Fix $\delta\in(0,\pi/4)$. For $0\le\varepsilon\le\alpha/2$, let $V^{(\varepsilon)}$ be the operator with kernel $(2\pi)^{-d}\hat V(\xi-\eta)e^{\varepsilon(\jap{\xi}-\jap{\eta})}$ in Fourier variables, so that $Ve^{-\varepsilon\jap{p}}=e^{-\varepsilon\jap{p}}V^{(\varepsilon)}$. Since $|\jap{\xi}-\jap{\eta}|\le|\xi-\eta|$, this kernel is dominated by the convolution kernel $(2\pi)^{-d}|\hat V(\xi-\eta)|e^{\varepsilon|\xi-\eta|}$, and Schur's test gives $\|V^{(\varepsilon)}\|\le M_\alpha$ and $\|V^{(\varepsilon)}-V\|\le C\varepsilon$. Let $d_0>0$ be the distance from $E_0$ to the rest of the spectrum of $H$, and $P_0$ the projection onto $\phi_0$. For $|z-E_0|=d_0/2$ and $\varepsilon$ small, $z-H_0-V^{(\varepsilon)}=(z-H)-(V^{(\varepsilon)}-V)$ is invertible with uniformly bounded inverse, and applying $z-H$ to $e^{-\varepsilon\jap{p}}(z-H_0-V^{(\varepsilon)})^{-1}g$ gives $(z-H)^{-1}e^{-\varepsilon\jap{p}}=e^{-\varepsilon\jap{p}}(z-H_0-V^{(\varepsilon)})^{-1}$. Integrating over the circle, $P_0e^{-\varepsilon\jap{p}}=e^{-\varepsilon\jap{p}}P_0^{(\varepsilon)}$ with $P_0^{(\varepsilon)}$ bounded, and applying this to $\phi_0$ gives $\phi_0\in\mathrm{Dom}(e^{\varepsilon\jap{p}})$, a Combes--Thomas type bound \cite{S:CT73} in momentum space. We now fix $\varepsilon\le\alpha/2$ so small that $\phi_0\in\mathrm{Dom}(e^{2\varepsilon\jap{p}})$ and that the conditions on $\varepsilon$ in Steps 4 and 6 hold. Below, constants depend on $\delta$ and $\varepsilon$, but not on $\Delta$.

\emph{Step 1: compressed scheme.} Let $K=(\delta/\Delta)^{1/2}$, and let $\Delta$ be so small that $\Delta\|V\|_\infty\le\delta$. On $L^2_K$, let $a_K=\Pi_Ka\Pi_K$, $H_K=\Pi_KH\Pi_K=H_0\Pi_K+2a_K$, and $U_K=e^{-i\Delta a_K}e^{-i\Delta H_0}e^{-i\Delta a_K}$. Then $U_K$ is unitary, $\|U_K-e^{-i\Delta H_0}\Pi_K\|\le\delta$, and $\sigma(\Delta H_K)\subset[-\delta,2\delta]$. Hence the spectra of $U_K$ and $e^{-i\Delta H_K}$ lie in an arc of length less than $\pi$ that avoids $-1$, and $\tilde H_K=\frac i\Delta\mathrm{Log}\,U_K$ (principal branch) is a bounded self-adjoint operator on $L^2_K$ with $e^{-i\Delta\tilde H_K}=U_K$.

\emph{Step 2: $\tilde H_K$ is close to $H_K$.} The kernel of $[H_0,a_K]$ in Fourier variables is $(|\xi|^2-|\eta|^2)(2\pi)^{-d}\hat a(\xi-\eta)$, and $||\xi|^2-|\eta|^2|\le2K|\xi-\eta|$ on $L^2_K$. Schur's test therefore gives $\|[H_0\Pi_K,[H_0\Pi_K,a_K]]\|\le CK^2$ and $\|[a_K,[a_K,H_0\Pi_K]]\|\le CK$. The Duhamel argument of Lemma~\ref{lem:defect}(b) applies to the bounded operators $H_0\Pi_K$ and $a_K$, except that the conjugations $e^{-isa_K}H_0\Pi_Ke^{isa_K}$ and $e^{-isa_K}[H_0\Pi_K,a_K]e^{isa_K}$ no longer terminate; their Taylor expansions to second and first order have remainders bounded by $s^2\|[a_K,[a_K,H_0\Pi_K]]\|$ and $s\|[a_K,[a_K,H_0\Pi_K]]\|$. With the bounds above, $\|D_S(s)\|\le Cs^2(K+K^2)$, which gives $\|U_K-e^{-i\Delta H_K}\|\le C\Delta^3K^2=C\delta\Delta^2$. For unitary operators with spectra in a fixed arc avoiding $-1$, the Riesz--Dunford integral shows that $\mathrm{Log}$ is Lipschitz. Since $\mathrm{Log}\,e^{-i\Delta H_K}=-i\Delta H_K$, we get $\|\tilde H_K-H_K\|\le C\delta\Delta$.

\emph{Step 3: an isolated eigenvalue.} As an operator on $L^2_K$, $\Pi_KV\Pi_K$ is compact (Fact~\ref{fact:compact}, whose proof applies to $b=\one_{|\xi|\le K}$), so $\sigma_{\rm ess}(H_K)=[0,K^2]$. The quadratic form of $H_K$ is the restriction of that of $H$, so by min-max the eigenvalues of $H_K$ below $0$ are at least those of $H$. Using the test vectors $\Pi_K\phi$, with $\phi$ eigenfunctions of $H$ at energies $\le E_0$, for which $(1-\Pi_K)\phi\to0$ in $\mathcal H^1$, they exceed them by $o(1)$ as $K\to\infty$. Hence $H_K$ has exactly one eigenvalue in $(E_0-d_0/2,E_0+d_0/2)$, which is simple. By Step 2, for $\Delta$ small, $\tilde H_K$ has exactly one eigenvalue $\tilde E$ in $(E_0-d_0/4,E_0+d_0/4)$, which is simple. Let $\tilde\psi$ be a normalized eigenvector with $\ip{\phi_0}{\tilde\psi}\ge0$. Since $\|(\tilde H_K-E_0)\Pi_K\phi_0\|\le C(\Delta+e^{-\varepsilon K})$, the gap gives $\|\tilde\psi-\phi_0\|\le C(\Delta+e^{-\varepsilon K})$.

\emph{Step 4: momentum decay of $\tilde\psi$.} The weight $\jap{p}$ commutes with $\Pi_K$ and $H_0$, so $e^{\varepsilon\jap{p}}U_Ke^{-\varepsilon\jap{p}}=e^{-i\Delta a_K^{(\varepsilon)}}e^{-i\Delta H_0}e^{-i\Delta a_K^{(\varepsilon)}}$ with $a_K^{(\varepsilon)}=\Pi_Ka^{(\varepsilon)}\Pi_K$, and $\|e^{\varepsilon\jap{p}}U_Ke^{-\varepsilon\jap{p}}-U_K\|\le C\varepsilon\Delta$. On $L^2_K$, $e^{\varepsilon\jap{p}}$ is bounded with bounded inverse, so $e^{\varepsilon\jap{p}}\tilde H_Ke^{-\varepsilon\jap{p}}=\frac i\Delta\mathrm{Log}(e^{\varepsilon\jap{p}}U_Ke^{-\varepsilon\jap{p}})$, and the Riesz--Dunford integral gives $\|e^{\varepsilon\jap{p}}\tilde H_Ke^{-\varepsilon\jap{p}}-\tilde H_K\|\le C\varepsilon$, so that for $\varepsilon$ small the Riesz projection $\tilde P$ of $\tilde H_K$ onto $\tilde\psi$ satisfies $\|e^{\varepsilon\jap{p}}\tilde Pe^{-\varepsilon\jap{p}}\|\le C$. Since $\tilde P\Pi_K\phi_0=\ip{\tilde\psi}{\Pi_K\phi_0}\tilde\psi$ with $\ip{\tilde\psi}{\Pi_K\phi_0}\ge\frac12$, we get $\|e^{\varepsilon\jap{p}}\tilde\psi\|\le2C\|e^{\varepsilon\jap{p}}\phi_0\|$, uniformly in $\Delta$.

\emph{Step 5: the defect is exponentially small.} Since $U_K\tilde\psi=e^{-i\Delta\tilde E}\tilde\psi$, it suffices to bound $(U-U_K)\tilde\psi$. For $f\in L^2_K$, $(a-a_K)f=(1-\Pi_K)af$. Hence, by Duhamel,
\[
\|(e^{-i\Delta a}-e^{-i\Delta a_K})f\|\le\Delta\sup_{0\le\sigma\le\Delta}\|(1-\Pi_K)e^{-\varepsilon\jap{p}}\|\,\|a^{(\varepsilon)}\|\,\|e^{\varepsilon\jap{p}}e^{-i\sigma a_K}f\|\le C\Delta e^{-\varepsilon K}\|e^{\varepsilon\jap{p}}f\| .
\]
Here we used $\|e^{\varepsilon\jap{p}}e^{-i\sigma a_K}e^{-\varepsilon\jap{p}}\|\le e^{\sigma\|a^{(\varepsilon)}\|}$. Applying this to $f=\tilde\psi$ and to $f=e^{-i\Delta H_0}e^{-i\Delta a_K}\tilde\psi$, whose weighted norms are bounded by Step 4, gives $\|(U-e^{-i\Delta\tilde E})\tilde\psi\|\le C\Delta e^{-\varepsilon K}$. The bound of Step 5 of the proof of Theorem~\ref{thm:meta}, with $f=\phi_0$, $\psi=\tilde\psi$ and $|\beta|\le\|\phi_0-\tilde\psi\|\le C(\Delta+e^{-\varepsilon K})$ by Step 3, gives \eqref{eq:Sexpmeta} with $c=\varepsilon\sqrt\delta$.

\emph{Step 6: $\tilde E$ and $E_N$.} The quasi-modes $\chi^N$ of Theorem~\ref{thm:meta} have exponential momentum decay. Indeed, $R_0=(H-E_0+P_0)^{-1}-P_0$. Since $\phi_0\in\mathrm{Dom}(e^{2\varepsilon\jap{p}})$, $e^{\varepsilon\jap{p}}P_0e^{-\varepsilon\jap{p}}$ is bounded and close to $P_0$ for $\varepsilon$ small, and the argument of Step 0, applied to $H-E_0+P_0$, shows that $e^{\varepsilon\jap{p}}R_0e^{-\varepsilon\jap{p}}$ is bounded. Each $h_j$ is a finite sum of products of factors $H_0$ and $V$, and $e^{\varepsilon'\jap{p}}H_0e^{-\varepsilon''\jap{p}}$, $e^{\varepsilon'\jap{p}}Ve^{-\varepsilon'\jap{p}}$ are bounded for $\varepsilon'<\varepsilon''$. By induction on $k$, $\chi_k\in\mathrm{Dom}(e^{\varepsilon'\jap{p}})$ for every $\varepsilon'<\varepsilon$. Fix such an $\varepsilon'$. Then $\|(1-\Pi_K)\chi^N\|\le C_Ne^{-\varepsilon'K}$, and the argument of Step~5, with $\varepsilon'$ in place of $\varepsilon$ and $f=\Pi_K\chi^N$, gives $\|(U-U_K)\Pi_K\chi^N\|\le C_N\Delta e^{-\varepsilon'K}$. With $\lambda_N=e^{-i\Delta E_N}$ and Step 4 of the proof of Theorem~\ref{thm:meta},
\[
\|(U_K-\lambda_N)\Pi_K\chi^N\|\le\|(U_K-U)\Pi_K\chi^N\|+\|(U-\lambda_N)\chi^N\|+2\|(1-\Pi_K)\chi^N\|\le C_N\Delta^{N+2},
\]
since $e^{-\varepsilon'K}=O(\Delta^\infty)$, and $\|\Pi_K\chi^N\|\ge\frac12$ for $\Delta$ small. By Step 3, there is $c'>0$ such that the spectrum of $U_K$ within distance $c'\Delta$ of $\lambda_N$ consists of $e^{-i\Delta\tilde E}$ alone. Since $U_K$ is unitary, $|e^{-i\Delta\tilde E}-\lambda_N|\le2C_N\Delta^{N+2}$, i.e., $\tilde E-E_N=O(\Delta^{N+1})$.
\end{proof}

\subsection{Absence of eigenvalues for general compositions}\label{sec:Snoeig}
We extend Theorem~\ref{thm:noeig}(a,b) to compositions whose kicks act at distinct times modulo the step.

For the composition \eqref{eq:comp}, let $\beta_j\Delta$, $\beta_j=b_1+\dots+b_{j-1}$ ($\beta_1=0$), be the kick times of section~\ref{sec:setting}. With $\tilde K_j=e^{i\beta_j\Delta H_0}K_je^{-i\beta_j\Delta H_0}$ one checks that $U=U_0\tilde K_s\cdots\tilde K_1$.

\begin{proposition}\label{prop:Snoeig}
Let $V\in L^\infty(\R^d;\R)$ satisfy condition (a) or (b) of Theorem~\ref{thm:noeig}, $\Delta>0$, and let $U$ be a composition whose kick times $\beta_1\Delta,\dots,\beta_s\Delta$ are pairwise distinct modulo $\Delta$, i.e., the $\beta_j$ are pairwise distinct modulo $1$. Then $U$ has no eigenvalues. The same holds for $PUP^{-1}$ with $P$ a processor.
\end{proposition}

The last statement covers the Strang step (kick times $0$ and $\Delta$, conjugate to Lie splitting) and the triple jump $U_{\rm Y}$ of section~\ref{sec:Snumerics}: conjugating by its last kick merges the kicks at times $0$ and $\Delta$, and the remaining kick times $0,w_1\Delta,(1-w_1)\Delta$ are distinct modulo $\Delta$ since $2w_1\notin\Z$.

\begin{proof}
Let $Uw=\zeta w$, $\zeta=e^{-i\Delta E}$, and $\nu_j=K_j-1$. Since $\tilde K_s\cdots\tilde K_1-1=\sum_j(\tilde K_j-1)\tilde K_{j-1}\cdots\tilde K_1$, the equation $\tilde K_s\cdots\tilde K_1w=\zeta U_0^{-1}w$ becomes, in Fourier variables,
\[
\big(\zeta e^{i\Delta|\xi|^2}-1\big)\hat w=\sum_{j=1}^se^{i\beta_j\Delta|\xi|^2}F^{(j)},\qquad F^{(j)}=\widehat{\nu_jg_j},\quad g_j=e^{-i\beta_j\Delta H_0}\tilde K_{j-1}\cdots\tilde K_1w .
\]
Since $\|g_j\|=\|w\|$, Step~2 of the proof of Theorem~\ref{thm:noeig} applies to each $F^{(j)}$, and Step~3 shows that the right-hand side vanishes on the spheres $|\xi|=\rho_m$, $\rho_m^2=E+2\pi m/\Delta$, $m\ge m_0$.

Fix $\theta\in S^{d-1}$ and write $F^{(j)}(t\theta)=A^+_j(t^2)+tA^-_j(t^2)$, with $A^\pm_j$ entire. In case (a), $|F^{(j)}(t\theta)|\le Ce^{R|t|}$, so $A^\pm_j$ are of order at most $\frac12$; in case (b) they are of order less than $1$, by the growth bound of Step~2 of the proof of Theorem~\ref{thm:noeig}. Vanishing at $t=\pm\rho_m$ gives $\sum_je^{i\beta_j\Delta\rho_m^2}A^\pm_j(\rho_m^2)=0$. Now $e^{i\beta_j\Delta\rho_m^2}=e^{i\beta_j\Delta E}e^{2\pi i\beta_jm}$. The points $\beta_j$ modulo $1$ lie in an arc of length $<1$, so there are $c\in\R$ and integers $n_j$ with $\sigma_j=\beta_j-c-n_j\in(-\frac12,\frac12)$, pairwise distinct, so that $e^{2\pi i\beta_jm}=e^{2\pi icm}e^{2\pi i\sigma_jm}$. Hence the entire functions
\[
\Upsilon^\pm(\mu)=\sum_jc_je^{2\pi i\sigma_j\mu}A^\pm_j\big(E+2\pi\mu/\Delta\big),\qquad c_j=e^{i\beta_j\Delta E}\ne0,
\]
vanish at all integers $\mu\ge m_0$. They are of exponential type, and $|\Upsilon^\pm(iy)|\le C_\epsilon e^{(2\pi\max_j|\sigma_j|+\epsilon)|y|}$ with $2\pi\max_j|\sigma_j|<\pi$. By Carlson's theorem \cite{S:Boa54}, applied to $\Upsilon^\pm(\mu+m_0)$, $\Upsilon^\pm\equiv0$.

It remains to show that $\sum_je^{2\pi i\sigma_j\mu}B_j(\mu)\equiv0$, with distinct real $\sigma_j$ and entire $B_j$ of order $<1$, forces $B_j\equiv0$. If not, let $\sigma_k$ be the largest $\sigma_j$ with $B_j\not\equiv0$. For $\mu=-iy$, $y\to\infty$,
\[
|B_k(-iy)|\le\sum_{\sigma_j<\sigma_k}e^{-2\pi(\sigma_k-\sigma_j)y}|B_j(-iy)|\le Ce^{-cy}
\]
for some $c>0$, since the $B_j$ have order $<1$, so that $|B_j(-iy)|\le C_\epsilon e^{\epsilon y}$. But an entire function $B_k\not\equiv0$ of order $\rho_B<1$ satisfies $\min_{|\mu|=r}|B_k(\mu)|\ge M(r)^{\cos\pi\rho_B-\epsilon}\ge e^{-r^{\rho_B+2\epsilon}}$ for a sequence $r\to\infty$, where $M(r)=\max_{|\mu|=r}|B_k|$, by the $\cos\pi\rho$ theorem \cite{S:Lev96}; this contradicts the exponential decay along $\mu=-iy$. Hence $A^\pm_j\equiv0$, so $F^{(j)}=0$ on every line through the origin, the right-hand side vanishes, and $\hat w=0$ as in Step~5.

For $PUP^{-1}$, $P^{-1}$ maps eigenvectors of $PUP^{-1}$ to eigenvectors of $U$.
\end{proof}

\begin{remark}[Other decay classes]
We have not examined case (c) of Theorem~\ref{thm:noeig} for general compositions. For compactly supported $V$ of Gevrey class $s>1$, $|\hat V(\zeta)|\le Ce^{-c|\zeta|^{1/s}}$, and the proof of Theorem~\ref{thm:Sexpmeta} should extend with the weights $e^{\varepsilon\jap{\xi}^{1/s}}$, giving survival up to times $e^{c\Delta^{-1/(2s)}}$ before the ionization of Corollary~\ref{cor:ion}. For exponentially decaying $V$, such as the P\"oschl--Teller well of section~\ref{sec:Sevidence}, $F$ is analytic only in a strip, where bounded functions can vanish at the points $\pm\rho_m$, and we have no proof that $\sigma_p(U)=\emptyset$. The rates of section~\ref{sec:Sevidence} are, however, unchanged to six digits when the well is smoothly truncated outside $|x|\le28$, an $O(e^{-56})$ change, and Theorem~\ref{thm:noeig}(a) applies to the truncated well.
\end{remark}

\subsection{A lower bound on the rate of ionization}\label{sec:Srate}

The vector $\tilde\psi$ of the proof of Theorem~\ref{thm:Sexpmeta} is an exact eigenvector of the scheme compressed to momenta $|\xi|\le K$, and the full scheme couples it to the continuum only through the residual $(U-e^{-i\Delta\tilde E})\tilde\psi$. The resonant spheres lie beyond $K$, and are reached only through the exponentially decaying tails of $\hat V$ and $\hat\phi_0$. The absence of eigenvalues provides no rate of ionization. A rate requires an upper bound on the spectral density of $\phi_0$ for $U$ near its quasi-energy, which we derive from a lower bound on a golden rule for the residual, under a non-cancellation hypothesis. That golden rule is approximated by \eqref{eq:Sfgr}, but the two differ by a factor that tends to about $1.13$ as $\Delta\to0$ (section~\ref{sec:Sevidence}).

\medskip\noindent{\bf Setting.} Let $U$ be the Lie, reversed Lie or Strang step, $E_0<0$ a simple eigenvalue of $H$ with eigenfunction $\phi_0$, and $\tilde\psi$ a normalized vector with $\|\tilde\psi-\phi_0\|\le C\Delta$ and $0<\varrho(\Delta)=\min_{\lambda'\in\T}\|(U-\lambda')\tilde\psi\|<\Delta$. Examples are the normalized quasi-modes $\chi^N$ of Theorem~\ref{thm:meta}, with $\varrho=O(\Delta^{N+2})$, and, under \eqref{ass:A}, the quasi-mode of the proof of Theorem~\ref{thm:Sexpmeta}, with $\varrho=O(\Delta e^{-c/\sqrt\Delta})$. With $a=\ip{\tilde\psi}{U\tilde\psi}$, the minimum $\varrho=(2-2|a|)^{1/2}$ is attained at $\lambda_0=a/|a|$. We put
\[
\lambda=\lambda_0e^{i\varrho},\quad r=(U-\lambda)\tilde\psi,\quad y=\frac{U^*r}{\|r\|},\quad \eta=\ip{y}{\tilde\psi},\quad \tau=1-\frac{\|r\|}\eta,\quad \varsigma=\bar\tau-1,
\]
and $U_1=U(1+(\tau-1)\ip{y}{\cdot}y)$. Let $P$ be the orthogonal projection onto $\tilde\psi$ and $Q=1-P$. Step 1 of the proof shows that $U_1$ is unitary with $U_1\tilde\psi=\lambda\tilde\psi$, so that $Q\HH$ reduces $U_1$ and $\hat U=U_1|_{Q\HH}$ is unitary. We put $q=QU^*r$ and $\mathcal M(z)=\ip{Qy}{\hat U(\hat U-z)^{-1}Qy}$ for $|z|<1$. For a unitary $T$ and a vector $f$, $P^T_f(z)=(1-|z|^2)\|(T-z)^{-1}f\|^2=\int_\T\frac{1-|z|^2}{|e^{i\theta}-z|^2}d\mu^T_f(\theta)$ is the Poisson integral of the spectral measure of $f$. Finally, $\Theta_\Delta=\{\zeta\in\T:|\zeta-\lambda_0|<\varrho/\Delta\}$, an arc containing $\lambda$, and $\mathcal Z_\Delta=\{t\zeta:\zeta\in\Theta_\Delta,\ 1-\varrho/\Delta<t<1\}$.

\paragraph{Hypothesis (NC$_\Delta$).} There is $\gamma_0(\Delta)\in(0,1]$ such that $P^{\hat U}_q\ge\gamma_0(\Delta)$ and $|\varsigma\mathcal M|\le\frac12$ on $\mathcal Z_\Delta$.

\begin{theorem}\label{thm:Srate}
Assume \eqref{ass:V} and (NC$_\Delta$) for $\Delta\in(0,\Delta_0]$. There is $C$ such that for these $\Delta$ and all $N\ge1$:
\\[1mm]
\textup{(a)} $\|P_{\rm pp}(U)\phi_0\|\le C\Delta$;
\\[1mm]
\textup{(b)} $\lim_N\frac1N\sum_{n<N}|\ip{\phi_0}{U^n\phi_0}|^2\le C\Delta^4$;
\\[1mm]
\textup{(c)} $\frac1N\sum_{n<N}|\ip{\phi_0}{U^n\phi_0}|^2\le C/(N\gamma_0(\Delta))+C\Delta^2$.
\end{theorem}

When $\sigma_p(U)=\emptyset$, as in Corollary~\ref{cor:ion}, parts (a) and (b) are superseded, and the content of the theorem is (c), which provides a lower bound on the rate of ionization by the golden rule for the residual, up to a universal factor. Upper bounds on the rate follow from the metastability results of Theorems~\ref{thm:meta} and~\ref{thm:expmeta}. Under \eqref{ass:A}, both are exponential in $\Delta^{-1/2}$, and even sharp exponents would differ by a factor of about $2$, because the quasi-mode bounds are linear in the residual while the decay rate is quadratic in its amplitude on the resonant spheres. Closing the gap would require upper bounds and H\"older regularity of $\mu^{\hat U}_q$ on the scale of the resonance width.

\begin{proof}
We write $\psi=\tilde\psi$.

\emph{Step 1: rank-one decoupling.} Let $w=\bar\lambda a$, so $\|r\|^2=2-2\,{\rm Re}\,w$. Since $\ip{r}{\psi}=\bar a-\bar\lambda$, $\|r\|\eta=\ip{r}{U\psi}=\ip{r}{\lambda\psi+r}=\bar w-1+\|r\|^2=1-w$, and $w\ne1$ since $r\ne0$. Hence $\eta\ne0$, and $\tau=(1-w-\|r\|^2)/(1-w)=(\bar w-1)/(1-w)$ has modulus $1$, so $1+(\tau-1)\ip{y}{\cdot}y$ is unitary with inverse $1+\varsigma\ip{y}{\cdot}y$. Since $Uy=r/\|r\|$, $U_1\psi=\lambda\psi+r(1+(\tau-1)\eta/\|r\|)=\lambda\psi$, and $U_1^*\psi=\bar\lambda\psi$ by unitarity. Moreover $U=U_1(1+\varsigma\ip{y}{\cdot}y)$ and $|\varsigma||\eta|=|\tau-1||\eta|=\|r\|$.

\emph{Step 2: a Poisson bound.} Let $z\in\mathcal Z_\Delta$, $v=(U-z)^{-1}\psi$, $\alpha=\ip{\psi}{v}$ and $c=\varsigma\ip{y}{v}$. The equation $(U-z)v=\psi$ reads $U_1v+cU_1y-zv=\psi$, and its $Q$-component is $(\hat U-z)Qv=-c\,\hat UQy$. Hence $\|Qv\|=|c|\|(\hat U-z)^{-1}Qy\|$. Since $Py=\bar\eta\psi$, $\ip{y}{v}=\eta\alpha-c\mathcal M(z)$, so $c(1+\varsigma\mathcal M(z))=\varsigma\eta\alpha$ and $|c|\ge\frac23\|r\||\alpha|$ by (NC$_\Delta$). Therefore $P^U_\psi(z)\ge(1-|z|^2)\|Qv\|^2\ge|\alpha|^2W$ with $W=\frac49P^{\hat U}_q(z)$. On the other hand $P^U_\psi(z)={\rm Re}\,\ip{\psi}{(U+z)(U-z)^{-1}\psi}=1+2\,{\rm Re}(z\alpha)\le1+2(P^U_\psi(z)/W)^{1/2}$, and solving for $P^U_\psi(z)^{1/2}$ gives $P^U_\psi(z)\le2+4/W\le B_0:=2+9/\gamma_0$.

\emph{Step 3: absolute continuity on $\Theta_\Delta$.} Let $g\ge0$ be continuous with compact support in $\Theta_\Delta$, and $P_t(\theta)=(1-t^2)/|1-te^{i\theta}|^2$. Then $\frac1{2\pi}\int g(\theta)P^U_\psi(te^{i\theta})\,d\theta=\int(P_t*g)\,d\mu^U_\psi\to\int g\,d\mu^U_\psi$ as $t\uparrow1$, and by Step 2 the left-hand side is at most $\frac{B_0}{2\pi}\int g\,d\theta$. Hence $\mu^U_\psi|_{\Theta_\Delta}$ is absolutely continuous with density at most $B_0/2\pi$.

\emph{Step 4: conclusion.} Since $\int|e^{i\theta}-\lambda_0|^2d\mu^U_\psi=\varrho^2$, $\mu^U_\psi(\T\setminus\Theta_\Delta)\le\Delta^2$. (a) By Step 3 the atoms of $\mu^U_\psi$ lie outside $\Theta_\Delta$, so $\|P_{\rm pp}(U)\psi\|\le\Delta$, and $\|P_{\rm pp}(U)\phi_0\|\le\Delta+\|\phi_0-\psi\|$. (b) By Wiener's theorem the limit is $\sum_j\ip{\phi_0}{E^U(\{\zeta_j\})\phi_0}^2\le\|P_{\rm pp}(U)\phi_0\|^4$. (c) Write $\mu^U_\psi=h\,d\theta+\mu_c$ with $h\le B_0/2\pi$ supported in $\Theta_\Delta$ and $\mu_c(\T)\le\Delta^2$. Then $|\ip{\psi}{U^n\psi}|^2\le2|\hat h_n|^2+2\Delta^4$ with $\hat h_n=\int e^{in\theta}h\,d\theta$ and $\sum_{n\in\Z}|\hat h_n|^2=2\pi\int h^2\le B_0$, while $|\ip{\phi_0}{U^n\phi_0}|^2\le2|\ip{\psi}{U^n\psi}|^2+C\Delta^2$.
\end{proof}

\paragraph{On the hypothesis.} $\Gamma_r=P^{\hat U}_q(\lambda)/\Delta$ is the golden rule, per unit time, for the residual of the quasi-mode, computed with the continuum $\hat U$ of the scheme itself. Indeed, if $\mu^{\hat U}_q$ has a continuous density near $\lambda$, a formal pole computation for $\ip{\tilde\psi}{(U-z)^{-1}\tilde\psi}=\big(\lambda-z+\varsigma\lambda|\eta|^2/(1+\varsigma\mathcal M(z))\big)^{-1}$ gives $|\ip{\tilde\psi}{U^n\tilde\psi}|^2\approx e^{-n\Delta\Gamma_r}$. Only a lower bound on the spectral density of the coupling vector $q$ is used in the proof, not its regularity. This is the unitary counterpart of the bound ${\rm Im}\,G\le1/(\pi^2|v|^2y)$ for the Green's function of the Friedrichs model \cite{S:Fri48}, \cite[Ch.~XII]{S:RS78}. The condition on $\mathcal M$ is mild: with $w=|a|e^{-i\varrho}$ and $|a|=1-\varrho^2/2$, one finds $\|r\|\le\sqrt2\varrho$, $|\eta|\ge0.47$ and $|\varsigma|=\|r\|/|\eta|\le3\varrho$. The density of $\mu^{\hat U}_q$ near $\lambda$ comes from the resonant spheres and varies on angular scales of order $\Delta^{1/2}$, while $\Theta_\Delta$ and the depth of $\mathcal Z_\Delta$ are of size $\varrho/\Delta$. So (NC$_\Delta$) is expected with $\gamma_0(\Delta)=\Delta\Gamma_r(1-o(1))$ as soon as the amplitude of the residual on the resonant spheres does not cancel. This does not follow from perturbation theory, since the amplitude on the spheres is much smaller than the off-shell part of the residual, and any approximation of $\hat U$ near the spheres, by $U_0$ or by the generalized eigenfunctions of $H$, moves off-shell amplitude onto the spheres with an error of the same order (section~\ref{sec:Sevidence}). At fixed $\Delta$, the argument with (NC$_\Delta$) on a single circle $|z|=t$ bounds the averages up to $N\le(1-t)^{-1}$, and $P^{\hat U}_q(z)\ge(1-|z|^2)|\ip{g}{q}|^2/\|(\hat U^*-\bar z)g\|^2$ for every $g\ne0$, so that such a version can be checked with finitely many applications of the scheme. A proof as $\Delta\to0$ would require an inner problem at the complex singularity of $V$ nearest to the real axis, as for Stokes constants in exponentially small splitting problems.

\subsection{Numerical evidence}\label{sec:Sevidence}
We test the results of this section on the P\"oschl--Teller well, for which the golden rule is explicit.

\paragraph{The P\"oschl--Teller well.} Let $d=1$ and $V(x)=-2\,{\rm sech}^2x$. Then $\phi_0=\frac1{\sqrt2}{\rm sech}\,x$, $E_0=-1$, and $H$ is reflectionless with generalized eigenfunctions $\psi_k=e^{ikx}(\tanh x-ik)/\sqrt{1+k^2}$. A direct computation gives $\ip{\psi_k}{V\phi_0}=-i\frac{\sqrt2\pi k}3\sqrt{1+k^2}\,{\rm sech}\frac{\pi k}2$. Hence
\begin{equation}\label{eq:fgrPT}
\Gamma_{\rm FGR}(\Delta)=\frac{4\pi^3}{9\Delta}\sum_{m\ge1}m\,k_m\,{\rm sech}^2\Big(\frac{\pi k_m}2\Big),\qquad k_m=\sqrt{\frac{2\pi m}\Delta-1},
\end{equation}
so that $\Gamma_{\rm FGR}(\Delta)\sim\frac{16\pi^3}{9\Delta}\sqrt{2\pi/\Delta}\,e^{-\pi\sqrt{2\pi/\Delta-1}}$. Here $\ip{\phi_0}{|\nabla V|^2\phi_0}=\frac{128}{105}$, so that \eqref{eq:quasiE} gives $e_2=32/315\approx0.1016$. Fitting the phase of $\ip{\phi_0}{U_{\rm S}^n\phi_0}$ on $0\le n\Delta\le10$ gives $(E(\Delta)-E_0)/\Delta^2=0.1056$ and $0.1025$ for $\Delta=0.2$ and $0.1$, in agreement with $e_2+O(\Delta^2)$, while $|\ip{\phi_0}{U_{\rm S}^n\phi_0}|\ge0.9996$ on this interval.

\begin{figure}[ht]
\centering
\includegraphics[width=\textwidth]{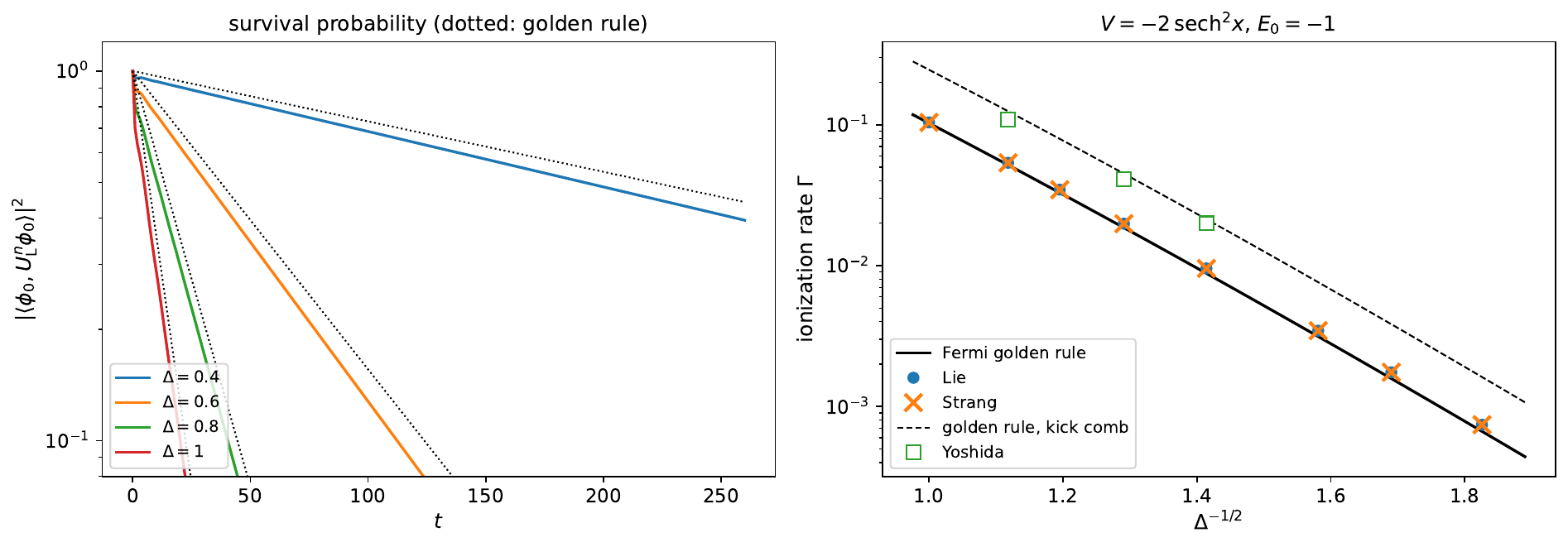}
\caption{Ionization of the P\"oschl--Teller bound state. Left: survival probability under Lie splitting (solid) and $e^{-\Gamma_{\rm FGR}t}$ (dotted). Right: measured rates against $\Delta^{-1/2}$, with the golden rules \eqref{eq:fgrPT} (solid) and \eqref{eq:fgr-comp} for the Yoshida triple jump (dashed).}
\label{fig:ionisation}
\end{figure}

\begin{figure}[ht]
\centering
\includegraphics[width=0.8\textwidth]{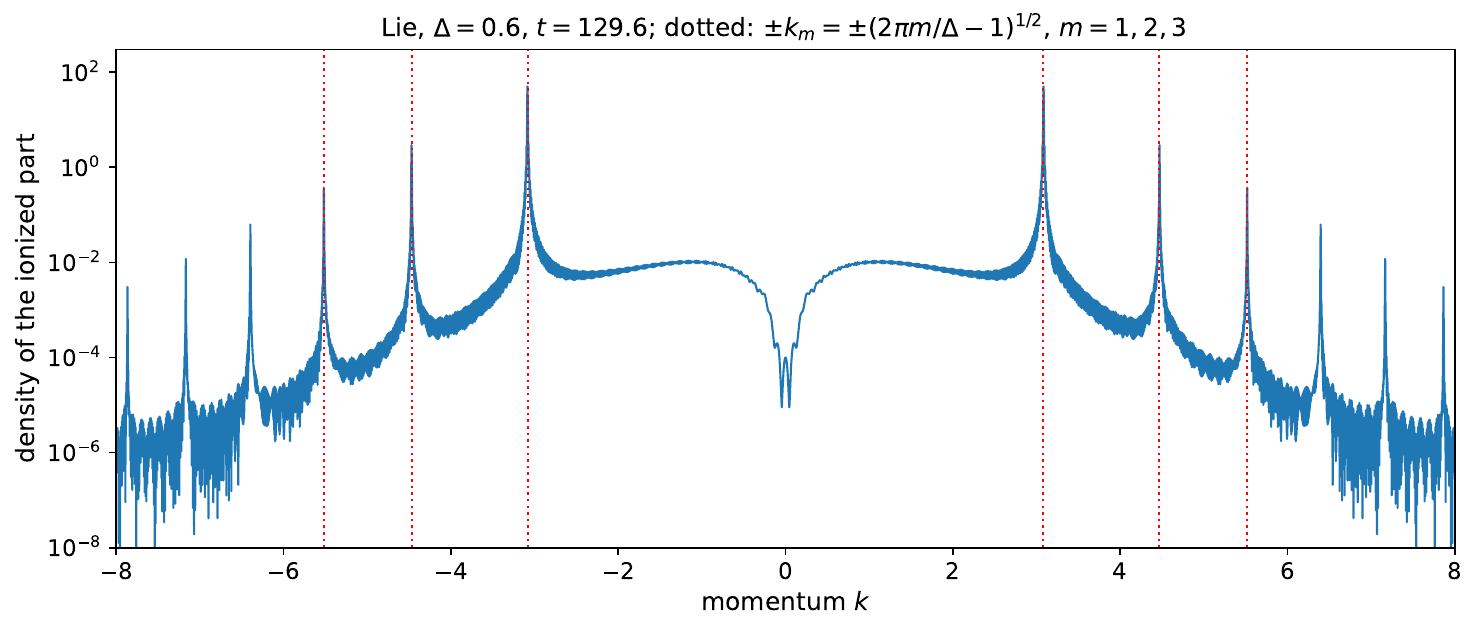}
\caption{Momentum density of the ionized part ($|x|>30$) of $U_{\rm L}^n\phi_0$ for $\Delta=0.6$ at $t=129.6$ (216 steps), without absorber. Dotted lines: $\pm k_m$, $m=1,2,3$; the lines at larger $|k|$ are the channels $m\ge4$.}
\label{fig:spectrum}
\end{figure}

\paragraph{Ionization rates.} We computed $P_n=|\ip{\phi_0}{U^n\phi_0}|^2$ by a Fourier split-step method on $[-1600,1600]$ with $2^{16}$ points and a smooth absorbing layer for $|x|>1000$, up to $T=\min(6/\Gamma_{\rm FGR},260)$, and fitted rates to $\log P_n$ for $t>\min(T/4,40)$, after the initial transient. As a control, Strang splitting with step $0.02$, for which \eqref{eq:fgrPT} gives $3.5\cdot10^{-20}$, gives $P=1-1.5\cdot10^{-8}$ at $t=1000$. Figures~\ref{fig:ionisation} and~\ref{fig:spectrum} display four main features. First, after the metastable regime of Theorem~\ref{thm:meta}, the decay is exponential, and \eqref{eq:fgrPT} has the correct dependence on $\Delta^{-1/2}$: the ratio $\Gamma/\Gamma_{\rm FGR}$ of the measured rate to the golden rule decreases slowly from $1.114$ at $\Delta=0.3$ to $1.01$ at $\Delta=1$ (Strang splitting). Second, Lie and Strang splitting give the same rate to four digits, in agreement with Proposition~\ref{prop:boundproved}(d): since $U_{\rm L}=K(\frac12)^*U_{\rm S}K(\frac12)$, $\ip{\phi_0}{U_{\rm L}^n\phi_0}=\ip{K(\frac12)\phi_0}{U_{\rm S}^nK(\frac12)\phi_0}$ is the Strang survival amplitude of a vector $O(\Delta)$-close to $\phi_0$, which changes the transient but not the decay rate. Third, the ionized part is emitted in sharp lines at $\pm k_m$ (Figure~\ref{fig:spectrum}); for $\Delta=0.6$ at $t=129.6$ the lines $m=1,2,3$ carry $0.919$, $0.037$ and $0.004$ of the ionized mass, the lines $m\ge4$ less than $10^{-3}$, and a broad background the remaining $4\%$. Within the lines the weights are $0.956$, $0.039$ and $0.004$, against $\Gamma_m/\Gamma_{\rm FGR}=0.962$, $0.036$ and $0.002$. Fourth, for $\Delta=0.3$ the lifetime of the bound state under the scheme is about $1300$ time units, and it grows exponentially in $\Delta^{-1/2}$ as $\Delta\to0$, like $1/\Gamma_{\rm FGR}$.

\paragraph{A fourth-order method.} For the Yoshida triple jump, the measured rates at $\Delta=0.5$, $0.6$ and $0.8$ are $2.00\cdot10^{-2}$, $4.11\cdot10^{-2}$ and $1.09\cdot10^{-1}$, that is, $0.94$, $0.91$ and $0.87$ times the golden rule \eqref{eq:fgr-comp}, and $2.09$, $2.07$ and $2.02$ times the Lie rates (Figure~\ref{fig:ionisation}). The fourth-order method thus ionizes about twice as fast as Lie splitting at the same step. At equal cost, the comparison is with Strang splitting at step $\Delta/3$, whose first channel sits at energy $2\pi\cdot3/\Delta$; for example, \eqref{eq:fgrPT} gives $\Gamma_{\rm Strang}(0.2)\approx4.5\cdot10^{-5}$, about $900$ times smaller than $\Gamma_{\rm Yoshida}(0.6)$. The aliasing frequency is set by the period of the method, not by its substeps.

\paragraph{The residual golden rule.} For the P\"oschl--Teller well of section~\ref{sec:Sevidence} and the Strang step, we computed $\tilde\psi$ as the eigenvector of the compressed step $U_K$ near $e^{-i\Delta E_0}$, as in Step 3 of the proof of Theorem~\ref{thm:Sexpmeta} but with the larger cutoffs $K=0.75k_1$ or $0.9k_1$, $k_1=(2\pi/\Delta-1)^{1/2}$, then $\lambda$, $\hat U$ and $q$ as above, and $P^{\hat U}_q$ near $\lambda$ from the correlations $\ip{q}{\hat U^nq}$, $n\le4000$, smoothed by a Gaussian window of width $0.02$ and $0.01$ and extrapolated to width $0$ ($2^{15}$ Fourier modes on $[-160,160]$, absorbing factor $\cos^2(\pi(|x|-110)/100)$ for $|x|>110$ after each step; script \texttt{residual\_rate.py}). $|\varsigma\mathcal M|$ is sampled on the closure of $\mathcal Z_\Delta$. Table~\ref{tab:Sresid} displays three main features. First, the residual golden rule reproduces the measured rates, to $0.4\%$ at $\Delta=0.3$, with a discrepancy of about $0.18\,|\varsigma|$ for $K=0.75k_1$. Second, for $K/k_1\in\{0.75,0.9\}$ it is independent of the cutoff to $0.2\%$ for $\Delta\le0.3$ ($1.22$ and $1.12$ for $K/k_1=0.5$ and $0.95$ at $\Delta=0.3$). Third, its ratio to the Fermi golden rule \eqref{eq:fgrPT} converges to about $1.126\pm0.002$ as $\Delta\to0$: \eqref{eq:fgrPT} has the right exponential rate but not the right prefactor. In addition, the condition on $\mathcal M$ in (NC$_\Delta$) holds with a wide margin, and at $\Delta=0.3$ the density of $\mu^{\hat U}_q$ varies by about $2\%$ over $\pm0.02$ rad around $\lambda$.

\begin{table}[ht]
\centering
\setlength{\tabcolsep}{4pt}
\begin{tabular}{cccccc}
\toprule
$\Delta$ & $\varrho$ & $\sup_{\mathcal Z_\Delta}|\varsigma\mathcal M|$ & $\Gamma_r/\Gamma_{\rm FGR}$, $K=0.75k_1$ & $\Gamma_r/\Gamma_{\rm FGR}$, $K=0.9k_1$ & measured $\Gamma/\Gamma_{\rm FGR}$\\\midrule
1.00 & $0.21$ & $0.34$ & 0.927 & 0.942 & 1.012\\
0.60 & $7.1\cdot10^{-2}$ & $0.12$ & 1.039 & 1.043 & 1.068\\
0.30 & $1.0\cdot10^{-2}$ & $1.8\cdot10^{-2}$ & 1.110 & 1.110 & 1.114\\
0.20 & $2.2\cdot10^{-3}$ & $3.9\cdot10^{-3}$ & 1.124 & 1.124 & --\\
0.10 & $7.7\cdot10^{-5}$ & $1.3\cdot10^{-4}$ & 1.127 & 1.128 & --\\
0.05 & $7.0\cdot10^{-7}$ & $1.0\cdot10^{-6}$ & 1.125 & 1.128 & --\\\bottomrule
\end{tabular}
\caption{Residual golden rule $\Gamma_r=P^{\hat U}_q(\lambda)/\Delta$ for two cutoffs, against the Fermi golden rule \eqref{eq:fgrPT} and the measured rates (Strang splitting). $\varrho$ and $\sup_{\mathcal Z_\Delta}|\varsigma\mathcal M|$ are for $K=0.9k_1$.}
\label{tab:Sresid}
\end{table}

These values also show why perturbation theory is not enough to verify (NC$_\Delta$). Approximating the continuum near the resonant spheres by plane waves gives the amplitudes $|\hat r(\pm k_m)|^2/(2\Delta^2k_m)$, which yield $1.05$, $1.24$ and $1.32$ times $\Gamma_{\rm FGR}$ at $\Delta=1$, $0.6$ and $0.3$ ($K=0.75k_1$), and the generalized eigenfunctions $\psi_k$ of $H$ give values between $1.02$ and $1.17$ depending on the cutoff, for $K/k_1\in\{0.5,0.75,0.9,0.95\}$, against $0.93$, $1.04$ and $1.11$ for $\hat U$ ($K=0.75k_1$).

%=====================================================================
\section{Long-range potentials}\label{sec:Slongrange}
%=====================================================================

For $\rho\le1$ the wave operators $\Omega^c_\pm$ and $\Omega^U_\pm$ do not exist in general, but the relative wave operators $\Xi_\pm$ and the relative scattering operator $\Sigma_U=\Xi_+^*\Xi_-$ exist for every $\rho>0$, with the rates of Theorems~\ref{thm:conv}(d) and~\ref{thm:S}(a). This section composes them with modified wave operators of $H$. Throughout, $U$ is a composition of order $r$, and $G=g(H)\jap{x}^{-s}$ with $g\in C_c^\infty((0,\infty))$ and $s>1$.

\subsection{Modified splitting wave operators}\label{sec:Slrresults}

Modified wave operators of $H$ are defined by a modifier, which we use in the following abstract form.

\paragraph{Hypothesis (M).}
$\Lambda_\pm(t,\xi)$ are real functions, smooth in $\xi\in\R^d\setminus\{0\}$. The limits
\[
\Omega^c_{\Lambda,\pm}=\slim_{t\to\pm\infty}e^{itH}e^{-itH_0}e^{-i\Lambda_\pm(t,p)}
\]
exist and are complete, $\Ran\Omega^c_{\Lambda,\pm}=\HH_{\rm ac}(H)$. Moreover, for every $\tau$ and $\xi\ne0$, $\Lambda_\pm(t+\tau,\xi)-\Lambda_\pm(t,\xi)\to0$ as $t\to\pm\infty$.

The last condition implies the intertwining relation $e^{-i\tau H}\Omega^c_{\Lambda,\pm}=\Omega^c_{\Lambda,\pm}e^{-i\tau H_0}$. We write $W^c_{\Lambda,\pm}=(\Omega^c_{\Lambda,\pm})^*$ and $S^c_\Lambda=(\Omega^c_{\Lambda,+})^*\Omega^c_{\Lambda,-}$.

\emph{Examples.}
\emph{Dollard modifier.} $\Lambda_\pm(t,\xi)=\Phi(t,\xi)=\int_0^tV(2s\xi)\,ds$. For $V$ satisfying \eqref{ass:V} with $\rho>\frac12$, existence and completeness are classical (\cite{S:Dol64,S:RS79}, \cite[Theorem~5.2]{S:JMP84}, the latter for $-\frac12\Delta_x$, which amounts to a rescaling of time). The last condition in (M) holds because $V(2s\xi)\to0$. \emph{Modifiers from a Hamilton--Jacobi equation.} For $0<\rho\le\frac12$, the Dollard phase is replaced by a solution of a Hamilton--Jacobi equation \cite{S:Hor76}; asymptotic completeness under \eqref{ass:V} is proved in \cite{S:DG97}.

\begin{theorem}[Modified splitting wave operators]\label{thm:modWO}
Assume \eqref{ass:V} with $\rho>0$ and (M). For every $\Delta\in(0,1]$:
\\[1mm]
\textup{(a)} The limits $\Omega^U_{\Lambda,\pm}(\Delta)=\slim_{n\to\pm\infty}U^{-n}U_0^ne^{-i\Lambda_\pm(n\Delta,p)}$ exist and equal $\Xi_\pm(\Delta)\Omega^c_{\Lambda,\pm}$. They are isometric, $U\Omega^U_{\Lambda,\pm}=\Omega^U_{\Lambda,\pm}U_0$, and $\Ran\Omega^U_{\Lambda,\pm}=\Ran\Xi_\pm\subset\HH_{\rm ac}(U)$.
\\[1mm]
\textup{(b)} $\|(\Omega^U_{\Lambda,\pm}(\Delta)-\Omega^c_{\Lambda,\pm})W^c_{\Lambda,\pm}G\|\le C\Delta^r$, and $\Omega^U_{\Lambda,\pm}(\Delta)\to\Omega^c_{\Lambda,\pm}$ strongly as $\Delta\to0$.
\\[1mm]
\textup{(c)} The modified splitting scattering operator $S^U_\Lambda=\Omega^U_{\Lambda,+}(\Delta)^*\Omega^U_{\Lambda,-}(\Delta)$, satisfies $S^U_\Lambda=(\Omega^c_{\Lambda,+})^*\,\Sigma_U\,\Omega^c_{\Lambda,-}$. If $U=P\tilde UP^{-1}$ with $\tilde U$ of order $q$, then
\[
|\ip{W^c_{\Lambda,+}Gu'}{(S^U_\Lambda-S^c_\Lambda)W^c_{\Lambda,-}Gu}|\le C\Delta^q\|u\|\,\|u'\| .
\]
\end{theorem}

The proofs are given in section~\ref{sec:Slrproofs}, which also collects the orders of convergence. Theorem~\ref{thm:modWO} follows from Theorem~\ref{thm:conv}(d) and contains no new error bound. With $W^U_{\Lambda,\pm}=(\Omega^U_{\Lambda,\pm})^*$, it gives $e^{i\Lambda_+(n\Delta,p)}U_0^{-n}U^n\psi\to W^U_{\Lambda,+}\psi$ for $\psi\in\Ran\Xi_+$, so that asymptotic states and $S^U_\Lambda$ are obtained from the scheme and $\Lambda_\pm$; this is needed for phases of scattering data, which depend on the modifier. Momentum distributions need no modifier: $|\widehat{U^n\psi}|^2\to|\widehat{W^U_{\Lambda,+}\psi}|^2$ in $L^1$, a limit independent of $\Lambda_+$. We do not know whether $\Ran\Xi_\pm=\HH_{\rm ac}(U)$ for $\rho\le1$.

\subsection{Proofs}\label{sec:Slrproofs}

\begin{proof}[Proof of Theorem~\ref{thm:modWO}]
(a) $U^{-n}U_0^ne^{-i\Lambda_\pm(n\Delta,p)}f=U^{-n}e^{-in\Delta H}\big(e^{in\Delta H}e^{-in\Delta H_0}e^{-i\Lambda_\pm(n\Delta,p)}f\big)$. The bracket converges in norm to $\Omega^c_{\Lambda,\pm}f\in\HH_{\rm ac}(H)$, which gives existence and the identity $\Omega^U_{\Lambda,\pm}=\Xi_\pm\Omega^c_{\Lambda,\pm}$. Intertwining follows from $U\Xi_\pm=\Xi_\pm e^{-i\Delta H}$ and $e^{-i\Delta H}\Omega^c_{\Lambda,\pm}=\Omega^c_{\Lambda,\pm}U_0$. The range statement follows from completeness of $\Omega^c_{\Lambda,\pm}$.
(b) For $\psi=Gu$ let $f=W^c_{\Lambda,\pm}\psi$; then $\Omega^c_{\Lambda,\pm}f=\psi$ and $(\Omega^U_{\Lambda,\pm}-\Omega^c_{\Lambda,\pm})f=(\Xi_\pm-1)\psi$. Apply Theorem~\ref{thm:conv}(d). Strong convergence follows since $W^c_{\Lambda,\pm}G\HH$ is dense in $\HH$.
(c) The identity follows from (a). Since $\Omega^c_{\Lambda,\pm}W^c_{\Lambda,\pm}=P_{\rm ac}(H)$, we get $\ip{W^c_{\Lambda,+}\psi'}{(S^U_\Lambda-S^c_\Lambda)W^c_{\Lambda,-}\psi}=\ip{\psi'}{(\Sigma_U-P_{\rm ac}(H))\psi}$, and Theorem~\ref{thm:S}(a) applies.
\end{proof}

Thus the choice of modifier only relabels asymptotic states.

\medskip\noindent{\bf Orders of convergence.} The following table collects the rates, all on energy-localized data and uniform in $\Delta\in(0,1]$. The rows for $\Sigma_U$ and $S^U_\Lambda$ refer to $U=P\tilde UP^{-1}$ with $\tilde U$ of order $q$.

\begin{center}
\renewcommand{\arraystretch}{1.25}
\begin{tabular}{lccccl}
quantity & Lie & Strang & order $r$ & eff.\ order $q$ & range\\\hline
$\sup_n\|(U^n-e^{-in\Delta H})G\|$ & $\Delta$ & $\Delta^2$ & $\Delta^r$ & & $\rho>0$\\
$\|(\Xi_\pm-1)G\|$ & $\Delta$ & $\Delta^2$ & $\Delta^r$ & & $\rho>0$\\
$\|(\Xi^{\rm L}_\pm-1-\frac{i\Delta}2V)G\|$ & $\Delta^2$ & & & & $\rho>0$\\
$\|P_{\rm pp}(U)G\|$ & $\Delta$ & $\Delta^2$ & $\Delta^r$ & & $\rho>0$\\
$\|G^*(\Sigma_U-P_{\rm ac}(H))G\|$ & $\Delta^2$ & $\Delta^2$ & & $\Delta^q$ & $\rho>0$\\
$\|(\Omega^U_{\Lambda,\pm}-\Omega^c_{\Lambda,\pm})W^c_{\Lambda,\pm}G\|$ & $\Delta$ & $\Delta^2$ & $\Delta^r$ & & (M)\\
modified $S^U_\Lambda$ (Thm.~\ref{thm:modWO}(c)) & $\Delta^2$ & $\Delta^2$ & & $\Delta^q$ & (M)\\
$\lim_n\|(U^n-e^{-in\Delta H})\psi\|$ & $\frac\Delta2\|V\psi\|+O(\Delta^2)$ & $O(\Delta^2)$ & $O(\Delta^r)$ & & $\rho>0$
\end{tabular}
\end{center}

The last line follows from $\lim_n\|(U^n-e^{-in\Delta H})\psi\|=\|(\Xi_+-1)\psi\|$ and Theorem~\ref{thm:conv}(d). Condition (M) holds with the Dollard modifier for $\rho>\frac12$, and with the modifiers constructed from a Hamilton--Jacobi equation for every $\rho>0$. Numerical results for $\rho=0.7$ and $\rho=0.3$ are given in section~\ref{sec:Snumerics}.

\subsection{Modifiers adapted to the scheme}\label{sec:Sextensions}

A splitting scheme suggests its own Dollard phase, namely the potential sampled along free trajectories at the kick times:
\[
\Phi^U_n(\xi)=\Delta\sum_{k=0}^{n-1}\sum_ja_j\,V\big(2(k+\beta_j)\Delta\xi\big),\qquad\beta_j=b_1+\dots+b_{j-1},
\]
for $n\ge0$, and $\Phi^U_{n}(\xi)=-\Delta\sum_{k=n}^{-1}\sum_ja_jV\big(2(k+\beta_j)\Delta\xi\big)$ for $n<0$.
Write $h_\xi(s)=V(2s\xi)$. For $\xi\ne0$, $|h_\xi^{(l)}(s)|\le C\jap{s}^{-\rho-l}$, and each step of $\Phi^U_n-\Phi(n\Delta,\cdot)$ is a quadrature error of size $O(\Delta^2\jap{k\Delta}^{-\rho-1})$, since $\sum_ja_j=1$. Hence the phase difference converges pointwise, $\Phi^U_n(\xi)-\Phi(n\Delta,\xi)\to c_\Delta(\xi)$, and by dominated convergence
\[
\slim_{n\to\pm\infty}U^{-n}U_0^ne^{-i\Phi^U_n(p)}=\Omega^U_{\Phi,\pm}(\Delta)\,e^{-ic_\Delta(p)}
\]
whenever the Dollard modifier satisfies (M). For $n\to+\infty$, the Euler--Maclaurin formula on $[0,\infty)$, whose boundary terms at infinity vanish, gives explicitly:
For Lie splitting, whose kicks sample the potential by the left rectangle rule, $c_\Delta(\xi)=\frac\Delta2V(0)-\frac{\Delta^2}6\,\xi\cdot\nabla V(0)+O(\Delta^4)$. For Strang splitting, which corresponds to the trapezoidal rule, $c_\Delta(\xi)=-\frac{\Delta^2}6\,\xi\cdot\nabla V(0)+O(\Delta^4)$.
The constant phase $\frac\Delta2V(0)$ for Lie splitting is the long-range counterpart of the first-order correction $1+\frac{i\Delta}2V$ in Theorem~\ref{thm:conv}(d). With the scheme-adapted modifier, the discrete modified wave operators differ from the continuous ones by the explicit phase $e^{-ic_\Delta(p)}$ in addition to the $O(\Delta^r)$ error of Theorem~\ref{thm:modWO}(b).

%=====================================================================
\section{Scalar elliptic operators}\label{sec:Sscalar}
%=====================================================================

We extend in this section the results of the paper from $H_0=-\Delta_x$ to semibounded elliptic differential operators with constant coefficients. The fixed-$\Delta$ theory (Theorems~\ref{thm:A} and~\ref{thm:B}) extends, with the changes in the proofs given below. The convergence theory extends under a local decay hypothesis, which for $-\Delta_x$ is the result of \cite{S:JMP84} used in the paper.

\paragraph{Hypothesis (E).} $H_0=\omega(p)$, where $\omega$ is a real polynomial on $\R^d$ of degree $2m\ge2$ whose principal part satisfies $\omega_{2m}(\xi)>0$ for $\xi\ne0$.

Under (E), $\omega(p)$ is self-adjoint on $\mathcal H^{2m}$ and bounded below, and so is $H=\omega(p)+V$ under \eqref{ass:V}. Euler's identity $\xi\cdot\nabla\omega_{2m}=2m\omega_{2m}$ gives $|\nabla\omega(\xi)|\ge c|\xi|^{2m-1}$ for $|\xi|\ge R$. The critical set $\{\nabla\omega=0\}$ is therefore compact, and the set $\mathcal T(\omega)$ of critical values of $\omega$ is finite, since $\omega$ is a polynomial. We call $\mathcal T(\omega)$ the thresholds of $\omega$, and set $\Theta_\omega=\varepsilon_\Delta(\mathcal T(\omega))\subset\T$, a finite set. For $\omega(\xi)=|\xi|^2$, $\mathcal T(\omega)=\{0\}$ and $\Theta_\omega=\{1\}$. The critical set is a proper algebraic set, hence a null set, and $\omega$ is a submersion outside it, so that preimages under $\omega$ of null sets are null sets. The spectral measures of $\omega(p)$ are therefore absolutely continuous, and $U_0=e^{-i\Delta\omega(p)}$ is purely absolutely continuous with $\sigma(U_0)=\T$.

\paragraph{Conjugate operator.} Let $v=\nabla\omega/(1+|\nabla\omega|^2)$ and $\gamma_\omega=v\cdot\nabla\omega=|\nabla\omega|^2/(1+|\nabla\omega|^2)$. Since the denominator is a polynomial bounded below by $c\jap{\xi}^{4m-2}$, the components of $v$ satisfy $|\partial^\beta v_j|\le C_\beta\jap{\xi}^{1-2m-|\beta|}$, so that $v_j$, $\dv v$ and $\gamma_\omega$ belong to the class $\Mfr$ of section~\ref{sec:scattering}, $\dv v\in C_0(\R^d)$, and $|v|\le\frac12$. Proposition~\ref{prop:group} and Propositions~\ref{prop:mult} and~\ref{prop:weights} hold with the same proofs, which use only these properties. Lemma~\ref{lem:free}(a) becomes $F_\tau^*[A,F_\tau]=\tau\Delta\,\gamma_\omega(p)$ for $F_\tau=e^{-i\tau\Delta\omega(p)}$, since $v\cdot\nabla e^{-i\tau\Delta\omega}=-i\tau\Delta \gamma_\omega e^{-i\tau\Delta\omega}$. For an arc $\Theta$ whose closure avoids $\Theta_\omega$, the set $Z=\{\xi:\varepsilon_\Delta(\omega(\xi))\in\overline\Theta\}$ is closed and contains no critical point, and $\gamma_\omega\to1$ at infinity, so that $\inf_Z\gamma_\omega>0$. This replaces Lemma~\ref{lem:free}(c).

\paragraph{Kicks.} The explicit computation in Proposition~\ref{prop:comm-a} uses $v=p\jap{p}^{-2}$. In general we use the SG calculus \cite{S:Cor95}, with symbol classes $SG^{\mu,\nu}$ defined by $|\partial_x^\alpha\partial_\xi^\beta a|\le C_{\alpha\beta}\jap{x}^{\nu-|\alpha|}\jap{\xi}^{\mu-|\beta|}$. Operators with symbols in $SG^{\mu,\nu}$ map $\Sch$ to $\Sch$, are bounded for $\mu,\nu\le0$, and are compact for $\mu,\nu<0$. Commutators satisfy $[\Op SG^{\mu_1,\nu_1},\Op SG^{\mu_2,\nu_2}]\subset\Op SG^{\mu_1+\mu_2-1,\nu_1+\nu_2-1}$. The Weyl symbol $x\cdot v(\xi)$ of $A$ belongs to $SG^{1-2m,1}$, and $\nu_a=K(a)-1\in SG^{0,-\rho}$ by Lemma~\ref{lem:kick}(a). Hence $\mathrm{ad}_A^k(\nu_a)\in\Op SG^{0,-\rho}$ for all $k$, so that $K(a)\in C^k(A)$ by Proposition~\ref{prop:core} and Fact~\ref{fact:reg}(c), and $[A,K(a)]\in\Op SG^{-2m,-\rho}$ and $[\gamma_\omega(p),K(a)]\in\Op SG^{-1,-\rho-1}$ are compact. This replaces Lemma~\ref{lem:kick}(b) and Propositions~\ref{prop:comm-a} and~\ref{prop:algebra}. Lemma~\ref{lem:product} then holds with $2\Delta(1-R)$ replaced by $\Delta \gamma_\omega(p)$.

\paragraph{Compactness after a momentum cut-off, and the Mourre estimate.} Lemma~\ref{lem:lowenergy} holds with $\eta(H_0)$ replaced by $\eta(p)$, $\eta\in C_0(\R^d)$, with the same proof, since $\nu_i\eta(p)\in\KK$ by Fact~\ref{fact:compact}. In Lemma~\ref{lem:mourre}, let $\Theta$ be an arc whose closure avoids $\Theta_\omega$, and $\varphi\in C(\T;[0,1])$ equal to $1$ on $\overline\Theta$ with support avoiding $\Theta_\omega$. Since $\{|\nabla\omega|\le\delta\}$ is compact and shrinks to the critical set as $\delta\to0$, and $\omega$ is continuous, there is $c'>0$ such that $\varepsilon_\Delta(\omega(\xi))\notin\supp\varphi$ whenever $\gamma_\omega(\xi)\le2c'$. Choose $\eta\in C_c(\R^d;[0,1])$ equal to $1$ on $\{\gamma_\omega\le c'\}$ and supported in $\{\gamma_\omega<2c'\}$. Then $\gamma_\omega\ge c'(1-\eta)$ pointwise and $\varphi(U_0)\eta(p)=0$. Steps~2--5 of the proof of Lemma~\ref{lem:mourre} then give
\[
E^U(\Theta)U^*[A,U]E^U(\Theta)\ge c_0\,E^U(\Theta)+\mathsf K,\qquad \mathsf K\in\KK,
\]
for every $c_0<\Delta c'$. Consequently, Theorems~\ref{thm:A} and~\ref{thm:B} hold under (E) with $\T\setminus\{1\}$ replaced by $\T\setminus\Theta_\omega$. In the proof of Theorem~\ref{thm:B}, the free decay estimate of Lemma~\ref{prop:decay} holds for $e^{-it\omega(p)}\chi(p)$ with $\chi\in C_c^\infty$ supported away from the critical set, by the same non-stationary phase argument with $|\nabla\Phi|=|z-t\nabla\omega(\xi)|\ge|t|c_\chi-|z|$, $c_\chi=\min_{\supp\chi}|\nabla\omega|$. The continuous wave operators exist and are complete for $\rho>1$ \cite[Ch.~XIV]{S:Hor83}.

\paragraph{The spectrum of $H$.} The Mourre estimate for $H$ with the same conjugate operator holds on every compact interval disjoint from $\mathcal T(\omega)$, since $[\omega(p),iA]=\gamma_\omega(p)$ and $[V,iA]$ is compact. Hence $H$ has no singular continuous spectrum, and its eigenvalues outside $\mathcal T(\omega)$ have finite multiplicity and accumulate at most at $\mathcal T(\omega)$ \cite[Ch.~7]{S:ABG96}. In contrast with $-\Delta_x$, embedded eigenvalues may occur for $m\ge2$; see \cite{S:DSS11} for the bilaplacian. The energy-localized data are therefore $g(H)\jap{x}^{-s}$ with $g\in C_c^\infty(\R\setminus(\mathcal T(\omega)\cup\sigma_p(H)))$, whose ranges are dense in $\HH_{\rm ac}(H)$.

\paragraph{Hypothesis (LD).} For every $g\in C_c^\infty(\R\setminus(\mathcal T(\omega)\cup\sigma_p(H)))$, $\sigma>0$, $\epsilon>0$ and $K\ge0$, $\|\jap{p}^K\jap{x}^{-\sigma}e^{-itH}g(H)\jap{x}^{-\sigma}\|\le C\jap{t}^{-\sigma+\epsilon}$.

For $\omega(\xi)=|\xi|^2$ this is Proposition~\ref{prop:localdecay}. For general $\omega$, we state the convergence results under (LD).

\paragraph{Convergence theory.} Under (E) and (LD), Theorems~\ref{thm:conv}, \ref{thm:wo}, \ref{thm:S} and~\ref{thm:modWO} hold with the same rates. The Sobolev orders depend on $m$, and Hypothesis~(M) is used with modifiers adapted to $\omega$, e.g., the Dollard phase $\int_0^tV(s\nabla\omega(\xi))\,ds$. The proofs change in three places.
\textup{(i)} \emph{Sobolev and weighted estimates.} In Fact~\ref{fact:sobolev}, $\mathrm{Dom}((H+c)^j)=\mathcal H^{2mj}$, with the same induction. In Fact~\ref{fact:weights}(b,c), the rescaling of Lemma~\ref{lem:fcalc} applies to $p_h(x,\xi)=h^{2m}\omega(\xi/h)+h^{2m}V(x)$, which is bounded in $S(\jap{\xi}^{2m})$ and uniformly elliptic for $h\in(0,1]$ because $\omega$ is a polynomial, with $\Op_h(p_h)=h^{2m}H$. The exact Egorov theorem ($\Psi$3) is specific to quadratic symbols. It is replaced by the exact identity $e^{i\tau\omega(p)}x_je^{-i\tau\omega(p)}=x_j+\tau(\partial_j\omega)(p)$, whose right-hand sides commute. For integer $n$ it gives $e^{-i\tau\omega(p)}\jap{x}^{2n}=(1+|x-\tau\nabla\omega(p)|^2)^ne^{-i\tau\omega(p)}$, and the polynomial on the right is a sum of terms $x^\alpha Q_\alpha(p)$ with $|\alpha|\le2n$ and $\deg Q_\alpha\le(2n-|\alpha|)(2m-1)$. Hence $\jap{p}^k\jap{x}^{-2n}e^{-i\tau\omega(p)}\jap{x}^{2n}\jap{p}^{-k-2n(2m-1)}$ is bounded uniformly for $|\tau|\le T_0$; non-integer weights follow by complex interpolation. Each unit of weight passed through a free flow thus costs $2m-1$ derivatives instead of $1$.
\textup{(ii)} \emph{Defects.} $[\omega(p),V]$ is a differential operator of order $2m-1$ whose coefficients are derivatives of $V$ of order $\ge1$, so that the Lie defect again decays like $\jap{x}^{-\rho-1}$. For Strang, $e^{-isa}pe^{isa}=p+s\nabla a$ gives $e^{-isa}\omega(p)e^{isa}=\omega(p+s\nabla a)$ exactly, a polynomial in $s$ whose coefficient of $s^k$ contains at least $k$ factors of derivatives of $a$ of order $\ge1$. The first-order terms cancel as in \eqref{eq:DS}, and the remainder has coefficients $O(\jap{x}^{-\rho-2})$. For a general composition, Lemma~\ref{lem:defect-comp} holds with $\mathcal H^{2r+2}$ replaced by $\mathcal H^{2m(r+1)}$ and $[H_0,V]$ by $[\omega(p),V]$, so that Proposition~\ref{prop:defectbounds} holds with $\mu\le\rho+1$.
\textup{(iii)} \emph{Local decay.} Proposition~\ref{prop:localdecay} is replaced by (LD).

\paragraph{Bound states.} Proposition~\ref{prop:boundproved} and Theorem~\ref{thm:meta} extend with the same proofs, with $\T\setminus\{1\}$ replaced by $\T\setminus\Theta_\omega$ and the continuum energies $E_0+2\pi\ell/\Delta$, $\ell\ge1$, above $\min\omega$. In the proof of Theorem~\ref{thm:noeig}, the spheres $|\xi|^2=E+2\pi\ell/\Delta$ become the level sets $\{\omega=E+2\pi\ell/\Delta\}$. For $\ell$ large these are smooth compact hypersurfaces on which $\nabla\omega\ne0$, so that Step~3 applies, and each ray $\{t\theta:t>0\}$ meets them at points $t_\ell(\theta)$ with $\omega_{2m}(\theta)t_\ell(\theta)^{2m}=(2\pi\ell/\Delta)(1+o(1))$. In Step~4 these zeros contribute $c_\theta\Delta r^{2m}(1+o(1))$ to Jensen's formula, with $c_\theta=(\omega_{2m}(\theta)+\omega_{2m}(-\theta))/(4\pi m)$, instead of $\Delta r^2/(2\pi)$, while the growth bounds of Step~2 do not depend on $\omega$. Cases (a) and (b) therefore extend, and for $m\ge2$ case (c) holds for every $\Delta>0$; the threshold $\Delta>\pi/(4a)$ is specific to $\omega(\xi)=|\xi|^2$. We have not checked the extension of Theorem~\ref{thm:expmeta}, whose proof uses the compression to $|\xi|\le(\delta/\Delta)^{1/2}$ and the bound $||\xi|^2-|\eta|^2|\le2K|\xi-\eta|$ on $L^2_K$.

%=====================================================================
\section{Technical tools}\label{sec:Stools}
%=====================================================================

This section contains the proofs of the results stated in the appendix of the paper, as well as the number of derivatives of $V$ used in the proofs.

\paragraph{Regularity.} Each result uses only finitely many derivatives; we do not optimize their number, but the proofs give the following counts. The spectral and scattering theory at fixed $\Delta$ (Theorems~\ref{thm:A} and \ref{thm:B}) uses $|\alpha|\le2$. For energy-localized data the count in the convergence theory is set by two ingredients. The first is the local decay estimate of \cite[Theorem~5.1]{S:JMP84} (Proposition~\ref{prop:localdecay}), which needs $|\alpha|\le5$ when $\min(\rho+1,s)>7/6$, with $s$ as in Theorem~\ref{thm:conv}. Below that threshold the number grows like $(\min(\rho+1,s)-1)^{-1}$ as $\min(\rho+1,s)\downarrow1$. The second is the one-step defect. For a method of order $r$ the number of derivatives used by the proof of Proposition~\ref{prop:defectbounds} grows with $r$ and with the number of substeps. For Lie and Strang splitting, $|\alpha|\le4$ suffices for the defects. The proof of Fact~\ref{fact:weights} below uses symbol classes with all derivatives. For its parts (b) and (c), the Helffer--Sj\"ostrand formula \cite[Ch.~8]{S:DS99} combined with commutators $[(H-z)^{-1},\jap{x}^{\mu}]$, which contain no derivatives of $V$, gives the same bounds from finitely many derivatives, as assumed in these counts.

\paragraph{Proof of Proposition~\ref{prop:core}.} For $f\in D$ we have $Sf\in D\subset \mathrm{Dom}(A)$, so $\ip{Af}{Sf}-\ip{f}{SAf}=\ip{f}{(AS-SA)f}$, which is bounded by $c\|f\|^2$. For $f\in \mathrm{Dom}(A)$ choose $f_n\in D$ with $f_n\to f$ and $Af_n\to Af$. Then $\ip{Af_n}{Sf_n}\to\ip{Af}{Sf}$ and $\ip{f_n}{SAf_n}\to\ip{f}{SAf}$, so the bound extends to $\mathrm{Dom}(A)$ and Fact~\ref{fact:reg}(a) applies. By polarization the form is $\ip{f}{[A,S]h}$, and for $f,h\in D$ it equals $\ip{f}{(AS-SA)h}$, which identifies $[A,S]$ on $D$.

\paragraph{Proof of Proposition~\ref{prop:group}.} (a) $v$ is smooth with bounded derivatives of all orders, hence globally Lipschitz, so the flow is global and smooth, and $\phi_t\circ\phi_s=\phi_{t+s}$. Since $|v|\le1/2$, $|\phi_t(\xi)-\xi|\le|t|/2$. The derivative $D\phi_t$ solves $\frac{d}{dt}D\phi_t=Dv(\phi_t)D\phi_t$, so $\|D\phi_t\|\le e^{|t|\|Dv\|_\infty}$ by Gronwall; higher derivatives solve linear inhomogeneous equations with bounded coefficients, and induction on the order gives bounds uniform in $\xi$. Liouville's formula gives $J_t$, and its derivatives are bounded by the same argument.

(b) Unitarity follows from the change of variables $\eta=\phi_t(\xi)$, $d\eta=J_t(\xi)\,d\xi$. The group law follows from $\phi_t\circ\phi_s=\phi_{t+s}$ and $J_{s+t}(\xi)=J_s(\xi)J_t(\phi_s(\xi))$ (chain rule). If $\hat u\in\Sch$, then $\hat u\circ\phi_t$ is smooth, its derivatives are finite sums of products of derivatives of $\hat u$ at $\phi_t(\xi)$ with bounded derivatives of $\phi_t$, and $|\phi_t(\xi)|\ge|\xi|-|t|/2$; with (a) this gives $T_t\Sch\subset\Sch$, and equality follows from $T_{-t}=T_t^{-1}$. Strong continuity holds on $\Sch$ by dominated convergence and extends by density and unitarity.

(c) For $\hat u\in\Sch$, $\frac{d}{dt}T_tu|_{t=0}$ exists in $\HH$ and equals the inverse Fourier transform of $v\cdot\nabla\hat u+\frac12(\dv v)\hat u$. Hence $\Sch\subset \mathrm{Dom}(A)$ and $-iAu$ is given by this expression, which is the second formula in \eqref{eq:A}. Since $\Sch$ is dense, contained in $\mathrm{Dom}(A)$ and invariant under the group, it is a core by the invariant domain theorem \cite[Ch.~VIII]{S:RS80}. For the first formula, $\widehat{x_jw}=i\partial_{\xi_j}\hat w$ gives $\widehat{x\cdot v(p)u}=i\nabla\cdot(v\hat u)=i(\dv v)\hat u+iv\cdot\nabla\hat u$, and $[x_j,v_j(p)]=i(\partial_jv_j)(p)$ gives $\frac12(x\cdot v+v\cdot x)=x\cdot v(p)-\frac i2(\dv v)(p)$.

(d) $(T_tb(p)T_{-t}u)^\wedge(\xi)=J_t^{1/2}(\xi)\,b(\phi_t(\xi))\,\widehat{T_{-t}u}(\phi_t(\xi))=b(\phi_t(\xi))\hat u(\xi)$.

\paragraph{Proof of Fact~\ref{fact:compact}.} For $a,b\in L^2(\R^d)$, $a(x)b(p)$ is Hilbert--Schmidt with $\|a(x)b(p)\|_{\mathscr S_2}=(2\pi)^{-d/2}\|a\|_2\|b\|_2$ \cite[Thm.~4.1]{S:Sim05}. Approximate $a,b\in C_0$ uniformly by $C_c$ functions; norm limits of compact operators are compact.

\paragraph{Proof of Proposition~\ref{prop:mult}.} \eqref{eq:multcomm} follows from the second formula in \eqref{eq:A}: $[iv\cdot\nabla,b]=i(v\cdot\nabla b)$ and $\dv v$ commutes with $b$. The second statement follows from Proposition~\ref{prop:core} with $D=\Sch$, applied inductively to $(v\cdot\nabla)^kb$, together with Fact~\ref{fact:reg}(c).

\paragraph{Proof of Proposition~\ref{prop:weights}.} For $u\in\Sch$, using $x_jv_j(p)=v_j(p)x_j+i(\partial_jv_j)(p)$ and \eqref{eq:A},
\[
A\jap{x}^{-1}u=v(p)\cdot\frac{x}{\jap{x}}\,u+\frac i2(\dv v)(p)\jap{x}^{-1}u=:Bu,
\]
and $B$ is bounded. Since $\jap{x}^{-1}\Sch\subset\Sch$ and $A$ is closed, approximating $h\in\HH$ by $h_n\in\Sch$ shows $\jap{x}^{-1}h\in \mathrm{Dom}(A)$ and $A\jap{x}^{-1}h=Bh$. Hence $\mathrm{Dom}(\jap{x})\subset \mathrm{Dom}(\jap{A})$ and $\|\jap{A}u\|\le C\|\jap{x}u\|$ with $C^2=1+\|B\|^2$, because $\|\jap{A}w\|^2=\|w\|^2+\|Aw\|^2$. Fact~\ref{fact:heinz} with $T_1=\jap{A}$, $T_2=C\jap{x}$ and $\alpha=s$ gives $\|\jap{A}^su\|\le C^s\|\jap{x}^su\|$ on $\mathrm{Dom}(\jap{x}^s)$, i.e., $\jap{A}^s\jap{x}^{-s}$ is bounded. Finally, for $u,w\in\HH$, $\ip{w}{\jap{x}^{-s}u}=\ip{\jap{A}^s\jap{x}^{-s}w}{\jap{A}^{-s}u}$, and taking the supremum over $\|w\|=1$ gives the estimate.

\paragraph{Proof of Proposition~\ref{prop:comm-a}.} All operators below map $\Sch$ into $\Sch$, so the manipulations are legitimate on $\Sch$.

\emph{Step 1.} By \eqref{eq:A}, and since $x_j$ commutes with $a$,
\[
[A,a]=\sum_jx_j[v_j(p),a]-\tfrac i2[(\dv v)(p),a].
\]

\emph{Step 2.} Using $[x_j,v_j(p)]=i(\partial_jv_j)(p)$,
\[
[v_j(p),x_ja]=x_j[v_j(p),a]+[v_j(p),x_j]\,a=x_j[v_j(p),a]-i(\partial_jv_j)(p)\,a,
\]
so $x_j[v_j(p),a]=[v_j(p),h_j]+i(\partial_jv_j)(p)a$. Summing over $j$ and inserting in Step 1,
\begin{multline*}
[A,a]=\sum_j[v_j(p),h_j]+i(\dv v)(p)a-\tfrac i2\big((\dv v)(p)a-a(\dv v)(p)\big)\\=\sum_j[v_j(p),h_j]+\tfrac i2\{(\dv v)(p),a\}.
\end{multline*}

\emph{Step 3.} For smooth $h$ with polynomially bounded derivatives, $v_j(p)=p_jR$ and
\[
[p_jR,h]=[p_j,h]R+p_j[R,h],\quad [p_j,h]=-i\partial_jh,\quad [R,h]=R[h,1+p^2]R=R[h,p^2]R,
\]
\[
[h,p^2]=-\sum_l\big(p_l[p_l,h]+[p_l,h]p_l\big)=i\sum_l\big(p_l(\partial_lh)+(\partial_lh)p_l\big).
\]
Hence $[v_j(p),h]=-i(\partial_jh)R+ip_jR\sum_l\big(p_l(\partial_lh)+(\partial_lh)p_l\big)R$. With $h=h_j$ and $p_jRp_l=p_jp_lR$ this is the bracket in \eqref{eq:comm-a}.

\paragraph{Proof of Proposition~\ref{prop:algebra}.} (a) is clear, since $\overline b\in\Mfr$ and $\overline a\in\Ffr$. (b) By the Leibniz rule on $\Sch$ it suffices to treat generators. For $b\in\Mfr$, \eqref{eq:multcomm} gives $i(v\cdot\nabla b)(p)$ with $v\cdot\nabla b\in\Mfr$. For $a\in\Ffr$, \eqref{eq:comm-a} is a finite sum of products of $\dv v,\ R,\ p_jp_lR,\ p_jR\in\Mfr$ and $a,\ \partial_lh_j\in\Ffr$. (c) By (a), (b) and Proposition~\ref{prop:core}, $P\in C^1(A)$ with $[A,P]=\overline{P'}\in\Acal$; iterate and use Fact~\ref{fact:reg}(c).

\paragraph{Sobolev and weighted estimates.} We use the Weyl calculus of \cite[Ch.~7]{S:DS99} with Planck constant $h=1$, and its functional calculus \cite[Ch.~8]{S:DS99}, where a small fixed $h$ is used once (Lemma~\ref{lem:fcalc}). Page numbers refer to \cite{S:DS99}.

\emph{Quantization and symbol classes.} For $a\in\Sch'(\R^{2d})$ let $\Op(a)$ be the Weyl quantization \cite[(7.1)]{S:DS99} with $h=1$, $t=\frac12$:
\[
\Op(a)u(x)=(2\pi)^{-d}\iint e^{i(x-y)\cdot\xi}\,a\big(\tfrac{x+y}2,\xi\big)\,u(y)\,dy\,d\xi .
\]
If $a=a(x)$, then $\Op(a)$ is multiplication by $a$; if $a=b(\xi)$, then $\Op(a)=b(p)$. An \emph{order function} is a function $M:\R^{2d}\to(0,\infty)$ with $M(X)\le C_0\jap{X-Y}^{N_0}M(Y)$ \cite[Def.~7.4]{S:DS99}, and $S(M)$ is the Fr\'echet space of $a\in C^\infty(\R^{2d})$ with $|\partial^\gamma a|\le C_\gamma M$ for all $\gamma$ \cite[Def.~7.5]{S:DS99}. By Peetre's inequality $\jap{x+y}^s\le2^{|s|/2}\jap{x}^{|s|}\jap{y}^s$, every $M=\jap{x}^a\jap{\xi}^b$ with $a,b\in\R$ is an order function. If $M\le CM'$ then $S(M)\subset S(M')$ continuously. We use four results.

\textup{($\Psi$1)} \emph{Composition} \cite[Prop.~7.7, Lemma~7.8, Thm.~7.9, pp.~82--83]{S:DS99}. For $a\in S(M)$, $\Op(a)$ is continuous $\Sch\to\Sch$. For $a_j\in S(M_j)$, $\Op(a_1)\Op(a_2)=\Op(a_1\#a_2)$, where $\#$ is a continuous bilinear map $S(M_1)\times S(M_2)\to S(M_1M_2)$. In particular $\#$ maps bounded sets to bounded sets. \textup{($\Psi$2)} \emph{$L^2$-boundedness} \cite[Thm.~7.11, p.~85]{S:DS99}. For $a\in S(1)$, $\Op(a)$ extends to a bounded operator on $\HH$, whose norm is bounded by a seminorm of $a$ in $S(1)$ (the proof reduces to $h=1$ and gives this bound). \textup{($\Psi$3)} \emph{Exact Egorov theorem for the free flow} \cite[Ch.~7, Appendix, Thm.~A.2, p.~91]{S:DS99}. For $\tau\in\R$ and $a\in S(M)$, on $\Sch$,
\[
\Op(a)\,e^{-i\tau H_0}=e^{-i\tau H_0}\,\Op(a\circ\kappa_\tau),\qquad\kappa_\tau(x,\xi)=(x+2\tau\xi,\xi).
\]
For $\tau\ne0$, $e^{-i\tau H_0}$ has the kernel $(4\pi i\tau)^{-d/2}e^{i\Phi}$ with $\Phi(x,y)=|x-y|^2/(4\tau)$. This quadratic phase has no fibre variables, and $\{(x,\partial_x\Phi;y,-\partial_y\Phi)\}$ is the graph of $\kappa_\tau$. Hence $e^{-i\tau H_0}$ is a metaplectic operator associated with $\kappa_\tau$ in the sense of \cite[p.~90]{S:DS99}, and Theorem~A.2 there applies. As a check, for $a=x_j$ the formula reads $x_je^{-i\tau H_0}=e^{-i\tau H_0}(x_j+2\tau p_j)$. \textup{($\Psi$4)} \emph{Functional calculus} \cite[Thm.~8.1, Prop.~8.6, Thm.~8.7, pp.~94--102]{S:DS99}, in the following form.

\begin{klemma}\label{lem:fcalc}
Assume \eqref{ass:V}. Let $\mathcal F\subset C_c^\infty(\R)$ be bounded in $C^\infty$, with all supports in a fixed compact set. Then for every $N\in\N$ there are symbols $a_f\in S(\jap{\xi}^{-N})$, $f\in\mathcal F$, forming a bounded subset of $S(\jap{\xi}^{-N})$, with $f(H)=\Op(a_f)$.
\end{klemma}

\begin{proof}
\emph{Semiclassical rescaling.} For $h\in(0,1]$ put $p_h(x,\xi)=|\xi|^2+h^2V(x)$ and $m(x,\xi)=\jap{\xi}^2$. Then $p_h$ is real, $\{p_h\}_{h\in(0,1]}$ is bounded in $S(m)$, and $|p_h+i|\ge m/C$ with $C=2\|V\|_\infty+4$. Hence $p_h$ satisfies the standing assumptions of \cite[pp.~100--101]{S:DS99}. Let $h_0$ be the constant given there and fix $h_*\in(0,h_0]$. The Weyl $h_*$-quantization of $p_{h_*}$ is $P:=-h_*^2\Delta_x+h_*^2V=h_*^2H$. By \cite[Prop.~8.5]{S:DS99}, $P|_\Sch$ is essentially self-adjoint. Since $h_*^2H$ on $\mathcal H^2$ is a self-adjoint extension of $P|_\Sch$, the self-adjoint operator of \cite[Ch.~8]{S:DS99} is $P=h_*^2H$.

\emph{The symbol of $f(P)$.} By \cite[Thm.~8.7]{S:DS99}, for $f\in C_c^\infty(\R)$ and $k\in\N$, $f(P)=\Op_{h_*}(b_f)$ with $b_f\in S(m^{-k})$. We need these symbols to be bounded when $f$ ranges over a bounded set $\mathcal F'$ with supports in a fixed compact set. The proof of \cite[Thm.~8.7]{S:DS99} gives this. In \cite[Thm.~8.1]{S:DS99} any compactly supported $C^1$ extension $\tilde f$ with $\bar\partial\tilde f=O(|\operatorname{Im}z|)$ may be used. Take the finite-order extension $\tilde f_N(x+iy)=\chi(y)\sum_{j\le N}f^{(j)}(x)(iy)^j/j!$, with $\chi\in C_c^\infty(\R)$ equal to $1$ near $0$. Then $|\bar\partial\tilde f_N|\le C_N\sum_{j\le N+1}\|f^{(j)}\|_\infty\,|y|^N$, and the support of $\tilde f_N$ lies in a fixed compact set. Inserting this in \cite[(8.13)]{S:DS99} and using the resolvent symbol bounds \cite[(8.12)]{S:DS99} with $N\ge|\alpha|+|\beta|+2d+2$ shows that each seminorm of $b_f$ in $S(1)$ is bounded by a finite $C^j$-norm of $f$. The passage to $S(m^{-k})$ in that proof writes $f(P)=(P+i)^{-k}f_k(P)$ with $f_k(t)=(t+i)^kf(t)$ and $(P+i)^{-k}\in\Op_{h_*}(S(m^{-k}))$ \cite[p.~101]{S:DS99}. Here $\{f_k:f\in\mathcal F'\}$ is again bounded with supports in the same compact set, and ($\Psi$1) applies.

\emph{Back to $h=1$.} For $f\in\mathcal F$ put $f_*(t)=f(t/h_*^2)$, so that $f(H)=f_*(P)$, and $\{f_*\}$ is bounded with supports in a fixed compact set. The substitution $\eta=h_*\xi$ in \cite[(7.1)]{S:DS99} gives $\Op_{h_*}(b)=\Op(b(\cdot,h_*\,\cdot))$. Since $\jap{h_*\xi}^{-1}\le h_*^{-1}\jap{\xi}^{-1}$ and $\xi$-derivatives of $b(x,h_*\xi)$ carry factors $h_*$, the map $b\mapsto b(\cdot,h_*\,\cdot)$ is continuous $S(\jap{\xi}^{-2k})\to S(\jap{\xi}^{-2k})$. Take $2k\ge N$.
\end{proof}

\paragraph{Proof of Fact~\ref{fact:sobolev}.} Commuting derivatives with the equation and using energy estimates gives bounds on $\|e^{-itH}w\|_{\mathcal H^k}$ that grow with $|t|$. Uniformity in $t$ comes instead from the conserved quantities $\|(H+c)^je^{-itH}w\|=\|(H+c)^jw\|$, which we compare with Sobolev norms. Fix $c>\|V\|_\infty+1$, so that $H+c\ge\delta:=c-\|V\|_\infty>1$, and use the norm $\|w\|_{\mathcal H^k}=\|\jap{p}^kw\|$.

\emph{Claim.} For $j\in\N$, $\mathrm{Dom}((H+c)^j)=\mathcal H^{2j}$, and $C_j^{-1}\|\jap{p}^{2j}w\|\le\|(H+c)^jw\|\le C_j\|\jap{p}^{2j}w\|$ on $\mathcal H^{2j}$.

For $j=1$: $\mathrm{Dom}(H)=\mathcal H^2$ by the Kato--Rellich theorem, $\|(H+c)w\|\le\|(H_0+1)w\|+(c+\|V\|_\infty)\|w\|$, and $\|(H_0+1)w\|\le\|(H+c)w\|+(c+\|V\|_\infty)\|w\|$ with $\|w\|\le\delta^{-1}\|(H+c)w\|$.

Assume the claim for $j$. If $w\in \mathrm{Dom}((H+c)^{j+1})$, then $u:=(H+c)w\in \mathrm{Dom}((H+c)^j)=\mathcal H^{2j}$ and $w\in \mathcal H^{2j}$, so $(H_0+1)w=u-(V+c-1)w\in \mathcal H^{2j}$ by the Leibniz rule, i.e., $w\in \mathcal H^{2j+2}$. Conversely, if $w\in \mathcal H^{2j+2}$ then $u\in \mathcal H^{2j}=\mathrm{Dom}((H+c)^j)$. For the norms, the Leibniz rule gives $\|\jap{p}^{2j}(V+c-1)w\|\le C\|\jap{p}^{2j}w\|$. Hence
\[
\|\jap{p}^{2j+2}w\|\le\|\jap{p}^{2j}u\|+C\|\jap{p}^{2j}w\|\le C_j\|(H+c)^{j+1}w\|+CC_j\|(H+c)^jw\|\le C'\|(H+c)^{j+1}w\|,
\]
using $\|(H+c)^jw\|\le\delta^{-1}\|(H+c)^{j+1}w\|$. The reverse inequality follows in the same way from $\|(H+c)^{j+1}w\|=\|(H+c)^ju\|\le C_j\|\jap{p}^{2j}u\|$.

\emph{All $k$.} Let $k\le 2j$. Apply Fact~\ref{fact:heinz} with $\alpha=k/(2j)$, once to $T_1=\jap{p}^{2j}$, $T_2=C_j(H+c)^j$ and once with the roles exchanged. Both are positive self-adjoint with the common domain $\mathcal H^{2j}$. This gives $\mathrm{Dom}((H+c)^{k/2})=\mathcal H^k$ and $\|\jap{p}^kw\|\simeq\|(H+c)^{k/2}w\|$. Since $e^{-itH}$ is unitary and commutes with $(H+c)^{k/2}$, $\|e^{-itH}w\|_{\mathcal H^k}\le C\|(H+c)^{k/2}w\|\le C'\|w\|_{\mathcal H^k}$ uniformly in $t$.

\emph{Kicks and free flows.} $\partial^\gamma e^{-icV}$ is $e^{-icV}$ times a polynomial in $c$ and derivatives of $V$, so $|\partial^\gamma e^{-icV}|\le C_\gamma\jap{c}^{|\gamma|}$, and the Leibniz rule gives the bound for integer $k$. Free flows commute with $\jap{p}^k$.

\paragraph{Proof of Fact~\ref{fact:weights}.} All operators below map $\Sch$ to $\Sch$ (by ($\Psi$1) and since free flows and kicks preserve $\Sch$). We prove the bounds on $\Sch$; by density they define bounded operators on $\HH$.

(a) \emph{Free flow.} By ($\Psi$3) with $a=\jap{x}^{-\mu}$, $\jap{x}^{-\mu}e^{-i\tau H_0}=e^{-i\tau H_0}\Op(w_\tau)$ with $w_\tau(x,\xi)=\jap{x+2\tau\xi}^{-\mu}$. Here $\partial_x^\alpha\partial_\xi^\beta w_\tau=(2\tau)^\beta(\partial^{\alpha+\beta}\jap{\cdot}^{-\mu})(x+2\tau\xi)$, $|\partial^\gamma\jap{y}^{-\mu}|\le C_\gamma\jap{y}^{-\mu}$, and by Peetre's inequality $\jap{x+2\tau\xi}^{-\mu}\le2^{\mu/2}\jap{x}^{-\mu}\jap{2\tau\xi}^\mu$. Hence $\{w_\tau\}_{|\tau|\le T_0}$ is bounded in $S(\jap{x}^{-\mu}\jap{\xi}^\mu)$. Since $\jap{p}^k$ commutes with $e^{-i\tau H_0}$,
\[
\jap{p}^k\jap{x}^{-\mu}e^{-i\tau H_0}\jap{x}^\mu\jap{p}^{-k-\mu}=e^{-i\tau H_0}\,\Op(\jap{\xi}^k)\Op(w_\tau)\Op(\jap{x}^\mu)\Op(\jap{\xi}^{-k-\mu}).
\]
By ($\Psi$1) the product is $\Op$ of a bounded family in $S(\jap{\xi}^k\cdot\jap{x}^{-\mu}\jap{\xi}^\mu\cdot\jap{x}^\mu\cdot\jap{\xi}^{-k-\mu})=S(1)$, and ($\Psi$2) gives the bound. The factor $\jap{\xi}^\mu$ in the order of $w_\tau$ is the loss of $\mu$ derivatives.

\emph{Kicks.} $e^{-icV}=\Op(e^{-icV(x)})$, and $\{e^{-icV}\}_{|c|\le T_0}$ is bounded in $S(1)$ by the bound on $\partial^\gamma e^{-icV}$ in the proof of Fact~\ref{fact:sobolev}. The product has order $\jap{\xi}^k\cdot1\cdot\jap{\xi}^{-k}=1$.

\emph{Differential operators.} $Q=\sum_{|\beta|\le\ell}c_\beta(x)\partial^\beta=\sum_\beta\Op(c_\beta)\Op((i\xi)^\beta)$, so $Q=\Op(a_Q)$ with $a_Q\in S(\jap{\xi}^\ell)$ by ($\Psi$1). The product has order $\jap{\xi}^k\jap{x}^{-\mu}\jap{\xi}^\ell\jap{x}^\mu\jap{\xi}^{-k-\ell}=1$.

\emph{Consequence \eqref{eq:keyweight}.} Let $Q=\sum_{|\beta|\le\ell}c_\beta\partial^\beta$ with $|\partial^\alpha c_\beta|\le C\jap{x}^{-\mu}$. Then $Q\jap{x}^\mu=\sum_{\gamma\le\beta}\binom\beta\gamma c_\beta(\partial^{\beta-\gamma}\jap{x}^\mu)\partial^\gamma$ has $C_b^\infty$ coefficients. Since $e^{-icV}$ commutes with $\jap{x}^\mu$,
\begin{multline*}
Qe^{-i\tau H_0}e^{-icV}\jap{x}^\mu\jap{p}^{-k}=\big(Q\jap{x}^\mu\jap{p}^{-\ell}\big)\big(\jap{p}^\ell\jap{x}^{-\mu}e^{-i\tau H_0}\jap{x}^\mu\jap{p}^{-\ell-\mu}\big)\\\big(\jap{p}^{\ell+\mu}e^{-icV}\jap{p}^{-\ell-\mu}\big)\jap{p}^{\ell+\mu-k},
\end{multline*}
and each factor is bounded uniformly for $|\tau|,|c|\le T_0$ when $k\ge \mu+\ell$.

(b) By Lemma~\ref{lem:fcalc} with $\mathcal F=\{\tilde g\}$ and $N\ge k$, $\tilde g(H)=\Op(a)$ with $a\in S(\jap{\xi}^{-N})$. Then $\jap{p}^k\jap{x}^{-\mu}\tilde g(H)\jap{x}^\mu$ is $\Op$ of a symbol in $S(\jap{\xi}^k\jap{x}^{-\mu}\jap{\xi}^{-N}\jap{x}^\mu)=S(\jap{\xi}^{k-N})\subset S(1)$ by ($\Psi$1), and ($\Psi$2) applies. The bound is uniform for $\tilde g$ in any family $\mathcal F$ as in Lemma~\ref{lem:fcalc}.

(c) Put $f_\tau(\lambda)=e^{i\tau\lambda}\tilde g(\lambda)$. Then $e^{i\tau H}\tilde g(H)=f_\tau(H)$, $\supp f_\tau=\supp\tilde g$, and $\|f_\tau^{(j)}\|_\infty\le\sum_{i\le j}\binom ji|\tau|^i\|\tilde g^{(j-i)}\|_\infty$. So $\{f_\tau\}_{|\tau|\le T_0}$ is a family as in Lemma~\ref{lem:fcalc}, and (c) is the uniform version of (b) with $\mu$ replaced by $-\mu$ and $k=0$.

The only inputs from \cite{S:DS99} are the Weyl calculus at $h=1$ (($\Psi$1)--($\Psi$3)) and the functional calculus at one fixed small $h$ (Lemma~\ref{lem:fcalc}). The uniformity in Lemma~\ref{lem:fcalc} is not stated in \cite[Thm.~8.7]{S:DS99}, but it follows from its proof as indicated. Uniformity in $\tilde g$ is used only in (c).

%=====================================================================
\section{Index of notation}\label{sec:Snotation}
%=====================================================================

The index covers the paper and this supplement. The last column gives the place where each symbol is defined; references without the prefix S are to the paper. Letters that also have a local meaning, distinguished by context, are listed in the second table.

{\small
\begin{longtable}{p{0.22\textwidth}p{0.53\textwidth}p{0.17\textwidth}}
\toprule
Symbol & Meaning & Defined in\\\midrule
\endhead
\multicolumn{3}{l}{\emph{General notation}}\\
$L^2$, $x$, $p$, $\jap{\xi}$ & $L^2(\R^d)$, position, momentum $p=-i\nabla$, $(1+|\xi|^2)^{1/2}$ & section~\ref{sec:intro}\\
$\Sch$, $\mathcal H^k$, $\mathcal H^\infty$ & Schwartz space, Sobolev spaces, $\cap_k\mathcal H^k$ & section~\ref{sec:intro}\\
$b(p)$, $\hat u$ & Fourier multiplier, Fourier transform $\int e^{-ix\cdot\xi}u\,dx$ & section~\ref{sec:intro}\\
$\BB$, $\KK$, $\mathscr B(X,Y)$ & bounded and compact operators, bounded operators $X\to Y$ & section~\ref{sec:intro}\\
$\T$, $C_0$, $C_b^\infty$, $\one_B$ & unit circle, continuous functions vanishing at infinity, smooth functions with bounded derivatives, indicator of $B$ & section~\ref{sec:intro}\\
$E^T$, $E^T_{\rm sc}$ & spectral measure of $T$ and its singular continuous part & section~\ref{sec:intro}\\
$P_{\rm ac}$, $P_{\rm sc}$, $P_{\rm pp}$, $L^2_{\rm ac}$, \dots & spectral projections and subspaces & section~\ref{sec:intro}\\
$\Delta_x$, $\Delta$ & Laplacian; time step & section~\ref{sec:intro}\\
\multicolumn{3}{l}{\emph{Setting, methods and scattering operators}}\\
$H_0$, $H$, $V$, $\rho$ & $-\Delta_x$, $H_0+V$, potential, decay exponent & \eqref{ass:V}\\
$\Omega^c_\pm$, $W^c_\pm$, $S$ & wave operators, their inverses and scattering operator of $(H,H_0)$ & section~\ref{sec:setting}\\
$U$, $U(s)$, $U_0$ & one-step operator, step $s$, free step $e^{-i\Delta H_0}$ & \eqref{eq:comp}\\
$F_j$, $K_j$, $a_j$, $b_j$ & free flows, kicks and coefficients of a composition & \eqref{eq:comp}\\
$K(a)$, $\nu_a$ & kick $e^{-ia\Delta V}$ of strength $a$, $K(a)-1$ & \eqref{eq:UminusU0}, Lemma~\ref{lem:kick}\\
$U_{\rm L}$, $U_{\rm R}$, $U_{\rm S}$ & Lie, reversed Lie and Strang steps & section~\ref{sec:setting}\\
$r$, $q$ & order and effective order of a method & Def.~\ref{def:order}\\
$P$, $\tilde U$, $C_P$, $\Pi$ & processor, kernel $U=P\tilde UP^{-1}$, total free duration of $P$ over $\Delta$, $e^{iC_P\Delta H_0}$ & Def.~\ref{def:order}, Thm.~\ref{thm:B}\\
$\varepsilon_\Delta$, $\Theta_{[a,b]}$ & aliasing map $E\mapsto e^{-i\Delta E}$, arcs & section~\ref{sec:setting}\\
$\Omega^U_\pm$, $W^U_\pm$, $S_U$, $\Omega_\pm(\Delta)$ & splitting wave operators, their inverses, splitting scattering operator; $\Omega^U_\pm$ for Lie & section~\ref{sec:setting}\\
$g$, $s$, $\chi$, $J_\chi$ & energy cut-off and weight of the data $g(H)\jap{x}^{-s}$; asymptotic data $\chi(p)\jap{x}^{-\rho}$ & section~\ref{sec:setting}\\
$\Omega_\pm(U,U_1)$ & wave operators of a pair of unitary operators & section~\ref{sec:Sbg-scat}\\
$\Xi_\pm$, $\Sigma_U$ & relative wave operators and relative scattering operator & Thms.~\ref{thm:conv}(d), \ref{thm:S}(a)\\
\multicolumn{3}{l}{\emph{Scattering theory of the scheme (section~\ref{sec:scattering}, appendices)}}\\
$T_i$ & factors of $U=T_m\cdots T_1$ & section~\ref{sec:scattering}\\
$X_j$, $\nu_j$, $Y_j$, $U_{[j]}$ & terms of the telescoping identity & \eqref{eq:telescope-comp}\\
$A$, $v$, $\gamma$, $R$ & conjugate operator, $v(\xi)=\xi\jap{\xi}^{-2}$, $\gamma(\xi)=|\xi|^2\jap{\xi}^{-2}$, $R=\jap{p}^{-2}$ & section~\ref{sec:kick}\\
$\Mfr$, $\Ffr$, $\Acal$ & symbols of order $0$ in $\xi$, decaying symbols in $x$, the algebra they generate & section~\ref{sec:kick}, appendix~\ref{sec:conjugate}\\
$\phi_t$, $J_t$, $T_t$ & flow of $v$, its Jacobian, the group $e^{-itA}$ & Prop.~\ref{prop:group}\\
$C^k(A)$, $[A,S]$ & regularity classes and commutators & appendix~\ref{sec:abstract}\\
$\mathsf K$, $\mathsf K_1$, $c_0$ & compact remainders and constant in a Mourre estimate & \eqref{eq:mourre}, Lemma~\ref{lem:product}\\
$\mathrm{Dom}(\jap{A}^s)^*$ & completion of $L^2$ for $\|\jap{A}^{-s}u\|$ & appendix~\ref{sec:abstract}\\
$\mathcal O$, $\nu_{i,s}$ & set $\T\setminus(\{1\}\cup\sigma_p(U))$, weighted kicks $\jap{x}^{2s}\nu_i$ & section~\ref{sec:proofA}\\
\multicolumn{3}{l}{\emph{Convergence theory (section~\ref{sec:convergence})}}\\
$\mu$, $k$, $\epsilon$ & weight exponent, Sobolev order and loss in local decay & Prop.~\ref{prop:localdecay}\\
$R_+$ & remainder of outgoing scattering states & Lemma~\ref{lem:scattstates}\\
$D(s)$, $D_S(s)$, $C_a$ & one-step defects, $C_a=[H_0,[H_0,a]]$ & Lemmas~\ref{lem:defect}, \ref{lem:defect-comp}\\
$L_i$, $Z_i$, $c_i$, $c_{\max}$, $\mathrm{Ad}$ & partial products, generators and coefficients of a composition, $\mathrm{Ad}(T)X=TXT^{-1}$ & section~\ref{sec:defects}\\
\multicolumn{3}{l}{\emph{Bound states}}\\
$\phi_0$, $E_0$ & bound state of $H$ and its energy & section~\ref{sec:bound-results}\\
$e_k$, $E_N(\Delta)$ & coefficients of the quasi-energy and its truncation & Thm.~\ref{thm:meta}\\
$\hat H$, $h_k$, $\chi_k$, $\chi^N$, $\lambda_N$, $A_k$ & modified Hamiltonian, its terms, quasi-modes, $e^{-i\Delta E_N}$, Taylor coefficients of $U$ (not the conjugate operator $A$) & proof of Thm.~\ref{thm:meta}\\
$M_\alpha$ & exponential moment of $\hat V$ & \eqref{ass:A}\\
$\zeta$, $F$, $\rho_m$ & eigenvalue, Fourier transform of $(e^{-i\Delta V}-1)w$, radii of the resonant spheres & proof of Thm.~\ref{thm:noeig}\\
$\Gamma_m$, $\Gamma_{\rm FGR}$, $\Gamma$ & golden-rule rates and measured ionization rate & \eqref{eq:Sfgr}\\
$U_{\rm Y}$, $S(h)$, $w_1$, $w_0$ & Yoshida triple jump, Strang step with step $h$, its weights & section~\ref{sec:Snumerics}\\
$E_c$, in-band and aliased parts & mean energy $\ip{\psi}{H\psi}$ of a packet; parts of an error below and above the $H_0$-energy $E_c+\pi/\Delta$ & section~\ref{sec:Snumerics}\\
$\hat\kappa_m$ & Fourier coefficients of the kick comb & \eqref{eq:fgr-comp}\\
$\Pi_K$, $U_K$, $H_K$, $\tilde H_K$, $\tilde\psi$, $\tilde E$ & momentum cut-off, compressed step and Hamiltonian, logarithm of $U_K$, quasi-mode and quasi-energy & section~\ref{sec:Sexpmeta}\\
$\beta_j\Delta$ & kick times of a composition, $\beta_j=b_1+\dots+b_{j-1}$ & section~\ref{sec:setting}\\
$r$, $\varrho$, $\lambda_0$ & residual of an approximate eigenvector, its minimum over $\lambda$, the minimizing phase & section~\ref{sec:Srate}\\
$y$, $\eta$, $\tau$, $U_1$, $\hat U$, $\varsigma$, $q$, $\mathcal M$ & $U^*r/\|r\|$, $\ip{y}{\tilde\psi}$, the unimodular $1-\|r\|/\eta$, rank-one decoupled scheme, its restriction, coefficient, coupling vector, Cauchy-type transform & section~\ref{sec:Srate}\\
$P^T_f$ & Poisson integral of the spectral measure of $f$ for $T$ & section~\ref{sec:Srate}\\
$\gamma_0$, $\Gamma_r$, (NC$_\Delta$) & lower bound, golden rule for the residual, non-cancellation hypothesis & section~\ref{sec:Srate}\\
\multicolumn{3}{l}{\emph{Long-range potentials and elliptic operators (supplement)}}\\
$G$ & the data operator $g(H)\jap{x}^{-s}$ & section~\ref{sec:Slongrange}\\
$\Lambda_\pm$, $\Omega^c_{\Lambda,\pm}$, $W^c_{\Lambda,\pm}$, $S^c_\Lambda$ & modifiers, modified wave and scattering operators of $H$ & Hyp.~(M)\\
$\Omega^U_{\Lambda,\pm}$, $W^U_{\Lambda,\pm}$, $S^U_\Lambda$ & modified splitting wave operators, their adjoints, modified splitting scattering operator & Thm.~\ref{thm:modWO}\\
$\Phi$, $\Phi^U_n$, $c_\Delta$ & Dollard phase, scheme-adapted phase, limiting phase difference & Hyp.~(M), section~\ref{sec:Sextensions}\\
$\omega$, $\mathcal T(\omega)$, $\Theta_\omega$, $\gamma_\omega$ & symbol of $H_0=\omega(p)$, its critical values, their images under $\varepsilon_\Delta$, commutator symbol & section~\ref{sec:Sscalar}\\
$SG^{\mu,\nu}$, $\Op$, $S(M)$ & SG symbol classes, Weyl quantization, symbol classes & sections~\ref{sec:Sscalar}, \ref{sec:Stools}\\\bottomrule
\end{longtable}
}

{\small
\begin{longtable}{p{0.08\textwidth}p{0.86\textwidth}}
\toprule
Letter & Meanings\\\midrule
\endhead
$a$ & strength of a kick; $a=V/2$ in Lemma~\ref{lem:defect}(b) and section~\ref{sec:Sexpmeta}; functions in $\Ffr$; endpoints of arcs; Gaussian rate in Theorem~\ref{thm:noeig}\\
$A$ & conjugate operator; Taylor coefficients $A_k$ in Theorem~\ref{thm:meta}; $A^\pm_j$ in Proposition~\ref{prop:Snoeig}\\
$F$ & free flows $F_j$, $F_\tau$; Fourier transform of $(e^{-i\Delta V}-1)w$ in Theorem~\ref{thm:noeig}\\
$g$ & energy cut-off; auxiliary functions in Theorem~\ref{thm:noeig}, Corollary~\ref{cor:ion} and Proposition~\ref{prop:Snoeig}\\
$h$ & $h_j=x_ja$ in Proposition~\ref{prop:comm-a}; terms $h_k$ of the modified Hamiltonian\\
$k$ & Sobolev order; momenta $k_0$, $k_m$\\
$K$ & kick $K(a)$; momentum cut-off in section~\ref{sec:Sexpmeta}\\
$m$ & index of harmonics, written $\ell$ in section~\ref{sec:Sscalar}; degree $2m$ of $\omega$ there; an integer in Proposition~\ref{prop:localdecay}; an order function in Lemma~\ref{lem:fcalc}\\
$P$ & processor; with indices or decorations, projections, the survival probability $P_n$ and Poisson integrals $P^T_f$\\
$q$ & effective order; coupling vector in sections~\ref{sec:Srate} and~\ref{sec:Sevidence}\\
$Q$ & with indices, products of factors or cut-offs; a differential operator in Proposition~\ref{prop:defectbounds}; a projection in Theorem~\ref{thm:meta}\\
$r$ & order of a method; residual in section~\ref{sec:Srate}; radius in Jensen's formula in Theorem~\ref{thm:noeig}\\
$R$ & $\jap{p}^{-2}$; with indices, remainders and resolvents\\
$s$ & weight exponent; time variable within one step; number of stages in \eqref{eq:comp}\\
$v$ & vector field $v(\xi)$; the vector $(U-z)^{-1}\tilde\psi$ in section~\ref{sec:Srate}\\
$\beta$ & multi-indices; decay exponent in Theorem~\ref{thm:noeig}(b); kick times $\beta_j\Delta$\\
$\gamma_0$ & $E/(1+E)$ in Lemma~\ref{lem:free}; the lower bound in (NC$_\Delta$)\\
$\eta$ & cut-off $\eta(H_0)$; imaginary parts; a parameter of \cite{S:JMP84} in Proposition~\ref{prop:localdecay}; $\ip{y}{\tilde\psi}$ in section~\ref{sec:Srate}\\
$\rho$ & decay exponent; radii $\rho_m$ in Theorem~\ref{thm:noeig}; order $\rho_B$ of an entire function in Proposition~\ref{prop:Snoeig}\\
$\sigma$ & spectra; integration variables\\
$\Phi$ & Dollard phase; phase functions in stationary-phase arguments\\
$\chi$ & cut-off of the asymptotic data; quasi-modes $\chi_k$, $\chi^N$ in Theorem~\ref{thm:meta}\\
\bottomrule
\end{longtable}
}

%=====================================================================
\section{Code}\label{sec:Scode}
%=====================================================================

The Python scripts reproducing the figures, Table~\ref{tab:Sresid} and the numbers of sections~\ref{sec:Snumerics} and~\ref{sec:Sevidence} are provided as ancillary files with the arXiv version of this paper, in the directory \texttt{numerics}. They use \texttt{numpy}, \texttt{matplotlib} and, for one check, \texttt{scipy}; the script \texttt{run\_all.sh} runs them in order, and the computed data are stored in \texttt{data/}.

\emph{Modules.} \texttt{common.py} implements the periodic grids, the Fourier split-step compositions (Lie, Strang and the Yoshida triple jump, and their inverses), norms, momentum densities, and the split of an error into its in-band and aliased parts (section~\ref{sec:Snumerics}); \texttt{style.py} the plotting style; \texttt{fgr.py} the golden rules \eqref{eq:fgrPT} and \eqref{eq:fgr-comp}, with a quadrature check of the matrix element $\ip{\psi_k}{V\phi_0}$.

\emph{Scattering states (section~\ref{sec:Snumerics}).} \texttt{scat\_errors.py} (Figure~\ref{fig:errors} and the Yoshida errors), \texttt{lie\_prekick.py} (the Lie plateau removed by a half kick), \texttt{scat\_alias.py} (Figure~\ref{fig:alias}), \texttt{scat\_limits.py} (Figure~\ref{fig:rates}), \texttt{slow\_decay.py} (Figure~\ref{fig:slow}) and \texttt{scat\_energy.py} (Figure~\ref{fig:energy}) compute the data, and \texttt{fig\_errors.py}, \texttt{fig\_alias.py}, \texttt{fig\_limits.py}, \texttt{fig\_slow.py} and \texttt{fig\_energy.py} draw the figures.

\emph{Bound states (section~\ref{sec:Sevidence}).} \texttt{bound\_rates.py} computes the survival probabilities and ionization rates (Lie, Strang, Yoshida, and the truncated well), \texttt{spectrum.py} the momentum density of the ionized part and its line weights, \texttt{quasienergy.py} the quasi-energy fits, and \texttt{fig\_bound.py} draws Figures~\ref{fig:ionisation} and~\ref{fig:spectrum}. \texttt{residual\_rate.py} (arguments $\Delta$ and the cutoffs $K/k_1$) computes the residual golden rule of Table~\ref{tab:Sresid} and the approximations discussed after it: it diagonalizes the compressed step $U_K$ in the Fourier basis, forms $r$, $\hat U$, $\varsigma$ and $q$ as in section~\ref{sec:Srate}, evaluates the density of $\mu^{\hat U}_q$ from the correlations $\ip{q}{\hat U^nq}$, and samples $|\varsigma\mathcal M|$ on $\mathcal Z_\Delta$; \texttt{make\_table.py} assembles Table~\ref{tab:Sresid}.

\emph{Other.} \texttt{aliasing\_fig.py} draws Figure~\ref{fig:aliasing}. \texttt{check\_commutators.py} and \texttt{check\_strang.py} check numerically the commutator identities of Lemmas~\ref{lem:free}(a) and~\ref{lem:kick} and Proposition~\ref{prop:comm-a}, and the defect formulas \eqref{eq:Ca} and \eqref{eq:DS}, to relative accuracy $10^{-7}$ or better when exact derivatives of the potential are used.


{\small
\begin{thebibliography}{99}\addcontentsline{toc}{section}{References}
\bibitem{ABG96} W.~O.~Amrein, A.~Boutet de Monvel, V.~Georgescu, \emph{$C_0$-groups, commutator methods and spectral theory of $N$-body Hamiltonians}, Progr.\ Math.\ 135, Birkh\"auser, Basel, 1996.

\bibitem{ABCF06} M.~A.~Astaburuaga, O.~Bourget, V.~H.~Cort\'es, C.~Fern\'andez, Floquet operators without singular continuous spectrum, \emph{J.\ Funct.\ Anal.}\ 238 (2006), 489--517.

\bibitem{ABC15a} M.~A.~Astaburuaga, O.~Bourget, V.~H.~Cort\'es, Commutation relations for unitary operators I, \emph{J.\ Funct.\ Anal.}\ 268 (2015), 2188--2230.

\bibitem{BCF23} W.~Bao, Y.~Cai, Y.~Feng, Improved uniform error bounds of the time-splitting methods for the long-time (nonlinear) Schr\"odinger equation, \emph{Math.\ Comp.}\ 92 (2023), 1109--1139.

\bibitem{BGSvL26} S.~Becker, N.~Galke, R.~Salzmann, L.~van Luijk, Convergence rates for the Trotter splitting for unbounded operators, \emph{Found.\ Comput.\ Math.}\ 26 (2026), 2471--2524.

\bibitem{BCR99} S.~Blanes, F.~Casas, J.~Ros, Symplectic integration with processing: a general study, \emph{SIAM J.\ Sci.\ Comput.}\ 21 (1999), 711--727.

\bibitem{BC05} S.~Blanes, F.~Casas, On the necessity of negative coefficients for operator splitting schemes of order higher than two, \emph{Appl.\ Numer.\ Math.}\ 54 (2005), 23--37.

\bibitem{BFHJY24} D.~Burgarth, P.~Facchi, A.~Hahn, M.~Johnsson, K.~Yuasa, Strong error bounds for Trotter and Strang-splittings and their implications for quantum chemistry, \emph{Phys.\ Rev.\ Research} 6 (2024), 043155.

\bibitem{But69} J.~C.~Butcher, The effective order of Runge--Kutta methods, in \emph{Conference on the Numerical Solution of Differential Equations}, Lecture Notes in Math.\ 109, Springer, Berlin, 1969, 133--139.

\bibitem{CS24} R.~Carles, C.~Su, Scattering and uniform in time error estimates for splitting method in NLS, \emph{Found.\ Comput.\ Math.}\ 24 (2024), 683--722.

\bibitem{CM26} R.~Carles, G.~Maierhofer, On scattering for NLS: rigidity properties and numerical simulations via the lens transform, \emph{Math.\ Models Methods Appl.\ Sci.}\ (2026), doi:10.1142/S0218202526500302.

\bibitem{Car26} R.~Carles, Time splitting and error estimates for nonlinear Schr\"odinger equations with a potential, \emph{Found.\ Comput.\ Math.}\ 26 (2026), 2207--2258.

\bibitem{CLR00} O.~Costin, J.~L.~Lebowitz, A.~Rokhlenko, Exact results for the ionization of a model quantum system, \emph{J.\ Phys.\ A} 33 (2000), 6311--6319.

\bibitem{CCLR01} O.~Costin, R.~D.~Costin, J.~L.~Lebowitz, A.~Rokhlenko, Evolution of a model quantum system under time periodic forcing: conditions for complete ionization, \emph{Comm.\ Math.\ Phys.}\ 221 (2001), 1--26.

\bibitem{DG97} J.~Derezi\'nski, C.~G\'erard, \emph{Scattering Theory of Classical and Quantum $N$-Particle Systems}, Texts and Monographs in Physics, Springer, 1997.

\bibitem{DS99} M.~Dimassi, J.~Sj\"ostrand, \emph{Spectral Asymptotics in the Semi-Classical Limit}, London Math.\ Soc.\ Lecture Note Ser.\ 268, Cambridge University Press, 1999.

\bibitem{Dol64} J.~D.~Dollard, Asymptotic convergence and the Coulomb interaction, \emph{J.\ Math.\ Phys.}\ 5 (1964), 729--738.

\bibitem{DF07} G.~Dujardin, E.~Faou, Normal form and long time analysis of splitting schemes for the linear Schr\"odinger equation with small potential, \emph{Numer.\ Math.}\ 108 (2007), 223--262.

\bibitem{FGP10a} E.~Faou, B.~Gr\'ebert, E.~Paturel, Birkhoff normal form for splitting methods applied to semilinear Hamiltonian PDEs. Part I. Finite-dimensional discretization, \emph{Numer.\ Math.}\ 114 (2010), 429--458.

\bibitem{FJ15} E.~Faou, T.~J\'ez\'equel, Resonant time steps and instabilities in the numerical integration of Schr\"odinger equations, \emph{Differential Integral Equations} 28 (2015), 221--238.

\bibitem{FRT13} C.~Fern\'andez, S.~Richard, R.~Tiedra de Aldecoa, Commutator methods for unitary operators, \emph{J.\ Spectr.\ Theory} 3 (2013), 271--292.

\bibitem{FH82} R.~Froese, I.~Herbst, Exponential bounds and absence of positive eigenvalues for $N$-body Schr\"odinger operators, \emph{Comm.\ Math.\ Phys.}\ 87 (1982), 429--447.

\bibitem{HLW06} E.~Hairer, C.~Lubich, G.~Wanner, \emph{Geometric Numerical Integration}, 2nd ed., Springer, 2006.

\bibitem{Hei51} E.~Heinz, Beitr\"age zur St\"orungstheorie der Spektralzerlegung, \emph{Math.\ Ann.}\ 123 (1951), 415--438.

\bibitem{HHZ19} M.~Heyl, P.~Hauke, P.~Zoller, Quantum localization bounds Trotter errors in digital quantum simulation, \emph{Sci.\ Adv.}\ 5 (2019), eaau8342.

\bibitem{Hor76} L.~H\"ormander, The existence of wave operators in scattering theory, \emph{Math.\ Z.}\ 146 (1976), 69--91.

\bibitem{How79} J.~S.~Howland, Scattering theory for Hamiltonians periodic in time, \emph{Indiana Univ.\ Math.\ J.}\ 28 (1979), 471--494.

\bibitem{Ign11} L.~I.~Ignat, A splitting method for the nonlinear Schr\"odinger equation, \emph{J.\ Differential Equations} 250 (2011), 3022--3046.

\bibitem{JL00} T.~Jahnke, C.~Lubich, Error bounds for exponential operator splittings, \emph{BIT} 40 (2000), 735--744.

\bibitem{JMP84} A.~Jensen, E.~Mourre, P.~Perry, Multiple commutator estimates and resolvent smoothness in quantum scattering theory, \emph{Ann.\ Inst.\ H.\ Poincar\'e Phys.\ Th\'eor.}\ 41 (1984), no.~2, 207--225.

\bibitem{Kat66} T.~Kato, Wave operators and similarity for some non-selfadjoint operators, \emph{Math.\ Ann.}\ 162 (1966), 258--279.

\bibitem{Kat04} Y.~Katznelson, \emph{An Introduction to Harmonic Analysis}, 3rd ed., Cambridge University Press, Cambridge, 2004.

\bibitem{KW01} E.~Kirr, M.~I.~Weinstein, Parametrically excited Hamiltonian partial differential equations, \emph{SIAM J.\ Math.\ Anal.}\ 33 (2001), 16--52.

\bibitem{Lev96} B.~Ya.~Levin, \emph{Lectures on Entire Functions}, Transl.\ Math.\ Monogr.\ 150, AMS, Providence, 1996.

\bibitem{Lub08} C.~Lubich, On splitting methods for Schr\"odinger--Poisson and cubic nonlinear Schr\"odinger equations, \emph{Math.\ Comp.}\ 77 (2008), 2141--2153.

\bibitem{MQ02} R.~I.~McLachlan, G.~R.~W.~Quispel, Splitting methods, \emph{Acta Numer.}\ 11 (2002), 341--434.

\bibitem{Mou81} E.~Mourre, Absence of singular continuous spectrum for certain self-adjoint operators, \emph{Comm.\ Math.\ Phys.}\ 78 (1981), 391--408.

\bibitem{RS79} M.~Reed, B.~Simon, \emph{Methods of Modern Mathematical Physics III: Scattering Theory}, Academic Press, 1979.

\bibitem{RST19} S.~Richard, A.~Suzuki, R.~Tiedra de Aldecoa, Quantum walks with an anisotropic coin II: scattering theory, \emph{Lett.\ Math.\ Phys.}\ 109 (2019), 61--88.

\bibitem{SW98} A.~Soffer, M.~I.~Weinstein, Nonautonomous Hamiltonians, \emph{J.\ Stat.\ Phys.}\ 93 (1998), 359--391.

\bibitem{SW99} A.~Soffer, M.~I.~Weinstein, Resonances, radiation damping and instability in Hamiltonian nonlinear wave equations, \emph{Invent.\ Math.}\ 136 (1999), 9--74.

\bibitem{Tie20} R.~Tiedra de Aldecoa, Stationary scattering theory for unitary operators with an application to quantum walks, \emph{J.\ Funct.\ Anal.}\ 279 (2020), 108704.

\bibitem{Yaf92} D.~R.~Yafaev, \emph{Mathematical Scattering Theory: General Theory}, Transl.\ Math.\ Monogr.\ 105, AMS, Providence, 1992.

\bibitem{Yaj77} K.~Yajima, Scattering theory for Schr\"odinger equations with potentials periodic in time, \emph{J.\ Math.\ Soc.\ Japan} 29 (1977), 729--743.

\bibitem{Yok98} K.~Yokoyama, Mourre theory for time-periodic systems, \emph{Nagoya Math.\ J.}\ 149 (1998), 193--210.

\bibitem{Yos90} H.~Yoshida, Construction of higher order symplectic integrators, \emph{Phys.\ Lett.\ A} 150 (1990), 262--268.

\end{thebibliography}
}


{\small
\begin{thebibliography}{99}\addcontentsline{toc}{section}{References for the supplementary material}
\bibitem{S:ABG96} W.~O.~Amrein, A.~Boutet de Monvel, V.~Georgescu, \emph{$C_0$-groups, commutator methods and spectral theory of $N$-body Hamiltonians}, Progr.\ Math.\ 135, Birkh\"auser, Basel, 1996.

\bibitem{S:ABCF06} M.~A.~Astaburuaga, O.~Bourget, V.~H.~Cort\'es, C.~Fern\'andez, Floquet operators without singular continuous spectrum, \emph{J.\ Funct.\ Anal.}\ 238 (2006), 489--517.

\bibitem{S:BCF23} W.~Bao, Y.~Cai, Y.~Feng, Improved uniform error bounds of the time-splitting methods for the long-time (nonlinear) Schr\"odinger equation, \emph{Math.\ Comp.}\ 92 (2023), 1109--1139.

\bibitem{S:BGSvL26} S.~Becker, N.~Galke, R.~Salzmann, L.~van Luijk, Convergence rates for the Trotter splitting for unbounded operators, \emph{Found.\ Comput.\ Math.}\ 26 (2026), 2471--2524.

\bibitem{S:BC05} S.~Blanes, F.~Casas, On the necessity of negative coefficients for operator splitting schemes of order higher than two, \emph{Appl.\ Numer.\ Math.}\ 54 (2005), 23--37.

\bibitem{S:BCR99} S.~Blanes, F.~Casas, J.~Ros, Symplectic integration with processing: a general study, \emph{SIAM J.\ Sci.\ Comput.}\ 21 (1999), 711--727.

\bibitem{S:Boa54} R.~P.~Boas, \emph{Entire Functions}, Academic Press, New York, 1954.

\bibitem{S:BFHJY24} D.~Burgarth, P.~Facchi, A.~Hahn, M.~Johnsson, K.~Yuasa, Strong error bounds for Trotter and Strang-splittings and their implications for quantum chemistry, \emph{Phys.\ Rev.\ Research} 6 (2024), 043155.

\bibitem{S:But69} J.~C.~Butcher, The effective order of Runge--Kutta methods, in \emph{Conference on the Numerical Solution of Differential Equations}, Lecture Notes in Math.\ 109, Springer, Berlin, 1969, 133--139.

\bibitem{S:Car26} R.~Carles, Time splitting and error estimates for nonlinear Schr\"odinger equations with a potential, \emph{Found.\ Comput.\ Math.}\ 26 (2026), 2207--2258.

\bibitem{S:CM26} R.~Carles, G.~Maierhofer, On scattering for NLS: rigidity properties and numerical simulations via the lens transform, \emph{Math.\ Models Methods Appl.\ Sci.}\ (2026), doi:10.1142/S0218202526500302.

\bibitem{S:CS24} R.~Carles, C.~Su, Scattering and uniform in time error estimates for splitting method in NLS, \emph{Found.\ Comput.\ Math.}\ 24 (2024), 683--722.

\bibitem{S:Cor95} H.~O.~Cordes, \emph{The Technique of Pseudodifferential Operators}, London Math.\ Soc.\ Lecture Note Ser.\ 202, Cambridge University Press, Cambridge, 1995.

\bibitem{S:CCLR01} O.~Costin, R.~D.~Costin, J.~L.~Lebowitz, A.~Rokhlenko, Evolution of a model quantum system under time periodic forcing: conditions for complete ionization, \emph{Comm.\ Math.\ Phys.}\ 221 (2001), 1--26.

\bibitem{S:CT73} J.-M.~Combes, L.~Thomas, Asymptotic behaviour of eigenfunctions for multiparticle Schr\"odinger operators, \emph{Comm.\ Math.\ Phys.}\ 34 (1973), 251--270.

\bibitem{S:DG97} J.~Derezi\'nski, C.~G\'erard, \emph{Scattering Theory of Classical and Quantum $N$-Particle Systems}, Texts and Monographs in Physics, Springer, Berlin, 1997.

\bibitem{S:DSS11} G.~Derks, S.~Maad Sasane, B.~Sandstede, Perturbations of embedded eigenvalues for the planar bilaplacian, \emph{J.\ Funct.\ Anal.}\ 260 (2011), 340--398.

\bibitem{S:DS99} M.~Dimassi, J.~Sj\"ostrand, \emph{Spectral Asymptotics in the Semi-Classical Limit}, London Math.\ Soc.\ Lecture Note Ser.\ 268, Cambridge University Press, Cambridge, 1999.

\bibitem{S:Dol64} J.~D.~Dollard, Asymptotic convergence and the Coulomb interaction, \emph{J.\ Math.\ Phys.}\ 5 (1964), 729--738.

\bibitem{S:DF07} G.~Dujardin, E.~Faou, Normal form and long time analysis of splitting schemes for the linear Schr\"odinger equation with small potential, \emph{Numer.\ Math.}\ 108 (2007), 223--262.

\bibitem{S:Fao12} E.~Faou, \emph{Geometric Numerical Integration and Schr\"odinger Equations}, Zurich Lectures in Advanced Mathematics, European Mathematical Society, Z\"urich, 2012.

\bibitem{S:Fri48} K.~O.~Friedrichs, On the perturbation of continuous spectra, \emph{Comm.\ Pure Appl.\ Math.}\ 1 (1948), 361--406.

\bibitem{S:FRT13} C.~Fern\'andez, S.~Richard, R.~Tiedra de Aldecoa, Commutator methods for unitary operators, \emph{J.\ Spectr.\ Theory} 3 (2013), 271--292.

\bibitem{S:HLW06} E.~Hairer, C.~Lubich, G.~Wanner, \emph{Geometric Numerical Integration}, 2nd ed., Springer Ser.\ Comput.\ Math.\ 31, Springer, Berlin, 2006.

\bibitem{S:HHZ19} M.~Heyl, P.~Hauke, P.~Zoller, Quantum localization bounds Trotter errors in digital quantum simulation, \emph{Sci.\ Adv.}\ 5 (2019), eaau8342.

\bibitem{S:Hor76} L.~H\"ormander, The existence of wave operators in scattering theory, \emph{Math.\ Z.}\ 146 (1976), 69--91.

\bibitem{S:Hor83} L.~H\"ormander, \emph{The Analysis of Linear Partial Differential Operators II}, Grundlehren Math.\ Wiss.\ 257, Springer, Berlin, 1983.

\bibitem{S:How79} J.~S.~Howland, Scattering theory for Hamiltonians periodic in time, \emph{Indiana Univ.\ Math.\ J.}\ 28 (1979), 471--494.

\bibitem{S:Ign11} L.~I.~Ignat, A splitting method for the nonlinear Schr\"odinger equation, \emph{J.\ Differential Equations} 250 (2011), 3022--3046.

\bibitem{S:JL00} T.~Jahnke, C.~Lubich, Error bounds for exponential operator splittings, \emph{BIT} 40 (2000), 735--744.

\bibitem{S:JMP84} A.~Jensen, E.~Mourre, P.~Perry, Multiple commutator estimates and resolvent smoothness in quantum scattering theory, \emph{Ann.\ Inst.\ H.\ Poincar\'e Phys.\ Th\'eor.}\ 41 (1984), 207--225.

\bibitem{S:Kat04} Y.~Katznelson, \emph{An Introduction to Harmonic Analysis}, 3rd ed., Cambridge University Press, Cambridge, 2004.

\bibitem{S:Lev96} B.~Ya.~Levin, \emph{Lectures on Entire Functions}, Transl.\ Math.\ Monogr.\ 150, AMS, Providence, 1996.

\bibitem{S:Lub08} C.~Lubich, On splitting methods for Schr\"odinger--Poisson and cubic nonlinear Schr\"odinger equations, \emph{Math.\ Comp.}\ 77 (2008), 2141--2153.

\bibitem{S:MQ02} R.~I.~McLachlan, G.~R.~W.~Quispel, Splitting methods, \emph{Acta Numer.}\ 11 (2002), 341--434.

\bibitem{S:RS78} M.~Reed, B.~Simon, \emph{Methods of Modern Mathematical Physics IV: Analysis of Operators}, Academic Press, New York, 1978.

\bibitem{S:RS79} M.~Reed, B.~Simon, \emph{Methods of Modern Mathematical Physics III: Scattering Theory}, Academic Press, New York, 1979.

\bibitem{S:RS80} M.~Reed, B.~Simon, \emph{Methods of Modern Mathematical Physics I: Functional Analysis}, revised ed., Academic Press, New York, 1980.

\bibitem{S:RST19} S.~Richard, A.~Suzuki, R.~Tiedra de Aldecoa, Quantum walks with an anisotropic coin II: scattering theory, \emph{Lett.\ Math.\ Phys.}\ 109 (2019), 61--88.

\bibitem{S:Sim05} B.~Simon, \emph{Trace Ideals and Their Applications}, 2nd ed., AMS, Providence, 2005.

\bibitem{S:SW98} A.~Soffer, M.~I.~Weinstein, Nonautonomous Hamiltonians, \emph{J.\ Stat.\ Phys.}\ 93 (1998), 359--391.

\bibitem{S:SW99} A.~Soffer, M.~I.~Weinstein, Resonances, radiation damping and instability in Hamiltonian nonlinear wave equations, \emph{Invent.\ Math.}\ 136 (1999), 9--74.

\bibitem{S:Tie20} R.~Tiedra de Aldecoa, Stationary scattering theory for unitary operators with an application to quantum walks, \emph{J.\ Funct.\ Anal.}\ 279 (2020), 108704.

\bibitem{S:Yaf92} D.~R.~Yafaev, \emph{Mathematical Scattering Theory: General Theory}, Transl.\ Math.\ Monogr.\ 105, AMS, Providence, 1992.

\bibitem{S:Yaj77} K.~Yajima, Scattering theory for Schr\"odinger equations with potentials periodic in time, \emph{J.\ Math.\ Soc.\ Japan} 29 (1977), 729--743.

\bibitem{S:Yok98} K.~Yokoyama, Mourre theory for time-periodic systems, \emph{Nagoya Math.\ J.}\ 149 (1998), 193--210.

\bibitem{S:Yos90} H.~Yoshida, Construction of higher order symplectic integrators, \emph{Phys.\ Lett.\ A} 150 (1990), 262--268.

\end{thebibliography}
}
\end{document}